\documentclass[11pt,amsfonts]{article}
\usepackage{graphicx}
\usepackage{latexsym}
\usepackage{amsfonts}
\usepackage{amssymb}
\usepackage{amsmath}
\usepackage{amsfonts,amsmath,latexsym}
\usepackage[T1]{fontenc}
\usepackage{enumerate}
\usepackage{epsf}
\usepackage{amsthm}
\usepackage{color}
\usepackage{amssymb}
\usepackage{fancyhdr}

\newtheorem{prop}{Proposition}[section]
\newtheorem{rem}[prop]{Remark}
\newtheorem{df}{Definition}[section]  
\newtheorem{thm}[prop]{Theorem}
\newtheorem{lem}[prop]{Lemma}
\newtheorem{cor}[prop]{Corollary}
\newtheorem{ex}[prop]{Example}
\newtheorem{hyp}[prop]{Hypothesis}
\def\real{{\mathord{{\rm I\kern-2.8pt R}}}}        % Fake blackboard bold R.
\def\inte{{\mathord{{\rm I\kern-2.8pt N}}}}

\def\sZZ{{\rm Z\kern-2.8ptem{}Z}}

\def\z{{\mathchoice
  {\sZZ}
  {\sZZ}
  {\rm Z\kern-0.30em{}Z}
  {\rm Z\kern-0.25em{}Z} }}
\def\sQQ{{\kern 0.27em \vrule height1.45ex width0.03em depth0em
          \kern-0.30em \rm Q}}
\def\qu{{\mathchoice
    {\sQQ}
    {\sQQ}
  {\kern 0.225em \vrule height1.05ex width0.025em depth0em \kern-0.25em \rm Q}
  {\kern 0.180em \vrule height0.78ex width0.020em depth0em \kern-0.20em \rm Q}
        }}
\def\sCC{{\kern 0.27em \vrule height1.45ex width0.03em depth0em
          \kern-0.30em \rm C}}
\def\complex{{\mathchoice
    {\sCC}
    {\sCC}
  {\kern 0.225em \vrule height1.05ex width0.025em depth0em \kern-0.25em \rm C}
  {\kern 0.180em \vrule height0.78ex width0.020em depth0em \kern-0.20em \rm C}
        }}

\newcommand{\R}{\mathbb{R}}
\newcommand{\N}{\mathbb{N}}

\newcommand{\ba}{\begin{array}}
\newcommand{\ea}{\end{array}}
\newcommand{\be}{\begin{equation}}
\newcommand{\ee}{\end{equation}}
\newcommand{\bea}{\begin{eqnarray}}
\newcommand{\eea}{\end{eqnarray}}
\newcommand{\beaa}{\begin{eqnarray*}}
\newcommand{\eeaa}{\end{eqnarray*}}

\newcommand{\eps}{\varepsilon}

\newcommand{\Xt}{(X_t)_{t \geq 0}}
\newcommand{\LOm}{L^2(\Omega)}
\newcommand{\mb}{\bar{\mu}}
\newcommand{\pR}{\frac{\pa ^2 R}{\pa s_1 \pa s_2}}
\newcommand{\RR}{\mathbb{R}^2_+}

\newcommand{\iin}{\int_0^\infty}
\newcommand{\nn}{\nonumber}
\newcommand{\fr}{\frac{1}{2}}
\newcommand{\he}{\chi_{\varepsilon}}
\newcommand{\pRe}{\frac{\pa ^2 R_{\varepsilon}}{\pa s_1 \pa s_2}}
\newcommand{\Lr}{\tilde{L}_R}

\newcommand{\hn}{\chi_{[0,n]}}
\newcommand{\hm}{\chi_{[0,m]}}
\newcommand{\Co}{C^{0,0}(\overline{\math R}_+)}
\newcommand{\Fc}{\mathcal{F}C_b^\infty}
\newcommand{\bL}{\bar{L}_R}
\newcommand{\DD}{|\mathbb D^{1,2}|}
\newcommand{\Dd}{\mathbb D^{1,2}}
\newcommand{\Dstar}{Dom(\delta)^*}
\newcommand{\jeden}{1_{[0,t]}}

\newcommand{\DDd}{|\mathbb D^{1,2}(L_R)|}
\newcommand{\Ddd}{\mathbb D^{1,2}(L_R)}
\newcommand{\cE}{\mathcal{E}}
\def\a{\alpha}

\def\g{\gamma}
\def\d{\delta}

\def\z{\zeta}
\def\k{\kappa}
\def\l{\lambda}
\def\m{\mu}
\def\n{\nu}
\def\s{\sigma}
\def\t{\tau}
\def\f{\varphi}

\def\o{\omega}
\def\r{\rho}

\font\tenmath=msbm10 \font\sevenmath=msbm7 \font\fivemath=msbm5
\newfam\mathfam \textfont\mathfam=\tenmath
\scriptfont\mathfam=\sevenmath \scriptscriptfont\mathfam=\fivemath
\def\math{\fam\mathfam}

\def \={{\buildrel {\rm (law)} \over =}}
\def \N{{\math N}}
\def \R{{\math R}}

\def\Om{\Omega}
\def\cE{{\cal E}}

\def\cH{{\cal H}}

\def\pa{\partial}

\def\qed{ \hfill \vrule width.25cm height.25cm depth0cm\smallskip}

\newcommand{\basa}{\begin{assumption}}
\newcommand{\easa}{\end{assumption}}

\newcommand{\bas}{\begin{assum}}
\newcommand{\eas}{\end{assum}}

\def\lra{\longrightarrow}

\def\pa{\partial}

\newcommand{\ignore}[1]{}
\begin{document}

\renewcommand{\thefootnote}{\fnsymbol{footnote}}

\title{Malliavin-Skorohod calculus and Paley-Wiener integral
for covariance singular processes}

\author{{\sc Ida KRUK 
$^*$}\footnote{Universit\'e Paris 13,
    Math\'ematiques LAGA, Institut Galil\'ee, 99 Av. J.B. Cl\'ement 93430
    Villetaneuse. E-mail:{\tt ida.kruk@googlemail.com
     }}\thanks{Banque Cantonale de Gen\`eve,  
Case Postale 2251, 
 CH-1211 Gen\`eve 2}
\ {\sc and}\ {\sc Francesco RUSSO }  
\footnote{ENSTA ParisTech, Unit\'e de Math\'ematiques 
appliqu\'ees, F-75739 Paris Cedex  15 (France)}
\thanks{INRIA Rocquencourt 
and Cermics Ecole des Ponts, Projet MathFi.
 E-mail: \quad {\tt  francesco.russo@ensta-paristech.fr}}
}

\date{November 29th 2010}
\maketitle
\thispagestyle{myheadings}
\markright{}

\begin{abstract}
We develop a stochastic analysis for a Gaussian process $X$
with singular covariance by an intrinsic procedure
focusing on several examples such as covariance measure 
structure processes, bifractional Brownian motion,
processes with stationary increments.
We introduce some new spaces associated with the 
self-reproducing kernel space and we define the Paley-Wiener
integral of first and second order even when $X$ is
only a square integrable process continuous in $L^2$. If $X$ has stationary 
increments, we provide necessary and sufficient conditions so that 
its paths belong to the self-reproducing kernel space.
We develop Skorohod calculus and its relation with
symmetric-Stratonovich type integrals and two types of It\^o's formula.
One of Skorohod type, which works under very general (even very singular)
 conditions
for the covariance; the second one of symmetric-Stratonovich type,
which works, when the covariance is at least as regular as the one of a
 fractional Brownian motion of Hurst index equal to $H = \frac{1}{4}$.
\end{abstract}

\vskip0.5cm



\vfill  \eject

\section{Introduction}

\setcounter{equation}{0}

Classical stochastic calculus is based on It\^o's integral. It operates
when the integrator $X$ is 
 a semimartingale. The present  paper concerns some specific
aspects of calculus for Gaussian non-semimartingales,
with some considerations about Paley-Wiener integrals for non Gaussian processes.
Physical modeling, hydrology,  telecommunications,
 economics and  finance  has generated the  necessity to 
make stochastic calculus with respect to more general processes
than semimartingales.
% in this class there are Gaussian and non-Gaussian
%processes. 
In the family of Gaussian processes, the most 
 adopted and    celebrated example is of course Brownian motion.
In this paper we will consider Brownian motion $B$ 
(and the processes having the same type of path regularity)
as the frontier of two large classes of processes: those whose path regularity 
is more regular than  $B$ and those whose path regularity is 
more singular than $B$. 
It is well-known that $B$ is a finite quadratic variation process
in the sense of \cite{Follmer} or \cite{RV1995, RV2000}. 
Its quadratic variation process is given by $[B]_t = t.$
If we have a superficial look to 
 path regularity, we can
{\it macroscopically}  distinguish three classes.
\begin{enumerate}
\item Processes $X$ which are more regular than Brownian motion $B$
are such that $[X] = 0$.
\item Processes $X$ which are as regular as $B$ are finite non-zero
quadratic variation processes.
\item  Processes $X$ which are more singular than $B$ are processes
which do not admit a  quadratic variation.
\end{enumerate}
Examples of processes $X$ summarizing previous classification
are given by the so called fractional Brownian motions $B^H$ with 
Hurst index $ 0 < H < 1$. If $H > \frac{1}{2}$
(resp. $H = \frac{1}{2}, \quad   H < \frac{1}{2}$),
 $B^H$ belongs to class 1. (resp. 2., 3.). \\
Calculus with respect to integrators which are not semimartingales
is more than thirty-five  years old, but a real activity acceleration 
was produced since  the mid-eighties. A significant starting lecture
note in this framework is \cite{KU}: well-written lecture
articles from H. Kuo, D. Nualart, D. Ocone constitute excellent 
pedagogical articles on Wiener analysis  (Malliavin calculus,
white noise calculus) and its application to anticipative calculus
via Skorohod integrals.

Since then, a  huge amount of papers have been produced, and it
 is impossible to list them here, so we will essentially
quote monographs.
% however
%we are still not so close from having a truly efficient
% approach for applications.

There are two mainly techniques for studying non-semimartingales integrators.
\begin{itemize}
\item Wiener analysis, as we mentioned before.
It is based on the so-called Skorohod integral
(or divergence operator); it allows 
to ``integrate'' anticipative integrands with respect to
various Gaussian processes. 
We quote for instance \cite{Nualart}, \cite{BHOZ}, \cite{ustunel, Mall},
for Malliavin calculus and \cite{HKPS, HOUZ} for white noise calculus.
In principle, Malliavin calculus   can be abstractly implemented
on any Wiener space generated by any Gaussian process $X$,
see for instance \cite{Watanabe}, but in the general abstract 
framework, integrands may live in some abstract spaces.
In most of the present literature about Malliavin-Skorohod
calculus,   integrands are supposed to live in 
a space $L$ which is isomorphic to 
the self-reproducing kernel space.
One condition for an integrand $Y$ to belong to the classical domain of the
divergence operator is that its paths belong a.s. to $L$. 
More recent papers as \cite{CN, MV} allow integrands to live outside
the classical domain of the divergence operator. 
Some activity about Skorohod integration was also performed  in the framework of Poisson measures integrators,
see e.g. \cite{nunno}.
%white noise  calculus, EX BENDER, NUALART, VIENS-MOCIOALCA ) .
\item Pathwise and quasi-pathwise  related techniques,
as rough paths techniques \cite{QL}, regularization 
(see for instance \cite{RV1993, RVsem})
or discretization techniques \cite{Follmer}, but also fractional
integrals techniques \cite{Z1, Z2} and also \cite{HZ}
for connections between rough paths and
fractional calculus.
%\item Dirichlet forms.
\end{itemize}
%Dirichlet forms techniques concern general Markov process
%integrator and it will not concern this paper
%which is situated in the first category with some
%connection to the second one.
This paper is the continuation of \cite{KRT}, which focused on the processes belonging to 
categories 1. and 2. Here we are mainly interested
in processes of category 3.
We develop some intrinsic stochastic analysis with respect to
those processes. The qualification {\it intrinsic} is related  
to the fact that in opposition to \cite{ALN, AMN, CN, MV}, where
there is no underlying Wiener process. 
We formulate a class of general assumptions 
Assumptions (A), (B), (C($\nu$)), (D) under which the calculus runs.
Many properties hold only under the two three first hypotheses.
Sometimes, however, we also make use of a supplementary assumption that
we believe to be technical, i.e. 
\begin{equation} \label{ConstAfterT}
X_t = X_T, \quad t > T.
\end{equation}
In some more specific situations, we also introduce Hypothesis
\eqref{ELTR} which intervenes for instance to guarantee, that $X$
itself belongs to the {\it classical} domain of the divergence operator. 
At Section \ref{s5.1} we also define a
 suitable Hilbert space  $L_R$
for integrand processes,
which is related to self-reproducing kernel space.
We describe the content 
of $L_R$ in many situations.
We implement the analysis and we verify the assumptions
in the following examples: the case when 
the process $X$ is defined through a kernel integration
with respect to a Wiener process, the case of processes 
with a covariance measure structure, the processes with stationary
increments, the case of bifractional Brownian motion. 
We provide a stochastic analysis framework starting from
Paley-Wiener integral for second order processes.
The Wiener integral with respect to the subclass of processes with 
stationary increments was studied with different techniques in \cite{Jolis}.
In our paper, we also define the notion to
a multiple Paley-Wiener integral,
involving independent processes $X^1, \cdots, X^n$.
If $n=2$ those integrals  have a natural relation with
 the notion of L\'evy's area in rough path theory.

Starting from Section \ref{s6}, we concentrate on Malliavin-Skorohod
 calculus,
with It\^o's type formulae and connections with
symmetric-Stratonovich integrals via regularization.
%%% COMPARE WITH CLASSICLA SOBOLEV SPACES
%makes some explorations between Wiener, Skorohod type integrals
%and integrals via regularizations when the integrator is singular,
%in particular it belongs to the category 3. above.
%APRES PRECISER EXACTEMENT CE QUE L'ON FAIT \\
%PROPRIETES TRAJECTORIELLES ETC\\
%An analysis in the same spirit concerned, but for ``regular
%integrators'' of the category 2., was carried on in \cite{KRT}.
Calculus  via regularizations was started by
 F.  Russo and P.  Vallois \cite{RVCRAS} 
 developing  a regularization procedure, 
whose philosophy is similar to the discretization.
 They introduced a forward (generalizing It\^o) integral,
and the symmetric (generalizing Stratonovich) integral.
%We recall the process $X$ it self is forward integrable
%(in symbols $\int_0^\cdot X d^-X$ exists 
%if and only if  $X$ has a finite quadratic variation,
%see for instance \cite{GRV}.
%Since processes of category 3. are not finite quadratic
%variation processes, than Skorohod type integrals
%will be linked with the symmetric integral;
%in \cite{KRT} the connection was related with forward type 
%integrals.

As we said, in the first six sections, we redefine 
a Paley-Wiener type integral 
 with respect to an $L^2$-continuous square integrable  process.
 We aim at showing some  interesting 
features and difficulties, which are encountered if one wants to
 define the integral in a natural function space avoiding distributions. 

 This allows in particular, but not only, 
to settle  the basis of Malliavin-Skorohod calculus 
for Gaussian processes with singular covariance.

 %An approach to stochastic calculus for "singular" processes was 
%performed by \cite{Latorre}. They had in mind to provide some framework,
% which could recover some techniques for Stratonovich type It\^o
% formula, as the one obtained for instance by \cite{RT} when $X$
% is a fractional Brownian motion with parameters $H$, $K$
% with $HK>\frac{1}{6}$.
% The main objective was however not to 
%obtain a Skorohod type calculus but more to use some Malliavin calculus 
%ideas to recover pathwise type techniques. \\
As we said, Malliavin calculus, according to \cite{Watanabe},
can be  developed abstractly for any Gaussian process
$X=(X_{t})_{t\in [0,T]}$. The Malliavin derivation can be naturally
defined on a general Gaussian abstract Wiener space. A Skorohod  integral (or divergence)
can also be  defined as the adjoint of the Malliavin derivative. 

The
 crucial ingredient is the canonical Hilbert space ${\cal{H}}$
 (called also, improperly,  by some authors  reproducing kernel Hilbert
 space) of the Gaussian process $X$ which is  defined as the closure
 of the linear space generated by the indicator functions $\{
 1_{[0,t]}, t\in [0,T]\}$ with respect to the scalar product
 \begin{equation}
 \label{ps} \langle 1_{[0,t]}, 1_{[0,s]}\rangle _{{\cal{H}}}= R(t,s),
 \end{equation}
 where $R$ denotes the covariance of $X$.
Nevertheless, this
calculus remains more or less abstract if the structure of the
elements of the Hilbert  space  is not known. When we
say abstract, we refer to the fact that, for example, it is
difficult to characterize the processes which are integrable with
respect to $X$, or to establish It\^o formulae.
% to estimate the $L^{p}$-norms of the Skorohod
%integrals or to push further this calculus to obtain an
%It\^o typ

In this paper, as we have anticipated, 
 we formulate some natural assumptions 
(A), (B), (C($\nu$)), (D),
 that the underlying process has to fulfill,
which let us efficiently define a Skorohod intrinsic calculus and
It\^o formulae, when integrators belong to categories 2. and 3. In particular, Assumption (D) truly translates
the singular character of the covariance.  
%A significant paper is \cite{MV} which aims at treating
%Malliavin calculus with singular covariance of the
%type  $X_t = \int_0^t g(t-s) dW_s$ where $g$
%is a locally square integrable function and $W$ is a
%classical Brownian motion and the underlying Wiener
%space is reated to $W$.

We link Skorohod integral with 
 integrals via regularization (so of almost pathwise type) 
 similarly to \cite{KRT}, where the connection was established with forward integrals.
We recall that the process $X$ is forward integrable
(in symbols $\int_0^\cdot X d^-X$ exists) 
if and only if  $X$ has a finite quadratic variation,
see for instance \cite{GRV} or \cite{GNRV}. 
Therefore if $X$ is a fractional
Brownian motion with Hurst index $H$,
the forward integral $\int_0^\cdot X d^-X$ exists
if and only if
$H \geq \frac{1}{2}$; on the other hand
the symmetric integral
$\int_0^\cdot X d^\circ X$ always exists.
% if and only if $H \ge \frac{1}{4}$.
Since we are mainly interested in singular covariance
processes (category 3.), which are not of finite
quadratic variation,
 Skorohod type integrals
will be linked with the symmetric integrals.

%Since we aim at treating integrators with singular 
%covariance, the connection can only be reasonably done with symmetric
%type integrals. In fact, see e.g. \cite{GNRV}, if $X$ is a fractional
%Brownian motion with Hurst index $H$
%the forward integral $\int_0^\cdot X d^-X$ exists if and only if
%$H < \frac{1}{2}$; on the other hand
%the symmetric integral exists if and only if $H \ge \frac{1}{4}$.

As we have mentioned before, a particular case was deeply analyzed
  in the literature.
 We refer here to
the situation when the covariance $R$ can be explicitly written  as
\begin{equation*}
R(t,s)= \int _{0}^{t\wedge s} K(t,u)K(s,u)du,
\end{equation*}
where $K(t,s)$, $0<s<t<T$, is a deterministic kernel satisfying some
regularity conditions. Enlarging, if needed, our probability space, we
can express the process $X$ as
\begin{equation}
\label{rep} X_{t}=\int _{0}^{t}K(t,s)dW_{s},
\end{equation}
where $(W_{t})_{t\in [0,T]}$ is a standard Wiener process and the
above integral is understood in the Wiener sense. In this case, more
concrete results can be proved, see \cite{AMN, De, MV}.
In this framework the underlying Wiener process $(W_t)$ is strongly used for 
developing anticipating calculus. 

For illustration, we come back to the case, when $X$ is a fractional
Brownian motion $B^H$ and $H$ is the Hurst index.
%Of course, the most studied case is  the fractional
 %Brownian motion (fBm), due to the multiple applications of this
 %process in various area.
The process $B^{H}$ admits the Wiener integral representation
(\ref{rep}) and the kernel $K$ together with the  space $\cal H$ can be
characterized by the mean of fractional integrals and derivatives,
see \cite{AMN, AN, DU, PT, CN, bender} among others. As a
consequence, one can prove for any $H > \frac{1}{4}$ (to guarantee
that $B^H$ is in the domain of the divergence), the following It\^o's formula:
\begin{equation*}
f(B^{H}_{t}) = f(0) + \int _{0}^{t} f'(B^{H}_{s})\delta B^{H}_{s} + H\int
_{0}^{t} f''(B^{H}_{s}) s^{2H-1}ds.
\end{equation*}

\cite{MV} 
puts emphasis on the case $K(t,s) = g (t-s)$, when the variance scale 
of the process is as general as possible, including logarithmic scales.
%That paper is very significant and it truly concentrates
%on the ``singular case'': in terms of fractional Brownian motion
%intuition, this refers to the case $H < \frac{1}{4}$.

%MARIA JOLIS. TROUVER PREPRINT CONCERNANT
%PROCESSUS DE CARRE' INTEGRABLE 

% The canonical space ${\cal{H}}$ can be written as
% \begin{equation*}
% {\cal{H}}= \left( K^{\ast }\right) ^{-1}\left( L^{2}([0,T])\right)
% \end{equation*}
% where the "transfer operator" $K^{\ast }$ is defined on the set of
% elementary functions as
% \begin{equation*}
% K^{\ast}(\varphi) (s) = K(T,s) \varphi (s)+ \int _{s}^{T} \left(
% \varphi (r)-\varphi (s) \right) K(dr,s)
% \end{equation*}
% and extended (if possible) to
% ${\cal{H}}$ (or a set of functions contained in ${\cal{H}}$).
% , a stochastic process $u$ will be Skorohod integrable
% with respect to $X$ if and only if $K^{\ast } u$ is Skorohod
% integrable with respect to $W$ and $\int u\delta X= \int
% (K^{\ast}u)\delta W$. Depending on the regularity of $K$ (in
% principal the H\"older continuity of $K$ and $\frac{\partial K
% }{\partial t}(t,s)$ are of interest) it becomes possible to have
% concrete results.

In section \ref{s4.5}, we establish some connections between the "kernel approach" discussed in the literature and the "covariance intrinsic approach" studied here.

%One can also study the relation between ''pathwise type'' integrals and the
%divergence integral, the regularity of the Skorohod integral process
%or the It\^o formula for indefinite integrals.
\smallskip

As we mentioned, if the deterministic kernel $K$ in the
representation (\ref{rep}) is not explicitly known, then the
Malliavin calculus with respect to the Gaussian process $X$ remains
in an abstract form and there are of course many situations when
this kernel is not explicitly known. As a typical example, we have in
mind the case of the  {\em bifractional Brownian motion (BFBM) }
$B^{H,K}$, where $H \in ]0,1[, K \in ]0,1]$; a kernel representation is known in 
a particular case, but with respect to a space-time white noise: in fact 
% in the particular case $H = \frac{1}{2} = K$ then
 the solution $F(t) = u(t,x)$ of a classical stochastic heat equation
driven by a white noise with zero initial condition,  is  distributed as $B^{H,K}$, 
for any fixed $x$, see \cite{swanson}. Another interesting representation is
provided by \cite{LeiNualart}, which shows the existence of real constants
$c_1, c_2$ and an absolutely continuous process $X(H,K)$
independent of $B^{H,K}$, such that 
$c_1 B^{H,K} + c_2 X(H,K)$ is distributed as a fractional Brownian motion
with parameter $HK$. 
In spite of those considerations, finding a kernel $K$,
such that $B^{H,K}_t = \int_0^t K(t,s) dW_s$,
is still an open problem.
%This process, denoted by $B^{H,K}$, is defined as a centered
%Gaussian process starting from zero with covariance
%\begin{equation}
%\label{covbiFBM} R(t,s)= \frac{1}{2^{K}}\left( \left( t^{2H} +
%s^{2H}\right) ^{K} -\vert t-s\vert ^{2HK}\right)
%\end{equation}
%where $H\in (0,1) $ and $K\in (0,1]$. When $K=1$, then we have a
%standard fractional Brownian motion.
Bifractional Brownian motion
 was introduced in \cite{HV} and a {\it quasi-pathwise type}
of regularization (\cite{RVsem}) type  approach to stochastic calculus
 was provided in \cite{RT}.
It is possible for instance to obtain
 an It\^o formula of the Stratonovich type (see \cite{RT}), i.e.
\begin{equation} \label{ERT}
f(B^{H,K}_{t}) = f(0)+ \int_{0}^{t} f'(B^{H,K}_{s}) d^{\circ}
B^{H,K}_{s}
\end{equation}
for any parameters $H\in (0,1)$ and $K\in (0,1]$
such that $HK > \frac{1}{6}$.
 An interesting property of
$B^{H,K}$ consists in the expression of its quadratic variation, 
defined as usual, as a limit of Riemann sums, or
in the sense of regularization.
 The following properties hold true.
\begin{description}
\item{$\bullet$ } If $2HK>1$, then the quadratic variation of
$B^{H,K}$ is zero and $B^{H,K}$ belongs to category 1.
\item{$\bullet$ } If $2HK<1$ then the quadratic variation of
$B^{H,K}$ does not exist and $B^{H,K}$ belongs to
category 3.
\item{$\bullet$ } If $2HK=1$ then the quadratic variation of
$B^{H,K}$ at time $t$ is equal to $2^{1-K}t$
 and $B^{H,K}$ belongs to
category 2.
\end{description}
The last property is remarkable; indeed, for $HK=\frac{1}{2}$  we
have a Gaussian process which has the same quadratic variation as
the Brownian motion. Moreover, the processes is not a semimartingale
(except for the case $K=1 $ and $H=\frac{1}{2}$), it is
self-similar, has no stationary increments.
%It is possible for instance to obtain
% an It\^o formula of the Stratonovich type (see \cite{RT})
%\begin{equation*}
%f(B^{H,K}_{t}) = f(0)+ \int_{0}^{t} f'(B^{H,K}_{s}) d^{\circ}
%B^{H,K}_{s}
%\end{equation*}
%for any parameters $H\in (0,1)$ and $K\in (0,1]$. \\
%The purpose of this work is to continue the work started in
%\cite{KRT} where a first step to expressing a Malliavin-Skorohod
%calculus for Gaussian process only in terms of the covariance.

Motivated by the consideration above, one developed in \cite{KRT} a Malliavin-Skorohod  calculus with
respect to Gaussian processes $X$ having a {\em covariance  measure structure} 
in  sense that  the covariance  is the distribution function of
 a (possibly signed) measure $\mu_R$ on ${\cal{B}}([0,T]^{2})$. We denote by $D_t$ the diagonal set
 \[
 \left\{(s,s)| s \in [0,t] \right\}
 \]
%This approach is particularly suitable for processes 
%whose  representation form  (\ref{rep}) is not explicitly given.
The processes having a {\em covariance  measure structure} belong
to the category 1. (resp. 2.), i.e. they are 
more regular than  Brownian motion (resp. as regular as Brownian motion)
if $\mu_R$  restricted to the diagonal $D_T$
vanishes (resp. does not vanish). In particular, it was shown that in this case $X$ is a finite quadratic variation process and $[X]_t = \m(D_t)$. 
This paper continues the spirit of \cite{KRT}, but it concentrates
on the case when $X$ is  less regular or equal than Brownian motion.

A significant paper, is \cite{Latorre}, which 
 establishes a It\^o-Stratonovich (of quasi-pathwise type,
in the discretization spirit) for  processes belonging to class 3.,
 i.e. less regular than
Brownian motion. In particular, the paper rediscovers 
It\^o's formula of \cite{RT} for $B^{H,K}$, if 
$HK > \frac{1}{6}$; for this purpose 
the authors implement innovating Malliavin calculus
techniques. 
Their main objective was however not to 
obtain a Skorohod type calculus, but more to use some Malliavin calculus 
ideas to recover pathwise type techniques. \\

In this paper, for simplicity of notations
and without restriction of generality,  we consider 
processes indexed by the whole first quarter of the completed plane 
 ${\bar \R}^2_+ = [0,\infty]^2 $. In particular we suppose that 
$X$ is a continuous  process in $L^2$, such that
$\lim_{s \to \infty} X_s $ exists, and it is denoted $X_\infty$
and under some circumstances we suppose even \eqref{ConstAfterT}.
%Sometimes however, we will need to introduce
%the more restrictive assumption.
%\begin{equation} \label{ConstAfterT}
%X_t = X_T, \quad t > T.
%\end{equation}
Let $R(s_1,s_2), s_1, s_2 \in [0, \infty]$ be the covariance function
of $X$.
As we said, we introduce a class of natural assumptions 
which have to be fulfilled in most of the results
in order to get an efficient Skorohod calculus. \\
The  processes of class 2. and 3., so essentially less
regular than Brownian motion, will fulfill the following:
\begin{itemize}
\item $R(ds,\infty)$ is a non-negative real measure.
\item If $D$ is the first diagonal of $\R^2_+$,
the Schwartz distribution $\partial^2_{s_1, s_2} R$
 restricted to $\R^2_+ \backslash D$ is  a 
non-positive $\sigma$-finite measure.
\end{itemize}
This will constitute the convenient
Assumption (D).
%PEUT ETRE DIRE QUE DANS UN PREMIER TEMPS NOUS NE 
%FORMULONS PAS D'HYPOTHESE DE ``GAUSSIANITE''

% We will see  that under this assumption, we can define suitable spaces
%on which the construction of the Malliavin derivation/Skorohod
%integration is coherent. 
%In fact, our initial purpose is more ambitious; we start to
% construct a stochastic analysis for general (non-Gaussian) processes
% $X$ having a covariance measure $\mu$. We define Wiener integrals
% for a large enough class of deterministic functions and we define a
% Malliavin derivative and a Skorohod integral with respect to it; we
% can also prove certain relations and properties  for these
% operators. However,  if one wants 
% to produce a consistent theory,   then the Skorohod integral applied to
% deterministic integrands should coincide with the Wiener integral.
% This property is based on  
% {\em integration by parts } on  Gaussian spaces which is proved
%  in Lemma  l5.4t PROBLEME REFERENCE. As it can be seen, that proof 
%  is completely based on the Gaussian character and
% it seems difficult to prove it for general processes.
% Consequently, in the sequel, we concentrate  our study on the Gaussian case and
% we show various results as  the continuity of the integral
% processes, the chaos expansion of local times, 
%As we anticipated, this work aims at formulating a natural
%Skorohod calculus with its relations to pathwise type integrals,
%including the situation of singular covariance matrix.
The basic space of integrands for which
Paley-Wiener integral is defined is $L_R$.
This space,  under Assumption (D), plays the role of self-reproducing kernel 
space. A necessary condition for the  process $X$ itself to be in the natural
domain ($Dom \delta$)
of the divergence operator is that it belongs a.s. to $L_R$.
Following the ideas of \cite{CN, MV} one defines for our general
class of processes an extended domain called $Dom^* \delta$ which allows to proceed
when $X$ does not always belong a.s. to $L_R$.
%If $Y$ does not belong to $Dom \delta$,  

Other products of this paper are the following.
\begin{itemize}
\item We define a corresponding appropriate Paley-Wiener integral
with respect to second order processes in Section \ref{s5}. In particular, see Section \ref{r3.8}, we extend some significant considerations of \cite{PT} made in the context of fractional Brownian motion; \cite{PT} illustrates that the natural space where Wiener integral is defined, is  complete if the Hurst index is smaller or
equal to $\frac{1}{2}$. 
\item The link between symmetric and  Skorohod integrals,
see Theorem \ref{PA2}, is given by  suitable trace of Malliavin derivative of the process.
\item %One of the basic spaces of integrands for which
%Wiener integral is defined is $L_R$, defined at Section \ref{s5.1}.
%This, under Assuption (D) plays the role of self-reproducing kernel 
%space.
 If the process $X$ is  continuous, Gaussian and has stationary increments,
 we provide
necessary and sufficient conditions such that the paths of $X$
 belong to $L_R$, see Corollary \ref{crDom3}.
\item We establish an It\^o type formula for Skorohod integrals
for very singular covariation when the underlying process $X$
is quite general, continuing the work of \cite{CN} and \cite{MV}.
This is done in Proposition  \ref{p34}:
 if $f \in C^\infty$ with bounded derivatives
\begin{equation}\label{ItoSko}
f(X_t) = f(X_0) + \int_0^t f'(X_s) \delta X_s + \frac{1}{2} \int_0^t
f''(X_s) d\gamma (s), 
\end{equation}
where $\gamma(t)$ is the variance of $X_t$. \\
We recall that if $X$ is a bifractional Brownian motion with indexes $H,K$
 such that $H K = \frac{1}{2}$, then $\gamma(t) = t$ and so
equation \eqref{ItoSko} looks very similar to the one related to classical
Wiener process.
      
Formula \eqref{ItoSko} implies the corresponding formula with respect to the symmetric
integral, see Corollary \ref{TA4}, i.e.
 $$ f(X_t) = f(X_0) + \int_0^t f'(X_s) d^\circ X_s. $$.
\end{itemize}

% the ''pathwise'' and the Skorohod integrals and finally we derive the following It\^o formula,
% see Corollary 
% PROBLE REFE  cor8.11c,
% for $f \in C^2(\R)$ such that $f''$ is bounded:
% $$ f(X_t) = f(X_0) + \int_0^t f'(X_s) \delta X_s + \frac{1}{2} \int_0^t f''(X_s) d\gamma (s), $$
% where $\gamma(t) = Var (X_t)$.
%  Our main examples include the Gaussian
% semimartingales, the fBm with $H \ge  \frac{1}{2}$, the bi-fBm with
% $HK\geq \frac{1}{2}$ and processes with stationary increments.
%  In the bi-fbm case, when $2HK=1$, we find a very
% interesting fact, that is, the bi-fBm with $2HK=1 $ satisfies the
% same It\^o formula as the standard Wiener process, that is
% \begin{equation*}
% f(B^{H,K}_{t})=  f(0) + \int _{0}^{t} f'(B^{H,K}_{s}) \delta
% B^{H,K}_{s}+ \frac{1}{2}\int _{0}^{t} f''(B^{H,K}_{s})ds
% \end{equation*}
% where $\delta $ denotes the Skorohod integral.

% We will also like to mention certain aspects that could be the
% object of a further study:
% \begin{description}
% \item{$\bullet$ } the proof of the Tanaka formula involving weighted
% local times; for the fBm case, this has been proved in \cite{CNT}
% but the proofs necessitates the expression of the kernel $K$.
% \item{$\bullet$ } the two-parameter settings, as developed in e.g.
% \cite{TV}.
% \item{$\bullet$ }the proof of the Girsanov transform and the use of
% it to the study of stochastic equations driven by Gaussian noises,
% as e.g. in \cite{NO}.
% \end{description}

\smallskip

We organize our paper as follows. 
After some preliminaries
stated at Section \ref{s0}, we introduce the basic assumptions
(A), (B), (C) (or only its restricted version
(C($\nu$))),  (D) in Section \ref{s1} followed by the motivating
examples in Section \ref{s2}, including the case when
the process has a covariance measure structure treated
in \cite{KRT}.
  Section \ref{s4.5} discusses 
the link with the case that   the process is of the type
$X_t = \int_0^t K(t,s) dW_s$ for a suitable kernel $K$.
At Section \ref{s5} we define the Wiener integral for second 
order processes together with a multiorder version; in the same 
section we discuss some path properties of the underlying process
and the relation with the integrals via regularization.
Starting from  Section \ref{s6} until \ref{s8bis},  we
 introduce and discuss the basic notions of Malliavin calculus. 
At Section \ref{s9} we introduce Skorohod integrals,
at Section \ref{s10} we discuss It\^o formula 
in the very singular case. Section \ref{s4.7}
shows that Skorohod integral is truly an extension of 
Wiener integral. Finally Section \ref{s4.8} 
provides the link with integrals via regularization and It\^o's
formula with respect to symmetric integrals.

\vskip0.5cm

\section{Preliminaries} \label{s0}

\setcounter{equation}{0}

\par Let $J$ be a closed set of the type   $\R_+$,  $\R^m$ or $\R^2_+ = [0, +\infty[  \times [0, +\infty[$,  and  $k \ge 1$.
In this paper $C^{\infty}_0(J)$ (resp. $C^{\infty}_b(J), C^k_0(J), C^k_{pol}(J), C_b(J)$)
 stands for the set of functions $f: J \rightarrow \R$ which are
 infinitely differentiable with compact support (resp.  smooth with all bounded partial derivatives,
of class $C^k$ with compact support,   of class $C^k$ such that the partial derivatives of order smaller or equal to $k$
have polynomial growth, bounded functions). 
\par If $g_1,g_2: \R_+ \lra \R$, we denote $g=g_1 \otimes g_2$, the function $g: \RR \lra \R$ defined by $g(s_1,s_2) = g_1(s_1)g_2(s_2)$.
 \par Let $I$ be a subset of $\mathbb{R}^2_+$
of the form
\[
I=]a_1,b_1] \times ]a_2,b_2]
\]
Given  $g:\mathbb{R}^2_+ \rightarrow
\mathbb{R}$ we will denote
\[
\Delta_I g=g(b_1,b_2)+g(a_1,a_2)-g(a_1,b_2)-g(b_1,a_2).
\]
It constitutes the \textbf{planar increment} of $g$.
\vskip0.5cm

\begin{df} \label{ad1.1}
$g:\RR \rightarrow \mathbb{R}$ will be said to have a {\bf
bounded planar variation} if
\begin{equation}
\sup_{ \tau} \sum_{i,j=0}^{n-1}\left| \Delta_{]t_i,t_{i+1}]\times
 ]t_j,t_{j+1}]}g \right| < \infty. \nn
\end{equation} 
where $\tau  = \{0 = t_0 < \ldots < t_n<\infty \}$, $n \ge 1$,
 i.e. $\tau$ is a subdivision
of $\R_+$. Previous quantity will be denoted by $\|g \|_{pv}$.
\end{df}
If $g$ has bounded planar variation and is vanishing on the axes, then there exists a signed measure $\chi$ (difference of two positive measures) such that
\be
g(t_1,t_2)= \chi (]0,t_1]\times ]0,t_2]). \label{EPBV}
\ee
For references, see a slight adaptation of Lemma 2.2 in \cite{KRT} and Theorem 12.5 in \cite{Bill}.
\par Some elementary calculations allow to show the following.
\begin{prop} \label{pBPV} Let $g:\RR \lra \R$ of class $C^2$, $g$ has a bounded planar variation if and only if
\[
\| g\|_{pv}:= \int_{\RR} \left|\frac{\pa^2 g}{\pa t_1 \pa t_2} \right| dt_1 dt_2 < \infty.
\]
In particular, if $g$ has compact support, then $g$ has a bounded planar variation.
\end{prop}
\par 
Let $X=\Xt$, be a zero-mean continuous process in $\LOm$ such that $X_0=0$ a.s.
For technical reasons we will suppose that 
\be
\lim_{t \to \infty} X_t = X_\infty \textrm{ in } \LOm. \label{1.1}
\ee
(\ref{1.1}) is verified if for instance
\bea
 X_t &= &X_T, \ t \geq T \label{1.2a} 
\eea
for some $T > 0$.
\par We denote by $R$ the covariance function, i.e. such that:
\[
R(s_1,s_2)=Cov(X_{s_1},X_{s_2})= E(X_{s_1}X_{s_2}), \ s_1,s_2 \in \overline{\R}^2_+ = [0,\infty]^2.
\]
In particular $R$ is continuous and vanishes on the axes.
\smallskip

We convene that all the continuous functions on $\R_+$ are extended by continuity to $\R_{-}$. A continuous function $f: \RR \lra \R$ such that $f(s)=0$, if $s$ belongs to the axes, will also be extended by continuity to the whole plane.

\smallskip
In this paper $D$ will denote the diagonal $\{(t,t)| t\geq 0 \}$ of the first plane quarter $\R^2_+$.
\begin{df} \label{d0.1}
$X$ is said to have a {\bf covariance measure structure} if $\pR$ is a finite Radon measure $\mu$ on $\RR$ with compact support. We also say that $X$ has a covariance measure $\m$.
\end{df}
A priori $\frac{\pa R}{\pa s_1}, \frac{\pa R}{\pa s_2}, \frac{\pa ^2 R}{\pa s_1 \pa s_2}$ are Schwartz distributions.
In particular for $s = (s_1,s_2)$ we have
\[
R(s_1,s_2) = \mu([0,s_1] \times [0,s_2])
\]
\begin{rem} \label{r0.2}
The class of processes defined in Definition \ref{d0.1} was introduced in \cite{KRT}, where the parameter set was $[0,T]$, for some $T > 0$, instead of $\R_+$. Such processes can be easily extended by continuity to $\R_+$ setting $X_t = X_T$, if $t \geq T$. In that case the support of the measure is $[0,T]^2$.
\end{rem}
\par The present paper %and its companion \cite{KR2}
constitutes a natural continuation of \cite{KRT} trying to extend Wiener integral and Malliavin-Skorohod calculus to a large class of more singular processes. 
\par
%The present section introduces some recalls and complements related to the class of processes with covariance measure structure. 
The covariance approach is an intrinsic way of characterizing square integrable processes. These processes include Gaussian processes defined through a kernel, as for instance \cite{AMN, MV}.
\par
We will see later that a process $X_t = \int_0^t K(t,s)dW_s$, where $(W_t)_{t \geq 0}$ is a classical Wiener process and $K$ is a deterministic kernel with some regularity, provide examples of processes with covariance measure structure. Other examples were given in \cite{KRT}.
\par
One relevant object of \cite{KRT} was Wiener integral with respect to $X$. Let $\f: \R_+ \lra \R$ has locally bounded variation with compact support. We set
\[
\iin \f dX = -\iin X d\f.
\]
If $\f$ is a Borel function such that
\be \label{0.1}
\int_{\RR}|\f \otimes \f|d|\mu| < \infty,
\ee
then similarly to Section 5 of \cite{KRT}, the Wiener integral $\iin \f dX$ can be defined through the isometry property
\[
E\left(\iin \f dX\right)^2 = \int_{\RR} \f \otimes \f d\mu.
\]
\begin{rem} \label{r0.1}
If $\f$ fulfills (\ref{0.1}), then the process $Z_t = \int_0^t \f dX$ has again a covariance measure structure with measure $\nu$ defined by
\[
d\nu = \f \otimes \f d\mu.
\]
\end{rem}
\section{Basic assumptions} \label{s1}

\setcounter{equation}{0}

\par In this section we formulate a class of fundamental hypotheses, which will be in force for the present paper. %and also for \cite{KR2}.
\newline
\textbf{Assumption (A)}
\begin{description}
	\item[i)] $\forall s \in \overline{\R}_+: R(s,dx)$ is a signed measure,
	\item[ii)]$s \longmapsto \int_0^\infty |R|(s,dx)$ is a bounded function.
\end{description}
\begin{rem} \label{RA}Since $R$ is symmetric, Assumption (A) implies the following.
\begin{description}
	\item[i)'] $\forall s \in \overline{\R}_+: R(dx,s)$ is a signed measure,
	\item[ii)']$s \longmapsto \int_0^\infty |R|(dx,s)$ is a bounded function.
\end{description}
\end{rem}

\begin{rem} \label{r1.1}
Suppose that $X$ has a covariance measure $\mu$. For $s \geq 0$ we have
\[
x \longmapsto R(s,x) = \int_{[0,s] \times [0,x]} d\mu.
\]
which is a bounded variation function whose total variation is clearly given by
\[
\int_0^\infty |R|(s,dx) = \int_{[0,s] \times \R_+} d|\mu|(s_1,s_2) \leq \int_{\R_+^2} d|\mu|(s_1,s_2).
\]
Hence Assumption (A) is fulfilled.
%We observe that
%\[
%|R|(s, A) = \int_{[0,s] \times A} d|\mu|.
%\]
\end{rem}
\textbf{Assumption (B)} We suppose that
\be
\label{1.2} \mb(ds_1,ds_2) :=  \frac{\pa ^2 R}{\pa s_1 \pa s_2}(s_1,s_2)(s_1-s_2)
\ee
is a Radon measure. In this paper by a Radon measure we mean the difference of two (positive Radon) measures.
\begin{rem} \label{r1.1bis}
\begin{description}
	\item[i)] The right-hand side  of (\ref{1.2}) is well defined being the product if a $C^\infty$ function and a Schwartz distribution.
	\item[ii)] If $D$ is the diagonal introduced before in Definition \ref{d0.1}, Assumption (B) implies that $\pR$ restricted to $\R^2_+ \backslash D$ is a $\sigma$-finite measure.
	Indeed, given $\f \in C_0^\infty (\R^2_+ \backslash D)$, we symbolize by $d$ the distance between $supp \,\f$ and $D$. Setting $g(s_1,s_2) = s_1-s_2$, since $\frac{\f}{g} \in C_0^\infty (\R^2_+ \backslash D)$, we have
	\[
	\left| \left\langle \pR , \f  \right\rangle \right|= \left| \left\langle \pR \cdot g , \frac{\f}{g} \right\rangle \right| = \left| \int_{\R_+^2} d\mb \frac{\f}{g}  \right| \leq \frac{1}{\inf_{|x| \geq d}|g|(x)} \left\| \f \right\|_{\infty} \left| \mb \right| (\R^2_+).
	\]
	\item[iii)] On each compact subset of $\RR \backslash D$, the total variation measure $|\mu|$ is absolutely continuous with density $\frac {1}{|g|}$ with respect to $|\mb|$.
\end{description}
\end{rem}
%MARGINALE DE MU BAR - PAR RAPPORT A QUELLE VARIABLE 
\textbf{Assumption (C($\nu$))} We suppose the existence of a positive Borel measure $\nu$ on $\R_+$ such that:
\begin{description}
	\item[i)] $R(ds, \infty) << \nu$,
	\item[ii)] The marginal measure of the symmetric measure $|\mb|$ is absolutely continuous with respect to $\nu$. 
\end{description}
If Assumption (C($\nu$)) is realized with $\nu(ds) = |R|(ds,\infty)$ then we will simply say that \textbf{Assumption (C)} is fulfilled.

\begin{prop} \label{pC}
Suppose that $X$ has a covariance measure structure and $\m = \frac{\pa^2 R}{\pa s_1 \pa s_2}$ has compact support, then Assumption (C($\nu$)) is fulfilled with $\nu$ being the marginal measure of $|\mu|$.
\end{prop}
\textbf{Proof}: Let $f: \R_+ \longrightarrow \R$ be a bounded non-negative Borel function.
\begin{description}
	\item[i)] $\left| \int_{\R_+} f(s) |R|(ds,\infty) \right| \leq \int_{\R_+ \times \R_+} f(s_1) d|\mu|(s_1,s_2) = \int_{\R_+} f(s) d\nu(s).$ Take $f$ being the indicator of a null set related to $\n$.
	\item[ii)] $\left| \int_{\RR}f(s_1) d|\mb|(s_1,s_2)\right| \leq k \left| \int_{\RR} f(s_1)d|\mu|(s_1,s_2) \right| = k \int_0^\infty f d\nu,$ where $k$ is the diameter of the compact support of $\m$.
\end{description}
\begin{cor} \label{cC} If $X$ has a covariance measure $\mu$, which is non-negative and with compact support, then Assumption (C) is verified. 
\end{cor}
\textbf{Proof:} This follows because $|R|(ds,\infty)= R(ds,\infty)$ is the marginal measure of $\m$. \qed
\newline 
Next proposition is technical but useful.
\begin{prop} \label{p1.2}
We suppose Assumptions (A), (B). Let $f: \R_+ \lra \R$ be a bounded variation and $\fr$-H\"older continuous function with compact support. Then
\be
\label{1.3} \int_{\RR} R(s_1,s_2) df(s_1)df(s_2) = \int_{\R_+} f^2(s) R(ds, \infty) - \frac{1}{2} \int_{\RR \backslash D}(f(s_1) - f(s_2))^2 d\mu(s_1,s_2).
\ee
\end{prop}
\begin{rem} \label{RP1.2} The statement holds of course if $f$ is Lipschitz with compact support.
\end{rem}
\textbf{Proof:} 
 a) We extend $R$ to the whole plane by continuity. We suppose  first $f \in C_0^\infty (\R_+)$. Let $\rho$ be a smooth, real function with compact support and $\int \rho(x)dx = 1$. We set $\rho_\eps(x) = \frac{1}{\eps} \rho(\frac{x}{\eps})$, for any $\eps > 0$. The left-hand side of (\ref{1.3}) can be approximated by
\be
\label{1.3a} \int_{\RR} R_{\eps}(s_1,s_2) df(s_1)df(s_2),
\ee
where
\[
R_{\eps} = (\rho_{\eps}\otimes \rho_{\eps})* R.
\]
We remark that $R_\eps$ is smooth and
\be
\label{1.3b} \pRe = (\rho_{\eps} \otimes \rho_{\eps}) * \pR,
\ee
where we recall that $\pR$ is a distribution.  By Fubini's theorem on the plane, (\ref{1.3a}) gives
\[
\int_{\R^2} \pRe (s_1,s_2)f(s_1)f(s_2) ds_1 ds_2.
\]
Let $\he \in C^\infty (\R)$ such that $\he = 1$ for $|x| \leq \frac{1}{\eps}$ and $\he(x) = 0$ for $|x| \geq \frac{1}{\eps} + 1$. Moreover we choose $\eps > 0$ large enough such that $[-\frac{1}{\eps}, \frac{1}{\eps}]$ includes the support of $f$. We have
\[
\int_{\R^2} \pRe (s_1,s_2)f(s_1)f(s_2) \he(s_1) \he(s_2) ds_1 ds_2 = \fr(-I_1(\eps) + I_2(\eps) + I_3(\eps)),
\]
where
\bea
I_1(\eps) & = & \int_{\R^2} \pRe(s_1,s_2)(f(s_1) - f(s_2))^2 \he(s_1) \he(s_2) ds_1 ds_2,\nn \\
I_2(\eps) & = & \int_{\R^2} \pRe(s_1,s_2)f(s_1)^2 \he(s_1) \he(s_2) ds_1 ds_2,\nn \\
I_3(\eps) & = & \int_{\R^2} \pRe(s_1,s_2)f(s_2)^2 \he(s_1) \he(s_2) ds_1 ds_2.\nn 
\eea
Since $\he = 1$ on supp$f$, $I_2(\eps)$ gives
\bea
&&-\iin ds_1 (f^{2})'(s_1) \iin ds_2 \frac{\pa R_{\eps}}{\pa s_2}(s_1,s_2) \he(s_2)\nn \\ &&= -\iin ds_1(f^2 \he)'(s_1) \iin R_{\eps}(s_1,ds_2)\he(s_2) = -(I_{2,1} + I_{2,2})(\eps), \nn
\eea
where
\bea
I_{2,1}(\eps) &=& \iin ds_1 (f^2)'(s_1) \iin R_{\eps}(s_1,ds_2)(\he-1)(s_2),\nn\\
I_{2,2}(\eps) &=& \iin ds_1 (f^2)'(s_1) \iin R_{\eps}(s_1,ds_2).\nn
\eea
$I_{2,1}(\eps)$ is bounded by
\be
\label{1.4d} \iin ds_1\left|\left(f^2\right)'\right|(s_1) \int_{\frac{1}{\eps}}^{\infty}|R_{\eps}|(s_1,ds_2).
\ee
For each $s_1 \geq 0$, we have
\be
\label{1.4c}
\int_{\frac{1}{\eps}}^{\infty}|R_{\eps}|(s_1,ds_2) = \int_{\frac{1}{\eps}}^{\infty}\left|\frac{\pa R_{\eps}}{\pa s_2}(s_1,s_2) \right| ds_2.
\ee
Now
\bea
\frac{\pa R_{\eps}}{\pa s_2}(s_1,s_2) &=& \int_{\R^2} dy_1 dy_2 R(y_1, y_2)\rho_{\eps}(s_1 - y_1)\rho_{\eps}'(s_2 - y_2)\nn \\&=& \iin dy_1 \rho_{\eps}(s_1 - y_1) \iin R(y_1, dy_2) \rho_{\eps}(s_2 - y_2).\nn
\eea
Using Fubini's theorem,(\ref{1.4c}) gives
\bea
&&\int_{\frac{1}{\eps}}^{\infty} ds_2 \iin dy_1 \rho_{\eps}(s_1 - y_1) \iin R(y_1,dy_2) \rho_{\eps}(s_2 - y_2)\nn\\ 
&& \ \ \ \ \ \ \ \ = \int_{0}^{\infty} dy_1 \rho_{\eps}(s_1 - y_1)\iin R(y_1,dy_2) \int_{\frac{1}{\eps}}^{\infty} ds_2 \rho_{\eps}(s_2 - y_2).\nn
\eea
But
\[
\int_{\frac{1}{\eps}}^{\infty} ds_2 \rho_{\eps}(s_2 - y_2) = \int_{\frac{1}{\eps}-y_2}^{\infty} ds_2 \rho_{\eps}(s_2).
\]
Let $M > 0$ such that $supp f \subset [-M,M]$.
Hence (\ref{1.4d}) is bounded by
\bea
&& \sup\left|\left(f^2\right)' \right| \int_0^M dy_1 \iin ds_1 \rho_{\eps}(s_1-y_1)\iin|R|(y_1,dy_2) \int_{\frac{1}{\eps}\left(\frac{1}{\eps} -y_2 \right)}^{\infty}
ds_2 \rho(s_2)\nn \\
&& \ \ \ \ \ \ \ \ \   \leq \sup\left|\left(f^2\right)'\right| \int_0^M dy_1 \iin|R|(y_1,dy_2) \int_{\frac{1}{\eps}\left(\frac{1}{\eps} -y_2 \right)}^{\infty}
ds_2 \rho(s_2) \nn
\eea
because $\int_{-\infty}^{\infty} \rho_{\eps}(y)dy = 1$. This is bounded by $\left(I_{2,1,1}(\eps) + I_{2,1,2}(\eps) \right) \sup \left| \left(f^2 \right)'\right|$ with
\bea
I_{2,1,1} (\eps)& =& \int_0^M dy_1\int_0^{\frac{1}{\eps}} |R|(y_1,dy_2) \int_{\left(\frac{1}{\eps}-y_2 \right)\frac{1}{\eps}}^{\infty} ds_2 \rho(s_2), \nn \\
I_{2,1,2} (\eps) &= & \int_0^M dy_1 \int_{\frac{1}{\eps}}^{\infty}|R|(y_1,dy_2). \nn
\eea
Both expressions above converge to zero because of Assumption (A) ii) and Lebesgue dominated convergence theorem.
Hence $I_{2,1}(\eps) \to 0$.
\par As far as $I_{2,2}(\eps)$ is concerned, when $\eps \to 0$ we get
\[
\iin df^2(s_1)R_{\eps}(s_1,\infty) \to \iin df^2(s_1)R(s_1,\infty) = 
-\iin f^2(s_1) R(ds_1,\infty)
\]
according to Assumption (A) i). Consequently 
$\lim_{\eps \to 0} I_2 (\eps) = \iin f^2(s)R(ds,\infty)$. Since $I_2(\eps) = I_3(\eps)$, we also have $\lim_{\eps \to 0}I_3 (\eps) = \iin f^2(s)R(ds,\infty)$. 
\par
It remains to prove that $\lim_{\eps \to 0}I_1(\eps)= \int_{\RR \backslash D} (f(s_1)-f(s_2))^2 d\m(s_1,s_2)$. By (\ref{1.3b}), transferring the convolution against $\rho_{\eps} \otimes \rho_{\eps}$ to the test function, $I_1(\eps)$ becomes the expression 
\[
\left\langle  \pR, (f^{\eps}(s_1)-f^{\eps}(s_2))^2 \he^{\eps}(s_1) \he^{\eps}(s_2)\right\rangle,
\]
where
\bea
f^{\eps} &=& f * \rho^{\eps},\nn \\
\he^{\eps}& = & \he * \rho^{\eps}. \nn
\eea
This gives
\be
\int_{\RR}d\mb(s_1,s_2)\frac{(f^{\eps}(s_1)- f^{\eps}(s_2))^2}{(s_1 - s_2)}\he^{\eps}(s_1)\he^{\eps}(s_2). \label{FApprH}
\ee
We observe that the functions
\be
 \label{FHoldergn} 
 g_{\eps} (s_1,s_2)= \left\{
\begin{array}{clrr}
\frac{(f^{\eps}(s_1)-f^{\eps}(s_2))^2}{s_1-s_2} &\ \ , s_1 \neq s_2\\
0 &\ \ , s_1=s_2.
\end{array} \right.
\ee
and 
\be
\label{FHolderg}
g(s_1,s_2)= \left\{
\begin{array}{clrr}
\frac{(f(s_1)-f(s_2))^2}{s_1-s_2} &\ \ , s_1 \neq s_2\\
0 &\ \ , s_1=s_2.
\end{array} \right.
\ee
are bounded by the square of the $\fr$-H\"older norm of $f$.
\newline Using Lebesgue's dominated convergence theorem, \eqref{FApprH} goes to
\[
\int_{\RR}d\mb(s_1,s_2)\frac{(f(s_1)- f(s_2))^2}{(s_1 - s_2)} = \int_{\RR\backslash D}d\mu(s_1,s_2)(f(s_1)-f(s_2))^2.
\]
This justifies the case $f \in C^\infty_0(\R)$. We consider now the general case. Let $\rho_n$ be a sequence of mollifiers converging to the Dirac delta function and we set $f_n = \rho_n * f$. Taking into account previous arguments, identity \eqref{1.3} holds for $f$ replaced with $f_n$. Therefore we have
\bea
\int_{\RR} R(s_1,s_2) df_n(s_1)df_n(s_2) &= &\int_{\R_+} f_n^2(s) R(ds, \infty) \nn \\
\label{1.3n}
\\ &-& \frac{1}{2} \int_{\RR \backslash D}(f_n(s_1) - f_n(s_2))^2 d\mu(s_1,s_2). \nn
\eea
 The total variation of $f_n$ is bounded by a constant times the total variation of $f$ and $f_n \lra f$ pointwise. So $df_n \lra df$ weakly and also $df_n \otimes df_n \lra df \otimes df$ by use of monotone class theorem. Since $R$ is continuous and bounded,
\[
\lim_{n \to \infty} \int_{\RR} R(s_1,s_2) df_n(s_1)df_n(s_2) = \int_{\RR} R(s_1,s_2) df(s_1)df(s_2).
\]
Therefore the left-hand side of \eqref{1.3n} converges to the left-hand side of \eqref{1.3}.
On the other hand $f_n(s_1) - f_n(s_2) \lra_{n\to \infty} f(s_1) - f(s_2)$ for every $(s_1,s_2) \in \RR$. Since $\bar{\mu}$ and $|R|(ds,\infty)$ are finite non-atomic measures and because of considerations around \eqref{FHoldergn} and \eqref{FHolderg}, the sequence of right-hand sides of (\ref{1.3n}) converges to the right-hand side of (\ref{1.3}).\qed  
\newline \textbf{Assumption (D)} 
\begin{description}
	\item[i)] $R(ds_1, \infty)$ is a non-negative, $\sigma$-finite measure,
	\item[ii)] $\left.\pR \right|_{\RR \backslash D}$ is a non-positive measure.
\end{description}
In the next section, we will expand some examples of processes for which Assumptions (A), (B), (C) and (D) are fulfilled.

\section{Examples} \label{s2}

\setcounter{equation}{0}

\subsection {Processes with covariance measure structure} 
	The first immediate example arises if $X$ has a covariance measure $\m$ with compact support, see Definition \ref{d0.1}. In this case $\pR$ is a measure $\mu$. In Remark \ref{r1.1} we proved that Assumption (A) is satisfied; obviously Assumption (B) is fulfilled and $\bar{\mu}$ is absolutely continuous with respect to $\mu$. We recall that $D$ is the diagonal of the first quarter of the plane. 
	\begin{rem} \label{r1.3}
	i) If $\m$ restricted to $D$ vanishes, Assumption (D) cannot be satisfied except if process $X$ is deterministic. Indeed, suppose that for some $a > 0$, $X_a$ is non-deterministic. Then, taking $f= 1_{[0,a]}$, we have
	\[
	\int_{\RR} f(x_1)f(x_2) d\mu(x_1,x_2) = Var(X_a),
	\]
	which is strictly positive. On the other hand previous integral equals
	\[
	\int_{\RR \backslash D} f(x_1)f(x_2)d\mu(x_1,x_2) = \m([0,a]^2) \leq 0.
	\]
	\newline
	ii) However, Assumption (D) is satisfied if $supp \;\mu \subset D$ and $R(ds_1,\infty)$ is positive. In this case $\left. \mu \right|_{\RR \backslash D}$ is even zero. Consider as an example the case of classical Brownian motion or a martingale.
	\end{rem}

\subsection{Fractional Brownian motion} \label{s4.2}
	\par 
	Let $X = B^H$, $0 < H < 1$, $H \neq \frac{1}{2}$ be a fractional Brownian motion with Hurst parameter $H$ stopped at some fixed time $T > 0$. Therefore we have $X_t = X_T$, $t \geq T$. Its covariance is
	\[
	R(s_1, s_2) = \frac{1}{2}(\tilde{s}_1^{2H}+\tilde{s}_2^{2H}-|\tilde{s}_2-\tilde{s_1}|^{2H} )
	\]
	with $\tilde{s}_i = s_i \wedge T$.
	Now
	\[
	s_1 \mapsto R(s_1, \infty) = R(s_1, T)
	\]
	has bounded variation and is even absolutely continuous since
	\[
	 \frac{\pa R}{\pa s_1} (s_1, \infty) = \left\{
\begin{array}{clrr}
H[s_1^{2H-1}+(T - s_1)^{2H-1}] & \textrm{ if } s_1 < T\\
 0 & \textrm{ if } s_1 > T.
\end{array} \right.
	\]
	So Assumption (A) is verified. Assumptions (B) and (C) are fulfilled because
	\bea
	\bar{\mu}(ds_1,ds_2) &=& (s_1-s_2) \frac{\pa^2 R(s_1,s_2)}{\pa s_1 \pa s_2} \nn \\ 
	&=&  H(2H-1)|s_1 - s_2|^{2H-1}1_{[0,T]^2}(s_1,s_2) sign(s_1 - s_2) ds_1 ds_2 , \nn\\
	\left. \pR \right|_{\R^2 \backslash D} & =& H(2H-1)|s_1 - s_2|^{2H-2} 1_{[0,T]^2}(s_1,s_2). \nn
	\eea
	\begin{rem} \label{r1.4}
	In this example $R(ds, \infty)$ is non-negative; $\mu$ is non-positive if and only if $H \leq \frac{1}{2}$. In that case Assumption (D) is fulfilled.
	\end{rem}
	
\subsection{Bifractional Brownian motion} \label{4.3}
	\par Suppose that $X = B^{H,K}$ is a bifractional Brownian motion with parameters $H \in ]0,1[$, $K \in ]0,1]$ stopped at some fixed time $T >0$. We recall that $B^{H,1}$ is a fractional Brownian motion with Hurst index $H$. Moreover its covariance function, see \cite{RT} and \cite{HV}, is given by
	\[
	R(s_1,s_2) = 2^{-K} \left[ (\tilde{s_1}^{2H}+\tilde{s_2}^{2H})^K - |\tilde{s_1} - \tilde{s_2}|^{2HK} \right]
	\]
	again with $\tilde{s_i} = s_i \wedge T$.
	\par
	We have,
	\be \label{1.3bis}
	\frac{\pa R}{\pa s_1}(s_1,s_2) = 2HK 2^{-K}\left[ \left[s_1^{2H} + s_2^{2H}\right]^{K-1}s_1^{2H-1} - |s_1 - s_2|^{2HK - 1} sign (s_1 - s_2) \right]
	\ee
	for $s_1,s_2 \in ]0,T[$.
	Hence
	\bea
	&&\left.\pR \right|_{[0,T]^2 \backslash D} \nn \\
	&&%\ \ \ \ \ \ \ 
	= 2^{-K} \left[(4H^2K(K-1)(s_1^{2H}+s_2^{2H})^{K-2}(s_1s_2)^{2H-1} %\nn \\&& \ \ \ \ \ \ \ \ \ \ \ \ \ \ %
	 + 2HK(2HK - 1)|s_1 - s_2|^{2HK - 2}\right].\nn
	\eea
	Consequently
	\[
	 \frac{\pa R}{\pa s_1} (s_1, \infty) = \left\{
\begin{array}{clrr}
2HK2^{-K} \left[ (s_1^{2H}+T^{2H})^{K-1}s_1^{2H-1} + (T-s_1)^{2HK-1})\right]& \textrm{ if } s_1 \in ]0,T[\\
 0 & \textrm{ if } s_1 > T.
\end{array} \right.
	\]
	Moreover
	\bea
	&&\bar{\mu}(ds_1,ds_2)  \nn \\
	&& \ \ \ \ \ \  =1_{[0,T]^2}(s_1,s_2) 2^{-K} \left[ 4H^2K(K-1)(s_1^{2H}+s_2^{2H})^{K-2}(s_1s_2)^{2H-1}(s_1 - s_2)^2 \right. 
	 \nn \\
	 && \ \ \ \ \ \ \ \ \ \ \ \ \ \ \ \ \ \ \ \ \ \ \ \ \ \ \ \ \ \ \ \ \ \left.+2HK(2HK - 1)|s_1 - s_2|^{2HK}ds_1 ds_2\right].\nn
	\eea
	\textbf{Conclusions}
	
\begin{description}
	\item[i)] Assumptions (A) and (B) are verified. Assumption (C) is verified because $R(ds,\infty)$ and the marginal measure of $|\bar{\mu}|$ are equivalent to Lebesgue measure on $]0,T]$ and they vanish on $]T,\infty[$. If $HK \geq \frac{1}{2}$, $X$ has even a covariance measure $\mu$, see \cite{KRT}, Section 4.4.
	\item[ii)] Assumption (D) is verified only if $HK \leq \frac{1}{2}$. Indeed, $R(ds, \infty)$ is non-negative and $\left. \pR\right|_{\RR \backslash D}$ is non-positive.
\end{description}
\begin{rem} \label{rrem}
If $HK = \fr$, then Assumption (D) is verified even if $K \neq 1$. In that case $B^{H,K}$ is not a semimartingale, see \cite{RT}, Proposition 3. This shows existence of a finite quadratic variation process which verifies Assumption (D) and it is not a local martingale. %Hence that Assumption is verified also for some finite quadratic variation process which is not a local martingale. 
\end{rem}

\subsection{Processes with weak stationary increments} \label{s4.4}

\begin{df} \label{d1.5}
A square $(\tilde{X}_t)_{t \geq 0}$, such that $\tilde{X}_0 = 0$, is said \textbf{with weak stationary increments} if for every $s, t , \t \geq 0$
\[
Cov(\tilde{X}_{s+\t}-\tilde{X}_{\t},\tilde{X}_{t+\t} - \tilde{X}_{\t}) = Cov (\tilde{X}_s,\tilde{X}_t).
\]
In particular setting $Q(t) = Var(\tilde{X}_t)$ we have
\[
Var(\tilde{X}_{t+\t}-\tilde{X}_{\t}) = Q(t), \  \forall t\geq 0.
\]
\end{df}
In general $\tilde{X}$ does not fulfill the technical assumption (\ref{1.1}) and therefore we will work with $X$, where $X_t = \tilde{X}_{t \wedge T}$. This is no longer a process with weak stationary increments and its covariance is the following:
\be
	R(s_1,s_2) = \label{ExprRQ}
	\left\{
\begin{array}{clrr}
\frac{1}{2}(Q(s_1)+Q(s_2)-Q(s_1-s_2)) &\ \ , s_1,s_2 \leq T, \\
\frac{1}{2}(Q(s_1)+Q(T)-Q(T-s_1)) &\ \ , s_2>T, s_1 \leq T, \\
\frac{1}{2}(Q(s_2)+Q(T)-Q(T-s_2)) &\ \ , s_2 \leq T, s_1 > T, \\
Q(T) &\ \ , s_1,s_2 > T.
\end{array} \right.
\ee
\begin{rem} \label{r1.6}
Without restriction of generality we will suppose
\be
	\label{1.4} Q(t) = Q(T), \ \ t \geq T 
\ee
so that $Q$ is bounded and continuous  and can be extended to the whole line. 
\end{rem}
\begin{prop} \label{p1.6}
Assumption (A) is verified if $Q$ has bounded variation.
\end{prop}
\textbf{Proof}: We have
\be
	R(\infty,s_2) \label{ExprRQ1}
	 = \left\{
\begin{array}{clrr}
\frac{Q(s_2)+Q(T)-Q(T-s_2)}{2} &\ \ , s_2 \leq T, \\
Q(T) &\ \ , s_2 > T, \\
\end{array} \right.
\ee
so that 
\be
	R(\infty,ds_2) = \label{ExprRQ2}
\left\{
\begin{array}{clrr}
\fr \left(Q(ds_2)-Q(T-ds_2)\right) &\ \ , s_2 \leq T, \\
0 &\ \ , s_2 > T; \\
\end{array} \right.
\ee
\eqref{ExprRQ}, \eqref{ExprRQ1} and \eqref{ExprRQ2} imply the validity of Assumption (A). \qed

\begin{prop} \label{p1.7}
We suppose the following.
\begin{description}
	\item[i)] $Q$ is absolutely continuous with derivative $Q'$.
	\item[ii)] $F_Q(s): = s Q'(s), \ s>0$ prolongates to zero by continuity to a bounded variation function, which is therefore bounded.
\end{description}
Then $\mb$ is the finite Radon measure
\[
1_{[0,T]^2}(s_1,s_2) \left( -Q'(s_1-s_2)ds_1ds_2 + F_Q (s_1) ds_1 \d_0(ds_2)-F_Q (s_1-ds_2)ds_1\right).
\]
Moreover Assumption (B) is verified as well as Assumption (C($\nu$)) with $\nu(ds_2) = 1_{]0,T[}(s_2)ds_2$.
\end{prop}

\begin{rem} \label{r1.8}
\begin{description}
		\item[i)] $F_Q$ can be prolongated to $\R$ by setting $F_Q(-s)= F_Q(s)$, $s \geq 0$.
		\item[ii)]In the sense of distributions we have 
		\[(Q'(s)s)' = Q''(s)s + Q'(s).
		\]
		Under i), ii) is equivalent to saying that $Q''(ds)\cdot s$ is a finite measure.
		\item[iii)] Consequently for any $\r$>0, $Q''|_{]-\infty,-\r] \cup [\r,+\infty[}$ is a finite, signed measure.		
		\end{description}
\end{rem}
\textbf{Proof} (of Proposition \ref{p1.7}):
We will evaluate
\be
\label{1.4b} \left\langle \frac{\pa R^2} {\pa s_1 \pa s_2}(s_1,s_2)(s_1-s_2), \f \right\rangle
\ee
for $\f \in C^\infty_0 (\RR)$. This gives 
\[
\int_{\RR} R(s_1,s_2) \frac{\pa^2}{\pa s_1 \pa s_2}(\f(s_1,s_2)(s_1 -s_2))ds_1 ds_2 = \fr(I_1 + I_2 - I_3),
\]
where
\bea
I_1 &=& \int_{\RR} Q(s_1) \frac{\pa^2}{\pa s_1 \pa s_2} (\f(s_1,s_2)(s_1 -s_2))ds_1 ds_2, \nn \\
I_2 &=& \int_{\RR} Q(s_2) \frac{\pa^2}{\pa s_1 \pa s_2} (\f(s_1,s_2)(s_1 -s_2))ds_1 ds_2, \nn \\
I_3 &=& \int_{\RR} \tilde{Q}(s_1,s_2) \frac{\pa^2}{\pa s_1 \pa s_2} (\f(s_1,s_2)(s_1 -s_2))ds_1 ds_2, \nn
\eea
where
\[
	\tilde{Q}(s_1,s_2) = \left\{
\begin{array}{clrr}
Q(s_1 -s_2) &\ \ , s_1,s_2 \leq T, \\
Q(T-s_1) &\ \ , s_2>T, s_1 \leq T, \\
Q(T-s_2) &\ \ , s_1 > T,s_2 \leq T,\\
0 &\ \ , s_1,s_2 > T.
\end{array} \right.
\]
First we evaluate $I_3$. Using assumption i) it is clear that for every $s_2 \geq 0$, $s_1 \longmapsto \tilde{Q}(s_1,s_2)$ is absolutely continuous. Similarly, for every $s_1 \geq 0$, $s_2 \longmapsto \tilde{Q}(s_1,s_2)$ has the same property. Therefore integrating by parts we obtain
\bea
	I_3 &=& - \int_0^\infty ds_1 \int_0^\infty ds_2 \frac{\pa \tilde{Q}}{\pa s_2}(s_1,s_2) \frac{\pa}{\pa s_1} \left( \f(s_1,s_2)(s_1 -s_2)\right) \nn \\
	&=& \int_0^T ds_1 \left( \int_0^T  ds_2 Q'(s_1 -s_2) \frac{\pa}{\pa s_1} \left(\f (s_1,s_2)(s_1 -s_2)\right) \right) \nn \\
	&+& \int_T^\infty ds_1 \left( \int_0^T ds_2 Q'(T -s_2) \frac{\pa}{\pa s_1} \left(\f (s_1,s_2)(s_1 -s_2)\right) \right). \nn
\eea
Using Fubini's theorem we get
\bea
&& \int_0^T ds_2 \int_0^T ds_1 Q'(s_1 -s_2)\left( \frac{\pa}{\pa s_1} \f(s_1,s_2) (s_1-s_2) + \f(s_1,s_2) \right) \nn\\
&& \ \ \ - \int_0^T ds_2 Q'(T-s_2) \f(T,s_2)(T-s_2). \nn
\eea
Therefore
\bea
I_3& =& - \int_0^T ds_2 \int_0^T ds_1 F_Q(ds_1-s_2)\f(s_1,s_2) + \int_0^T ds_2 \left\{F_Q(T-s_2) \f(T,s_2) - F_Q(-s_2) \f(0,s_2) \right\}\nn\\
&& \ + \int_0^T ds_2 \int_0^T ds_1 Q'(s_1 -s_2) \f(s_1,s_2)- \int_0^Tds_2 Q'(T-s_2) \f(T,s_2)(T-s_2). \nn
\eea
%Since $F_Q$ has bounded variation and $Q'$ is in $L^1_{loc}$ we have shown that
%\be
%\nn |I_3| \leq const(Q) \|\f \|_\infty.
%\ee
Consequently
\bea
I_3& =& - \int_0^T ds_2 \int_0^T ds_1 F_Q(ds_1-s_2)\f(s_1,s_2) \nn\\
\label{ExprI3}
\\
&& \ + \int_0^T ds_2 \int_0^T ds_1 Q'(s_1 -s_2) \f(s_1,s_2)- \int_0^T ds_2 F_Q(s_2) \f(0,s_2) \nn.
\eea
Concerning $I_1$ we obtain
\bea
&&\int_0^\infty ds_1 Q(s_1) \frac{\pa}{\pa s_1}\left.(\f(s_1,s_2)(s_1,s_2))\right|_{s_2=0}^\infty = -\int_0^\infty  ds_1 Q(s_1) \left(  \frac{\pa \f}{\pa s_1}\f(s_1,0)s_1+\f(s_1,0)\right)\nn\\
\label{ExprI1}
\\
&&= - \iin ds_1 Q(s_1) \frac{d}{ds_1} \left( \f(s_1,0) s_1 \right)= \int_0^T Q'(s_1) s_1 \f(s_1,0)ds_1 = \int_0^T F_Q(s_1)\f (s_1,0) ds_1.  \nn
\eea
%Again
%\be
%\nn |I_1| \leq const \|\f \|_\infty.
%\ee
%The case of $I_2$ can be discussed identically to $I_1$. \qed
Concerning $I_2$ we obtain
\be
I_2 = -\int_0^T F_Q(s_2) \f(0,s_2) ds_2 \label{ExprI2}.
\ee
Using \eqref{ExprI1}, \eqref{ExprI2} and \eqref{ExprI3} it follows
\bea
I_1+I_2-I_3 &=& \int_{\R_+} F_Q(s_1)\f(s_1,0) ds_1 \nn \\
  &+& \int_0^T ds_2\int_0^T F_Q(ds_1-s_2) \f(s_1,s_2) - \int_{\RR} ds_1 ds_2 \f(s_1,s_2) Q'(s_1-s_2).\nn
\eea
This allows to conclude taking into account the fact, by symmetry,
\[
ds_1 F_Q(ds_1-s_2) = ds_2 F_Q(s_1 - ds_2).
\]
At this point Assumptions (B) and (C($\nu$)) follow directly. \qed

%\begin{rem} \label{r1.9}
%Previous proof helps to provide an explicit expression of the measure $\tilde{\m}$ in terms of $Q$. This shows that its marginal measure restricted to $]0, \infty[$ is absolutely continuous with respect to Lebesgue measure. (VERIFIER)So Assumption (C) is verified with respect to $\nu = Lebesgue$ on $]0, \infty[$.
%\end{rem}
\begin{cor} \label{c1.9}
Under the assumptions of Proposition \ref{p1.7}, Assumption (D) is verified if $Q$ is increasing and $Q''$ restricted to $]0,T[$ is a non-positive measure.
\end{cor}
\begin{rem}
Previous Assumption (D) is equivalent to $Q$ increasing and concave.
\end{rem}
\textbf{Proof} (of Corollary \ref{c1.9}): By \eqref{ExprRQ2}, $R(ds_1,\infty)$ is a non-negative measure if $Q$ is increasing. $Q''$ restricted to $]0,\infty[$ is a Radon measure, hence $Q'$ restricted to $]0,\infty[$ is of locally bounded variation and
\[
\left. \pR \right|_{\RR \backslash D} = Q''(s_1 - ds_2) 1_{]0,T[^2}(s_1,s_2) ds_1
\]
and the result follows. \qed

\begin{cor} \label{c1.8ter} Under the same assumptions as Proposition \ref{p1.7}, if $Q'1_{]0,T[} > 0$ then Assumption (C) is verified.
\end{cor}
\textbf{Proof}: The result follows because in this case $\nu= |R|(\infty,ds)=R(\infty,ds)= \fr (Q'(s)+Q'(T-s))1_{]0,T[}ds$. \qed

\begin{ex}\label{e1.11}
\emph{Processes with weak stationary increments (particular cases).}
\newline Let $(\tilde{X}_t)_{t \geq 0}$ be a zero-mean second order process with weakly stationary increments. We set
\[
Q(t) = Var (\tilde{X}_t)
\]
We refer again to $X_t = \tilde{X}_{t \wedge T}$.
\begin{enumerate}
	\item Suppose
\[
Q(t) = \left\{
\begin{array}{clrr}
t^{2H} &\ \ , t < T\\
T^{2H} &\ \ , t \geq T.
\end{array} \right.
\]
Then Assumptions (A), (B) and (C) are verified. If $H \leq \fr$, Assumption (D) is fulfilled.
	\item We consider a more singular kernel than every fractional scale. We set $T = 1$.
\[
	Q(t) = \left\{
\begin{array}{clrr}
\frac{1}{\log(\frac{1}{t})} &\ \ , 0<t< e^{-2} \\
0 &\ \ , t = 0, \\
\fr &\ \ , t \geq e^{-2}.
\end{array} \right.
\]
Then Assumption (A) is verified since $Q$ is increasing. Moreover $Q$ is absolutely continuous with
\bea
Q'(t) &=& \left\{
\begin{array}{clrr}
\left(\log \frac{1}{t}\right) ^{-2}\frac{1}{t} &\ \ , \ 0<t< e^{-2}, \\
0 &\ \ , \ t > e^{-2}.
\end{array} \right. \nn
\eea
We observe that 
\[
F_Q(t) = t Q'(t)  = \left\{
\begin{array}{clrr}  \left(\log \frac{1}{t}\right)  ^{-2}&\ \ ,  0<t< e^{-2}, \\
0 &\ \ , \ t > e^{-2}.
\end{array} \right. 
\]
and $\lim_{t \to 0} F_Q (t) = 0$. It is not difficult to show that $F_Q$ has bounded variation, therefore Assumption (B) is fulfilled by Proposition \ref{p1.7}. Since $Q' >0$, on $]0,T[$ a.e., Assumption (C) is verified because of Corollary \ref{c1.8ter}. Finally
\bea
Q''(t) &=& \left\{
\begin{array}{clrr} -\left( \log \frac{1}{t} \right)^{-2}t^{-2} \left[1-2\log \frac{1}{t} \right]&\ \ ,  0<t< e^{-2}, \\
0 &\ \ , \ t > e^{-2}.
\end{array} \right. \nn
\eea
Since $Q''$ is negative, Assumption (D) is fulfilled. 
%On the other hand
%\be
%\nn Q''(t)t = -\left(\log\frac{1}{t} \right) ^{-2} \left[ \log \frac{1}{t} + \frac{1}{t}\right]
%\ee
%is integrable, so Assumption (B) is also satisfied.

\end{enumerate}
\end{ex}
\section{Comparison with the kernel approach} \label{s4.5}

We consider a process $X$ continuous in $\LOm$ of the type
\be \label{E1.1} 
X_t= \int_0^t K(t,s)dW_s, t \in [0,T],
\ee
where $K:\RR \lra \R$ is a measurable function, such that for every $t \geq 0$, $\int_0^t K^2(t,s)ds <\infty$.
\begin{rem} \label{rr1.1} $(X_t)_{t \in [0,T]}$ is a Gaussian process with covariance
\[
R(t_1,t_2) = \int_0^{t_1 \wedge t_2} K(t_1,s)K(t_2,s)ds.
\]
\end{rem}
We extend $(X_t)$ to the whole line, setting $X_t = X_T$, $t \geq T$, $X_t = 0$, $t < 0$  . 
\par We want to investigate here natural, sufficient conditions on $K$ so that $X$ has a covariance measure structure. We take inspiration from a paper of Alos-Mazet-Nualart \cite{AMN}, which discusses Malliavin calculus with respect to general processes of type \eqref{E1.1}. That paper distinguishes between the \textbf{regular} and \textbf{singular} case.

\par The aim of this section is precisely to provide some general considerations related to the approach presented in \cite{AMN} in relation to ours. In their regular context, we will show that the process has covariance measure structure. Concerning their singular case, we will restrict to the case that $K(t,s) = \k(t-s)$, $t\geq s\geq 0$, where $\k:\R_+ \lra \R$. We will provide natural conditions so that Assumptions (A) and (B) are verified. We formulate first two general assumptions on $K$.

\textbf{Assumption (K1)} For each $s \geq 0$, $\bar{K}(dt,s) = K(dt,s)(t-s)$ is a finite measure. This implies in particular, for $\eps >0$, 
\be
\label{K1} K(dt,s)1_{(s+\eps, \infty)}(t) \textrm{ is a finite measure.}
\ee

\bigskip

\textbf{Assumption (K2)}
\be
\nn \eps \sup_s K(s+\eps,s) \lra 0.
\ee

\bigskip

Let $T >0$. We extend $K$ to $\tilde{K}: \RR \lra \R$, so that
\be \label{1.10bis}
\tilde{K}(t,s) = \left\{
\begin{array}{clrr}  K(t,s)&\ \ , \ 0<s<t<T, \\
K(T,s) &\ \ , 0<s<T<t,\\
0 &\ \ , \textrm{otherwise}.
\end{array} \right.
\ee
Let $(W_t)_{t \geq 0}$ be a standard Brownian motion. Indeed
\be
 \label{1.10ter} \tilde{X}_t = \int_0^t \tilde{K}(t,s)dW_s, \ t \in \R_+
\ee
extends $X$ by continuity in $\LOm$ from $[0,T]$ to $\R_+$.
In the sequel $\tilde{K}$ and $\tilde{X}$ will often be denoted again by $K$ and $X$. For processes $X$ defined for $t \in [0,T]$, \cite{AMN} introduces two maps $G$ and $G^*$.
Similarly to \cite{AMN}, we define
\be
\nn G: L^2[0,T] \lra L^2[0,T]
\ee
by 
\be
G\f(t)=\int_0^t K(t,s) \f(s) ds, \ t \in [0,T]. \nn
\ee
Let $W^{1,\infty}([0,T])$ be the space of $\f \in L^2([0,T])$ absolutely continuous such that $\f' \in L^{\infty}([0,T])$. We set
\be
G^*:W^{1,\infty}([0,T]) \lra L^2[0,T], \nn
\ee
by
\be
\nn G^*\f(s) = \f(s)K(T,s) + \int_{[s,T]}(\f(t) - \f(s)) K(dt,s). \nn
\ee
We remark that $G^*$ is well defined because (K1) is verified.
\begin{rem} \label{r1.12} In order to better understand the definition of $G^*$, we consider the following "regular" case:
for $s \geq 0$, $t \longmapsto K(t,s)$, $0 \leq s \leq t \leq T$ has bounded variation and
\be
 \sup_{s \in [0,T] }|K|(dt,s) < \infty. \nn
\ee
Then integration by parts shows that
\be
\nn G^*\f(s) = \int_s^T \f(t)K(dt,s),
\ee
since $K(s-,s) =0$.
\end{rem}
\begin{lem} \label{l1.13}
Under Assumptions (K1) and (K2), for $\f \in C^1_0 (\R_+)$ we have
\be \label{1.11} 
 \int_0^T G^*\f dW = \f(T)X_T -\int_0^T X_s d\f_s.
\ee
\end{lem}
\textbf{Proof}: Let $\eps >0$. Since
\bea
\int_s^T|\f(t) -\f(s)||K|(dt,s) &=& \int_s^t \frac{|\f(t)-\f(s)|}{|t-s|}|\bar{K}|(dt,s)\nn\\
&& \ \ \ \ \leq \sup|\f'| |\bar{K}|([s,T],s) < \infty, \nn
\eea
Lebesgue's dominated convergence theorem gives
\be \nn
(G^* \f)(s) = \f(s)K(T,s) + \lim_{\eps \to 0} \int_{s+\eps}^{T} (\f(t) - \f(s))K(dt,s).
\ee
Integration by parts gives
\bea
&&\f(s)K(T,s) + \lim_{\eps \to 0} \left\{ (\f(T)-\f(s))K(T,s) + (\f(s+\eps) - \f(s))K(s+\eps,s)  \right\}\nn\\
&& \ \ \ \ \ \ -\int_{s+\eps}^T \f'(t)K(t,s) dt. \nn
\eea
Again Lebesgue's dominated convergence theorem implies
\be
\f(T)K(T,s) - \int_s^T \f'(t)K(t,s)dt - \lim_{\eps \to 0} (\f(s+\eps) - \f(s))K(s + \eps,s). \nn
\ee
Since $\f \in C^1_0$, Assumption (K2) says that the limit above is zero. Through stochastic Fubini's, the left member of (\ref{1.11}) gives
\bea
&&\f(T) \int_0^T K(T,s) dW_s -\int _0^T dt \f'(t) \int_0^T dW(s) K(t,s) \nn\\
&& \ \ \ \ \ \ =f(T)X_T - \int_0^T X_s d\f(s). \nn
\eea
So the result is proven. \qed

We leave now the general case and consider one assumption stated in \cite{AMN}.
\begin{rem} \label{r1.14} \cite{AMN} considers the following assumption
\be
\label{1.12} \int_0^T |K|(]s,T],s)^2 ds < \infty,
\ee
which characterizes their "regular" context. Proposition \ref{p1.15} below shows that, (\ref{1.12}) implies that $X$ has a covariance measure structure.
\end{rem}
\begin{rem} \label{r1.14bis} Under (\ref{1.12}), assumptions (K1) and (K2) are in particular fulfilled.
\end{rem}
\begin{prop} \label{p1.15}
Let $(X_t)_{t \in [0,T]}$ be a process defined by
\be 
\nn X_t = \int_0^t K(t,s) dW_s, t \in [0,T],
\ee
where $(W_t)_{t \geq 0}$ is a classical Wiener process. Then $X$ has a covariance measure structure if (\ref{1.12}) is verified.
\end{prop}
\textbf{Proof}: We recall that here $K$ (resp. $X$) is prolongated to $\RR$ (resp. $\R$) in conformity with (\ref{1.10bis}) and (\ref{1.10ter}). It is enough to show that there is a constant $C$, such that
\be
\nn \left|  \left\langle \frac{\pa^2 R}{\pa t_1 \pa t_2}, \f\right\rangle \right| \leq C \| \f\|_\infty, \forall \f \in C^\infty_0(\RR).
\ee
Let $\f \in C^\infty_0(\RR)$. We have
\[
\nn X_t = \int_0^\infty K(t,s) dW_s, \ t\geq 0,
\]
with
\[
R(t_1,t_2) = \int_0^\infty K(t_1,s) K(t_2,s) ds.
\]
Indeed, using Fubini's, we have
\bea
\left\langle \frac{\pa^2R}{\pa t_1 \pa t_2}, \f\right\rangle &=& \int_{\RR} R(t_1,t_2) \frac{\pa^2\f}{\pa t_1 \pa t_2}(t_1,t_2) dt_1 dt_2\nn\\
\label{0.2} && \\
\nn &&=\int_0^{\infty} ds \int_{\RR} \frac{\pa^2\f}{\pa t_1 \pa t_2}(t_1,t_2) K(t_1,s) K(t_2,s) dt_1 dt_2.
\eea
Now $K(dt_1,s)K(dt_2,s)$ is a Radon measure on $\RR$ because of Remark \ref{r1.12}. According to \cite{Bill}, Theorem 12.5, its total variation is the supremum over $s=t_0 <t_1 <\cdots<t_N=T$, of 
\bea
\sum_{i,j=0}^{N-1}&&\left|K(t_{i+1},s)K(t_{j+1},s)+K(t_{i},s)K(t_{j},s)-K(t_{i},s)K(t_{j+1},s)-K(t_{i+1},s)K(t_{j},s)\right|\nn\\
 \nn &&=\sum_{i,j=0}^{N-1} \left|\int_{]t_{i},t_{i+1}]\times]t_j,t_{j+1}]}K(dt_1,s)K(dt_2,s) \right| \leq \sum_{i,j=0}^{N-1} \int_{]t_{i},t_{i+1}]} |K|(dt_1,s)    \int_{]t_{j},t_{j+1}]} |K|(dt_2,s)\\
 \nn &&= \left( \int_{]s,T]}|K|(dt,s)\right)^2 = |K|(]s,T],s)^2.
\eea
Hence (\ref{0.2}) equals
\[
\int_0^\infty \int_{\RR}\f(t_1,t_2) K(dt_1,s) K(dt_2,s)
\]
and
\bea
&&\left| \left\langle  \frac{\pa^2 R}{\pa t_1 \pa t_2},\f \right\rangle \right| \leq \| \f \|_\infty \int_0^\infty ds \|K(dt_1,s)K(dt_2,s))\|_{var} \leq C\| \f\|_\infty. \nn
\eea
with $C = \iin |K|^2(]s,T],s)ds$ and $\| \cdot\|_{var}$ denotes the total variation norm. 
  \qed

In order to prepare the sequel, we specify $\frac{\pa^2R}{\pa s_1\pa s_2}$ if $X_t = \int_0^t \k(t-s)dW_s$, where $\k: \R_+ \lra \R$ has bounded variation, supposed cadlag by convention. So we remain for the moment in the regular case.
\begin{rem} \label{r1.16}
\begin{description}
	\item[a)] We prolongate $\k$ to $\k: \R \lra \R$ setting $\k(t) = 0$ if $t <0$.
	\item[b)] $R(t_1,t_2) = \int_0^{t_1 \wedge t_2} \k(t_1 -s)\k(t_2 -s) ds = \int_0^\infty \k(t_1-s) \k(t_2-s) ds.$ 
	\item[c)] If $\k$ has bounded variation then $\k 1_{]\eps,\infty[}$ has bounded variation for any $\eps >0$, which will constitute Assumption (K1') below. It is equivalent to Assumption (K1), when the kernel $K$ is not necessarily homogeneous.
\end{description}
\end{rem}
\begin{lem} \label{l1.19}
For $\phi \in C_0^\infty (\RR)$
\[
\left\langle \frac{\pa ^2 R}{\pa t_1 \pa t_2}, \phi \right\rangle = \left\langle I_1+I_2+I_3+I_4, \phi \right\rangle,
\]
where $I_1, I_2,I_3, I_4$ are the following Radon measures:
\bea
I_1 &=& \k^2(0) dt_1 \d(dt_2 -t_1), \nn\\
I_2 &=& \k(0) 1_{[t_2,\infty[}(t_1)\k(dt_1 -t_2), \nn\\
I_3 &=& \k(0) 1_{[t_1,\infty[}(t_2)\k(dt_2 -t_1), \nn\\
\left\langle I_4 ,\phi \right\rangle &=& \int_{\RR} \k(dt_1)\k(dt_2) \int_0^\infty \phi(t_1+s,t_2+s) ds.\nn
\eea
\end{lem}
\begin{rem} \label{r1.20} If $\k$ has bounded variation, Lemma \ref{l1.19} shows that $X$ has a covariance measure structure.
\end{rem}
In view of the verification of Assumption (B) we have the following result.
\begin{cor} \label{c1.21} Suppose that $\k(0) =0$, $\k$ with bounded variation. Let $\phi \in C_0^\infty (\RR)$. We have 
\[
\left\langle  \frac{\pa^2 R}{\pa t_1 \pa t_2}(t_1,t_2)(t_1 -t_2), \phi \right\rangle = \int_{\RR} \k(dt_1)\k(dt_2)(t_1 - t_2) \int_0^\infty ds\phi(t_1+s,t_2+s).
\]
\end{cor}
\textbf{Proof} (of Lemma \ref{l1.19}): By density arguments we will reduce to the case, where $\phi =\f \otimes \psi$, $\f, \psi \in C_0^\infty (\RR)$. The left-hand side equals
\[
\int_{\RR} R(t_1,t_2) \frac{\pa ^2\f}{\pa t_1 \pa t_2} (t_1,t_2)dt_1 dt_2 = \int_0^\infty dt_1\int_0^\infty dt_2 \f'(t_1) \psi'(t_2) \int_0^{t_1 \wedge t_2} \k(t_1 -s) \k(t_2 -s) ds.
\] 
We recall that by convention we extend $\k$ to $\R$ by setting zero on $]-\infty,0[$. Hence, by Fubini's theorem it equals
\bea
&&\int_0^\infty ds \int_s^\infty dt_1 \f'(t_1) \k(t_1-s) \int_s^\infty \psi'(t_2)\k(t_2 -s) dt_2 \nn \\
&&= \int_0^\infty ds \left\{ -\f(s)\k(0) - \int_s^\infty \f(t_1)\k(dt_1 -s)\right\}\left\{ - \psi(s)\k(0) - \int_s^\infty \psi(t_2) \k(dt_2 -s)\right\} \nn \\
&&= I_1 + I_2 + I_3 + I_4, \nn
\eea
where
\bea
I_1 &=& \k^2(0) \int_0^\infty \f(s)\psi(s) ds = \k(0)^2 \int_{\RR} ds_1 \d(ds_2-s_1) \f(s_1) \psi(s_2),\nn \\
I_2 &=& \int_0^\infty \psi(s)\k(0) \int_s^\infty \f(t_1) \k(dt_1-s) = \int_{\RR} \f(t_1) \psi(s) 1_{[s,\infty[}(t_1)\k(dt_1-s) ds \k(0), \nn\\
I_3 &=&  \int_{\RR} \f(s) \psi(t_2) 1_{[s,\infty[}(t_2)\k(dt_2-s) ds \k(0), \nn \\
I_4 &=& \int_{\RR} \f(t_1) \psi(t_2) \int_{0}^{t_1 \wedge t_2} \k(dt_1 -s)\k(dt_2 -s). \nn
\eea
By Fubini's theorem
\bea
I_4 &=& \int_0^\infty ds  \int_s^\infty \f(t_1) \k(dt_1-s) \int_s^\infty \psi(t_2) \k(dt_2-s) \nn \\
&=& \int_0^\infty ds  \int_0^\infty \f(t_1+s) \k(dt_1) \int_0^\infty \psi(t_2+s) \k(dt_2).\nn
\eea
This concludes the proof. \qed
\bigskip
\par We examine now some aspects related to the singular case. It is of course possible to give sufficient conditions on the kernel $K$, so that $X_t = \int_0^t K(t,s)dW_s$ fulfills Assumptions (A) and (B), however these conditions are too technical and not readable.

So we decided to consider the homogeneous case in the sense that $K(t,s)=\k(t-s)$, $\k: \R \lra \R$, where $\k|_{\R_-}=0$.
% will be supposed with compact support%. 
Clearly the minimal assumption, so that $X$ is defined, is $\k \in L^2([0,t])$, $\forall t \geq 0$. This is equivalent to $\k \in L^2(\R_+)$.

We formulate first an assumption on $\k$.
\smallskip
\smallskip
\newline \textbf{Assumption (K1')} $\k|_{]\eps,\infty[}$ is with bounded variation for any $\eps >0$.
\vskip0.3cm
We recall that this is equivalent to (K1), when $K$ is homogeneous.
\begin{prop} \label{p1.22} Let $\Xt$ be a process defined by
\[
X_t=\int_0^t \k(t-s) dW_s, \ t \geq 0.
\]
We suppose (K1'),(K2)  and moreover
\begin{description}
	\item[a)] $\k$ has compact support,
	\item[b)] 
	\be
	\label{1.17} \sup_{s \geq 0} \int_0^\infty  du \left| \int_0^u \left(\k(dx)\k(s+x-u) - \k(s-u)\right)\right| < \infty.
 	\ee
\end{description}
Then Assumption (A) is fulfilled.
\end{prop}
\begin{rem} \label{r1.18}
\begin{enumerate}
	\item If we assume (K1'), then $\k(dx)$ is a finite measure on $]\eps,\infty[$, so the left-hand side of (\ref{1.17}) is a priori not always finite. Indeed $|\k|(dx)$ on $[0,\infty[$ is only a $\sigma$-finite measure which may be infinite. $\int_0^u \k(dx)( \k(s+x-u)-\k(s-u))$ is evaluated as
	\[
	\lim_{\eps\to 0} \int_{\eps}^\infty \k(dx)( \k(s+x-u)-\k(s-u))
	\]
	\item Assumption (K2) implies here that $\k(\eps)\eps \xrightarrow[\eps \to 0+] {}0$.
\end{enumerate}
\end{rem}
\textbf{Proof} (of Proposition \ref{p1.22}): Let $\alpha \in C_0^\infty (\R_+)$. We want to show the existence of a constant $\mathcal{C}$ such that, for any $s \geq 0$
\be
\label{1.18} \left| \int_0^\infty dx \alpha'(x)R(s,x)\right| \leq \mathcal{C} \|\alpha\|_\infty.
\ee
This would establish the validity of Assumption (A).
The left-hand side of (\ref{1.18}) is given by
\bea
&& \iin dx \a'(x)\iin du \k(s-u)  \k(x-u) = \iin du \k(s-u) \int_u^\infty dx \a'(x)\k(x-u)\nn \\
&& = \int_0^s du \k(s-u) \iin dx \a'(x+u)\k(x) \nn \\
&& = \lim_{\eps \to 0} \int_0^s du \k(s-u) \int_{\eps}^\infty dx \frac{d}{dx}(\a(x+u)-\a(u))\k(x). \nn \\
&& = \lim_{\eps \to 0} \int_0^s du \k(s-u) \left\{ (\a(u+\eps) -\a(u)) \k(\eps) - \int_\eps^\infty \k(dx)(\a(x+u)-\a(u))\right\} \nn \\ 
&& \label{1.19} = \lim_{\eps \to 0} -\int_0^s du \k(s-u) \int_\eps^\infty \k(dx)(\a(x+u -\a(u)))
\eea
since $\lim_{\eps \to 0} \k(\eps) \eps =0$ and $\k \in L^2([0,s])$.

Now (\ref{1.19}) gives
\bea
&&\lim_{\eps \to 0} \int_\eps^\infty \k(dx) \iin du \k(s-u)(\a(x+u)-\a(u)) \nn \\
&& = \lim_{\eps \to 0} \int_\eps^\infty \k(dx) \left\{ \int_x^\infty d\tilde{u} \k(s+x-\tilde{u})\a(\tilde{u})-\iin du \k(s-u)\a(u) \right\}\nn \\
&& = \lim_{\eps \to 0} \int_\eps^\infty du \a(u) \int_\eps^u \k(dx) \k(s+x-u)- \int_\eps^\infty du \a(u) \int_\eps^\infty \k(dx) \k(s-u)\nn \\
&& = \lim_{\eps \to 0} \int_\eps^\infty du \a(u) \int_\eps^u \k(dx) (\k(s+x-u) - \k(s-u)) - \int_\eps^\infty du\a(u) \int_u^\infty \k(dx) \k(s-u)\nn \\ 
&& \ \ \ \ \ \ \ \xrightarrow[\eps \to 0] {} I_1(\a) - I_2(\a), \nn
\eea
where
\bea
I_1(\a)&=&  \int_0^\infty  du \a(u) \int_0^u \k(dx) (\k(s+x-u) - \k(s-u)),  \nn \\
I_2(\a)& = & \int_0^\infty du \a(u) \int_u^\infty \k(dx) \k(s-u). \nn 
\eea
We remark that
\[
|I_2(\a)| \leq \|\a\|_\infty \int_0^\infty du \k^2(u)
\]
because of Cauchy-Schwarz. Moreover
\[
|I_1(\a)| \leq \|\a\|_\infty \int_0^\infty du \left| \int_0^u \k(dx)(\k(s+x-u) -\k(s-u)) \right|.
\]
The right-hand side is bounded because of (\ref{1.17}). \qed
\begin{rem} \label{remm}We remark that (\ref{1.17}) is a quite general assumption. It is for instance verified if
\be
\label{1.19bis} \iin |\k|(dx) \iin |\k(x+u)-\k(u)| du < \infty
\ee
In particular, taking $\k(x)= x^{H-\fr}$, $H>0$, (\ref{1.19bis}) is always verified.
\end{rem}
We go on establishing sufficient conditions so that Assumption (B) is verified.
\begin{prop} \label{p1.20}
We suppose again (K1'). In particular $|\k|_{var}(x): = 
-\int_x^\infty d|\k|(y)$,
 $x>0$ exists. Suppose there is $\delta >0$ with
\be \int_0^\d |\k|^2_{var}(y)dy < \infty \label{ER120}
\ee
Then Assumption (B) is fulfilled.
\end{prop}
\begin{rem} \label{rem1.20} If $\k$ is monotonous
and  $\k(+\infty) = 0$, then \eqref{ER120} is always fulfilled
 since $|\k|_{var}(x)=-\k(x)$,
 which is square integrable.
\end{rem}
\textbf{Proof}(of Proposition \ref{p1.20}): Let $\f \in C^\infty_0(\RR)$. We need to show that
\be
\label{1.21} \left| \int_{\RR} R(t_1,t_2) \frac{\pa^2}{\pa t_1 \pa t_2}(\f(t_1,t_2)(t_1-t_2))\right| \leq const. \|\f \|_{\infty},
\ee
where
\[
R(t_1,t_2) = \int_0^{t_1 \wedge t_2} \k(t_1 -s) \k(t_2-s) ds.
\]
The left-hand side of (\ref{1.21}) is the limit when $\eps \to 0$ of
\be
\label{1.22} \int_{\RR} R_\eps(t_1,t_2) \frac{\pa^2}{\pa t_1 \pa t_2}(\f(t_1,t_2)(t_1-t_2))dt_1 dt_2,
\ee
where
\bea
&& R_\eps(t_1,t_2) = \int_0^{t_1 \wedge t_2} \k_\eps(t_1 -s) \k_\eps(t_2-s) ds, \nn \\
&& \k_\eps (u)= 1_{]\eps,\infty[} \k(u), \nn
\eea
$\k_\eps$ being of bounded variation. Applying Lemma \ref{l1.19} and the fact that $\k_{\eps}(0)=0$, expression (\ref{1.22})  gives
\be
\label{1.23} \int_{\RR} \k_\eps(dt_1) \k_\eps(dt_2) (t_1 -t_2) \iin ds \f(t_1+s, t_2+s)
\ee
We set $|\k_{\eps}|_{var}(x) : = - \int_{x \vee \eps} ^\infty d|\k|(y)$.
% the total variation of $\k_{\eps}$ on $[0,t]$. 
Let $M>0$ such that $supp \f \subset [0,M]^2$. Previous quantity is bounded by 
\[
\| \f \|_\infty M^2 \int_{\RR} |\k_\eps|_{var}(dt_1)|\k_\eps|_{var}(dt_2)|t_1-t_2|.
\]
We have
\[
\int_{\RR} |\k_\eps|_{var}(dt_1)|\k_\eps|_{var}(dt_2)(t_1-t_2) = 2\iin |\k_\eps|_{var}(dt_1)| \int_0^{t_1} |\k_\eps|_{var}(dt_2)(t_1-t_2).
\]
Integrating by parts, previous expression equals
\bea
&&2\iin|\k_\eps|_{var}(dt_1) \iin |\k_\eps|_{var}(t_2) dt_2 = 2 \iin dt_2 |\k_\eps|_{var} (t_2) \int_{t_2}^{\infty} |\k_\eps|_{var}(dt_1) \nn \\
&&= 2 \int_{\eps}^{\infty} dt_2 |\k|^2_{var}(t_2) \xrightarrow[\eps \to 0] {} 2 \iin dt_2 |\k|^2_{var} (t_2),
\eea
which is finite because of Assumption \eqref{ER120}. \qed
%Integrating by parts we get (A verifier!!!!!!!!!!!!!!!!!!!!!!!!!!!!!!!!)
%\bea
%&&2\iin |\k_\eps|(dt_1)| -\left\{ -\int_0^{t_1} |\k_\eps|(dt_2)(t_1-t_2)\right\} \nn \\
%&& \leq 4 \iin |\k_\eps|(dt_1) \int_0^{t_1} |\k_\eps|(dt_2)t_2 dt_2 + 4 \iin |\k_\eps|(dt_1) t_1 \int_0^{t_1}|\k_\eps|(dt_2)d_2 \nn \\
%&& = \iin|\k_\eps|^2(t_1)t_1(dt_1) + 4 \iin |\k_\eps|(dt_1)\left\{\int_0^{t_1}|\k_\eps|(t_2)dt_2 + t_1|\k_\eps|(t_1)\right\}dt_1 \nn \\
%&& = 8 \iin |\k_\eps|^2(t_1)t_1 dt_1 + 2 \left( dt_1 |\k_\eps|(t_1) \right)^2 \nn\\
%&& \leq 8 \iin |\k|^2(t)tdt + 2\left( \iin dt |\k|(t)\right)^2 \label{1.24}.
%\eea
%Let $\d >0$. Using (K1'), we have $|\k|(t) \leq Var([\d,\infty[)$, $t \geq \d$. Moreover $t \lra |\k|(t)$ has compact support. Hence (\ref{1.24}) is finite if there is $\d >0$ with
%\[
%\int_0^\d |\k|^2(t)t dt + \int_0^\d dt|\k|(t) <\infty.
%\]
%This is the case because of assumptions 1) and 2). \qed

%\begin{rem} \label{r1.21}
%\begin{description}
%	\item[a)] Assumptions 1) and 2) are verified if $\k(x)\tilde x^{2H-1}$, as $x \lra 0$, $0<H<\fr$.
%	\item[b)] If $\sup_\eps |\k|(\eps)\eps \lra 0$ (for instance if $\k$ is monotone at zero and (K2) is verified), then assumption 2) is fulfilled. 
%\end{description}
%\end{rem}

\section{Definition of the Paley-Wiener integral} \label{s5}

\setcounter{equation}{0}

\subsection{Functional spaces and related properties} \label{s5.1}
\par The aim of this section is to define a natural class of integrands for the so called Paley-Wiener integral (or simply Wiener integral). Let $X = (X_t)$ be a second order process, i.e. a square integrable process which is continuous in $\LOm$. We suppose moreover $X_0=0$ and $\lim_{t \to \infty} X_t = X_{\infty}$ in $\LOm$, as in \eqref{1.1}. As observed for instance in \cite{Jolis}, the natural strategy is to extend the linear map
\be \label{3.1}
I: C_0^1(\R_+) \longrightarrow \LOm
\ee
defined by $I(f)$:
\be \label{3.2a}
I(f):= \iin f dX = f(\infty)X_{\infty}-\iin Xdf.
\ee
\par In this section, we introduce a natural Banach space of integrands for which the Wiener integral is defined trough prolongation of operator $I$. In the whole section $X$ will be supposed to verify Assumptions (A), (B) and (C($\nu$)) for some $\nu$ by default. 
\par We recall that $\mu = \pR$ restricted to $\RR \backslash D$ is a $\sigma$-finite measure but $\pR$ is only a distribution.
\begin{df} \label{d3.1}
We denote by $\Lr$ the linear space of Borel functions $f: \R_+ \lra \R$, such that
\begin{description}
	\item[i)] $\iin f^2(s)|R|(ds,\infty) < \infty $,
	\item[ii)] $\int_{\RR \backslash D}(f(s_1) - f(s_2))^2d|\mu|(s_1,s_2) < \infty$,
\end{description}
where $|\mu|$ is the total variation measure  of the $\sigma$-finite measure $\m$.
\end{df}
\begin{rem} \label{r3.2}
The integral in ii) equals
\[
\int_{\RR}|f(s_1)-f(s_2)|^2d|\m|(s_1,s_2).
\]
\end{rem}
For $f \in \Lr$ we define
\bea
\|f\|^2_{\cH} &=& \iin f^2(s) R(ds, \infty) - \fr \int_{\RR}(f(s_1) - f(s_2))^2d\mu(s_1,s_2) \label{normH}\\
\|f\|^2_{R} &=& \iin f^2(s) |R|(ds, \infty) + \fr \int_{\RR}(f(s_1) - f(s_2))^2d|\mu|(s_1,s_2) \label{normR}
\eea
%Indeed, $\|\cdot\|_{\cH}$ is associated with the norm of a self-reproducing kernel space. We present here some suitable interpretation of it. 
Let $H_X$ be the Hilbert subspace of $\LOm$ constituted by the closure of ${I(f), f \in C_0^1(\R_+)}$.
\begin{rem}\label{r3.1}
If $f \in C_0^1(\R_+)$, then (\ref{3.2a}) and Proposition \ref{p1.2} give
\[
E(I(f)^2) = E\left(\int_0^\infty X_u df(u)\right)^2 = \int_{\RR}R(s_1,s_2) df(s_1)df(s_2) = \|f\|^2_{\cH}.
\]
So $C_0^1(\R_+)$ equipped with $\|\cdot\|_{\cH}$ is isometrically embedded into $\LOm$.
\end{rem}
Let $\cH$ be an abstract completion of $C_0^1(\R_+)$ with respect to $\|\cdot \|_{\cH}$. $\cH$ will be called "self-reproducing kernel space". The application I: $C_0^1(\R_+) \lra H_X$ uniquely prolongates to $\cH$. If $\cH$ were a space of functions, the prolongation $I$ would be candidate to be called Paley-Wiener integral.
\begin{rem} \label{r3.3}
i) For $f \in \Lr$, we have 
\[
\|f\|_{\cH} \leq	\|f\|_{R},
\]
ii)
$\|\cdot\|_{\cH}$ and $\|\cdot\|_R$ are seminorms on $\tilde{L}_R$ because they derive from the semi-scalar products
\bea
\left<f,g\right>_R &=& \iin (fg)(s)|R|(ds, \infty) + \int_{\RR} d|\mu|(s_1,s_2)(f(s_1)-f(s_2))(g(s_1)-g(s_2))\nn \\
\left<f,g\right>_{\cH}&=& \iin (fg)(s)R(ds, \infty) - \int_{\RR} d\mu(s_1,s_2)(f(s_1)-f(s_2))(g(s_1)-g(s_2))\nn
\eea
%$\Lr$ is a linear space equipped with the semi-scalar products $\left<\cdot,\cdot\right>_{\cH}$ and $\left<\cdot,\cdot\right>_R$.
\newline iii) $\left<\cdot,\cdot\right>_{\cH}$ is a semi-scalar product because it is a bilinear symmetric form; moreover it is positive definite on $C_0^1$ because $\left<\f,\f\right>_{\cH} = E(I(\f)^2)$.
% Moreover, it defines a true inner product in the obvious quotient composed by the classes of functions $\f$. $\left<\cdot,\cdot\right>_R$ is again a semiscalar product because it is bilinear, symmetric form, which is clearly positive defined. Moreover if $\left<\f,\f\right>_R = 0$ then $\left<\f,\f\right>_{\cH}=0$. 
\newline iv) We remark however that $\|f \|_{R} = 0$ implies in any case $f=0 \ |R|(ds, \infty)$ a.e.
\newline v) $\Lr / \sim$ is naturally equipped with a scalar product inherited by $\left<\cdot,\cdot\right>_R$, where $f \sim g$ if $f=g \ |R|(ds, \infty)$ a.e. and 
\[
f(s_1)-f(s_2)-(g(s_1)-g(s_2)) = 0\  |\m| \textrm{ a.e. }.
\]
In particular
\[
f(s_1)-f(s_2) - (g(s_1)-g(s_2))= 0\  |\mb| \textrm{ a.e. }.
\]
However it may not be complete.
\newline vi) The linear space of $\fr$-H\"older continuous functions with compact support $S$ included in $\Lr$. Indeed, if $f$ belongs to such space then
\[
\iin f^2(s)|R|(ds, \infty) \leq \|f\|^2_{\infty} \iin |R|(ds, \infty)
\]
Moreover, expression ii) in Definition \ref{d3.1} is bounded by
\[
k^2 \int_S |s_1 - s_2| d|\mu|(s_1,s_2) = k^2 |\mb|(S\times S),
\]
where $k$ is a H\"older constant for $f$.
\end{rem}
We denote by $L_R$ the closure of $C_0^1$ onto $\Lr$ with respect to $\|\cdot\|_R$.
\begin{rem} \label{r3.4}
$L_R$ is a normed linear space  (as $\Lr$) which is not necessarily complete.
\end{rem}
We repeat that we will not consider the (abstract) completion of $C^1_0(\R_+)$ with respect to norm $\|\cdot \|_{\cH}$ excepted if it is identifiable with a concrete space of functions.
\par Next proposition shows that in many situations $L_R$ is a rich subspace of $\Lr$.
\begin{prop} \label{p3.4bis} Suppose the existence of an even function $\phi: \R \lra \R_+$, such that $d|\mu|$ is equivalent to $\phi(x_1-x_2)dx_1 dx_2$. Let $f:\R_+ \lra \R$ be a bounded Borel function with compact support with at most countable jumps. Then $f \in \Lr \Rightarrow f\in L_R$.
\end{prop}
\textbf{Proof} (of Proposition \ref{p3.4bis}): Let $\rho: \R_+ \lra \R_+$ smooth with $\iin \rho(y)dy = 1$. We set $\rho_{\eps}=\frac{1}{\eps}\rho(\frac{1}{\eps})$, for $\eps > 0$. We show that $f$ is a limit with respect to $\|\cdot \|_R$ of a sequence of smooth bounded functions with compact support (which belong to $L_R$). We consider the sequence $f_n: \R_+ \lra \R$ defined by
\[
f_n(x) = \iin f(x-y)\rho_{\frac{1}{n}}(y)dy = \iin f\left(x-\frac{y}{n}\right) \rho(y)dy.
\]
Clearly $\|f_n\|_{\infty} \leq \|f\|_{\infty}$, $f_n(x) \lra f(x)$, for every continuity point $x$ of $f$.
Therefore 
\[
\iin|f_n(s)-f(s)||R|(ds,\infty) \xrightarrow[n \to \infty] {} 0,
\]
since $|R|(ds, \infty)$ is a non-atomic finite measure.
It remains to prove
\be \label{3.2}
\int_{\RR}(f_n(x_1)-f(x_1)-f_n(x_2)+f(x_2))^2d|\mu|(x_1,x_2) \xrightarrow[n \to \infty] {} 0.
\ee
The left-hand side of (\ref{3.2}) equals
\bea
%&&  \int_{\RR}(f_n(x_1)-f(x_1)-f_n(x_2)+f(x_2))^2d|\mu|(x_1,x_2) \nn\\
&& \int_{\RR}d|\mu|(x_1,x_2) \left[ \int_{\R_+}dy \rho(y) f\left(x_1 - \frac{y}{n}\right)-f(x_1) - f\left(x_2 - \frac{y}{n}\right)+f(x_2)\right]^2 \nn \\
&& \ \ \ \ \ \ \leq \int_{\RR}d|\mu|(x_1,x_2)  \int_{\R_+}dy \rho(y) \left[ f\left(x_1 - \frac{y}{n}\right)-f(x_1) - f\left(x_2 - \frac{y}{n}\right)+f(x_2)\right]^2. \nn
\eea
Last inequality comes from Jensen's. By Fubini's the right-hand side of previous expression equals
\bea
&& \iin dy \rho(y) \int_{\RR}d|\mu|(x_1,x_2) \left[ f\left(x_1 - \frac{y}{n}\right)-f(x_1) - f\left(x_2 - \frac{y}{n}\right)+f(x_2)\right]^2 \nn \\
&& \ \ \ \ \ \ \ = \iin dy \rho(y) I_n(y), \label{3.3}
\eea
where
\[
I_n(y) = \int_{\RR}d|\mu|(x_1,x_2) \left[ f\left(x_1 - \frac{y}{n}\right)-f(x_1) - f\left(x_2 - \frac{y}{n}\right)+f(x_2)\right]^2.
\]
$I_n(y)$ is bounded by
\bea
&& 2 \left[ \int_{\RR} dx_1 dx_2 \phi(x_1-x_2)\left( f\left(x_1 - \frac{y}{n}\right) - f\left(x_2 - \frac{y}{n}\right)\right)^2 \right.\nn \\
&& \label{3.4}\\
&& \ \ \ \ \ \ \ \ \  + \left. \int_{\RR} dx_1 dx_2 \phi(x_1-x_2)(f(x_1)-f(x_2)) ^2\right] \nn
\eea
Setting $\tilde{x}_i= x_i-\frac{y}{n}$, $i=1,2$, in the first integral, (\ref{3.4}) is upper bounded by
\[
 2 \int_{\RR} \phi(x_1-x_2)(f(x_1)-f(x_2))^2 dx_1 dx_2 =
 2 \int_{\RR} d|\mu|(x_1,x_2)(f(x_1)-f(x_2))^2.
\]
 On the other hand $I_n$ converges pointwise to zero. Lebesgue's dominated convergence theorem applied to (\ref{3.3}) allows to conclude (\ref{3.2}). \qed
 
\begin{rem} \label{r3.4ter} 1) We recall that for instance a cadlag function with compact support is a bounded function with at most countable jumps.
\newline 2) The assumption related to $\mu$ in the statement of Proposition \ref{p3.4bis} is for instance verified if $X$ is a fractional (or bifractional) Brownian motion.
\newline 3) If $|R|(ds,\infty)$ is absolutely continuous with respect to Lebesgue, it is easy to show that the statement of Proposition \ref{p3.4bis} holds for every bounded Borel function with compact support $f:\R_+ \lra \R$. In particular no jump condition is required.
\end{rem}

An easy consequence of the definition of the $L_R$-norm is the following.
\begin{prop} \label{p3.4000} Let $g: \R \lra \R$ be a Lipschitz function, $f \in L_R$. Then $g \circ f \in L_R$.
\end{prop}
\textbf{Proof:} $\n(ds) = |R|(ds,\infty)$. Calling $k$ the Lipschitz constant, since $|g|(s)\leq k(1+|s|)$ it follows
\bea
\iin &&g^2(f(s))\n(ds)+ \fr \int_{\RR}d|\mu|(s_1,s_2)(g(f(s_1))-g(f(s_2)))^2 \nn \\
&\leq &k^2 \n(\R_+)+k^2\iin f^2(s)\n(ds)+\frac{k^2}{2} \int_{\RR} d|\m|(s_1,s_2)(f(s_1)-f(s_2))^2\nn \\
&=& k^2 \n(\R_+)+k^2\|f \|^2_R,\nn
\eea
This shows that $g \circ f \in \Lr$.
$g \circ f$ is indeed in $L_R$ since, if $f_n \in C_0^1(\R_+)$ converges to $f \in L_R$, then performing similar calculations as before, we have
\[
\|g \circ f_n -g \circ f \|^2_R \lra 0.
\] \qed

By analogous arguments, we obtain the following result.
\begin{prop} \label{p3.401} 
Let $f,g \in L_R$ and bounded. Then $fg \in L_R$.
\end{prop}
\begin{prop}\label{p3.4.6}
Let $V:\R \lra \R_+$ increasing on $\R_+$ such that
\bea
&&\int_{\RR}V^2(x_1 - x_2)d|\mu|(x_1,x_2) <\infty  \label{3.5}
\eea
Then every Borel function $f:\R \lra \R_+$ with compact support 
such that
\be \label{3.6}
|f(x_1) - f(x_2)| \leq V(x_1 - x_2)
\ee
belongs to $L_R$.
\end{prop}
\begin{cor}  \label{c3.4.6}
Every bounded $\fr$-H\"older continuous function $f$ with compact support belongs to $L_R$.
\end{cor}
\textbf{Proof} of Corollary \ref{c3.4.6}: We apply Proposition \ref{p3.4.6} with $V(x)= |x|^{\fr}$. \qed
\begin{rem}\label{r3.4.7}
\begin{description}
	%\item[i)] If (\ref{3.5}) is verified then $V(0)=0$, at least in the significant case that $\mu$ does not extend to a measure on $\RR$.
	\item[i)] If $V$ is continuous, condition (\ref{3.6}) is equivalent to saying that $V$ is the continuity modulus of $f$.
	\item[ii)] If $V$ is continuous, then $f$ fulfilling (\ref{3.6}) is continuous.
%	\item[iii)] The function $1_{\R_+}$ belongs to $\Lr$. In Proposition \ref{p3.74} we will show that $1 \in L_R$ under Assumptions (A), (B), (C), (D).
\end{description}
\end{rem}
\textbf{Proof}(of Proposition \ref{p3.4.6}):
It follows the same scheme as the proof of Proposition \ref{p3.4bis}. Let $f_n$ be as in that proof.
\begin{description}
	\item[i)] $f_n \lra f$ pointwise and $|f_n| \leq \|f\|_\infty$,
	\item[ii)] $\iin |f_n - f|(s)|R|(ds,\infty)\lra 0$ because of Lebesgue's dominated convergence theorem.
	\item[iii)] (\ref{3.2}) is again a consequence of Lebesgue's; (\ref{3.3}) still holds with
	\[
	|I_n(y)| \leq 2 \int_{\RR} d|\mu|(x_1,x_2) V^2(x_1-x_2).
	\] \qed
\end{description}
\begin{rem} \label{cc3.5} If $X$ has a covariance measure with compact support, then every bounded function with compact support belongs to $L_R$. This follows by Proposition \ref{p3.4.6} taking $V(x) \equiv 2\sup|f|.$ However the property of compact support will be raised in Proposition \ref{r3.73}.
\end{rem}
\begin{prop} \label{r3.73}If $X$ has a covariance measure with compact support, any bounded function still belongs to $L_R$.
\end{prop}
\begin{prop} \label{p3.5} $\Lr$ is a Hilbert space if Assumption (C) is  fulfilled.
\end{prop}
\begin{rem} \label{r3.6}
Suppose that $X$ has a covariance measure with compact support, with corresponding signed measure $\mu $. We link below the present approach with the one in \cite{KRT}.
\newline 1) One alternative approach would be to consider the measure $\nu$ which was introduced in \cite{KRT}, Section 5, i.e. the marginal measure of $|\m|$. The space $L^2(d\nu): = L^2(\R_+,d\nu)$ was a natural space, where the Wiener integral could be defined.
\newline 2) In this paper $\| \cdot\|_R$ is the norm which allows to prolongate the Wiener integral operator $I$. Point 1) suggests that $\| \cdot\|_{R}$ should be somehow related to $\|\cdot \|_{L^2(\nu)}$. By Cauchy-Schwarz, we have
\[
\int_{\RR}|f(s_1)f(s_2)| d|\mu|(s_1,s_2) \leq \iin f^2(s)d\nu(s).
\]
This shows that the seminorm of $L^2(d\nu)$ is equivalent to $\|f \|_{R,\nu}$ where
\[
\|f\|^2_{R,\nu} = \iin f^2(s) \nu(ds) + \fr \int_{\RR}(f(s_1)-f(s_2))^2 d|\mu|(s_1,s_2).
\]
This looks similarly to $\| \cdot\|_R$ norm but they could be different. Indeed, we only have
\[
\| f\|_{L^2(d\nu)} \sim \| f\|_{R,\nu} \geq \| f\|_{R},
\]
which implies that
\be
L^2(d\n) \subset L_R. \label{EL2NUR}
\ee
This implies that $\|\cdot\|_R$ will provide a larger space, where the Wiener integral is defined.
\newline 3) In particular, if $\m$ is non-negative (as for the case $X$ being a fractional Brownian motion with Hurst index $H \geq  \fr$ stopped at some time $T$), we have 
\[
\| \cdot\|_R = \|\cdot \|_{R,\nu} \sim \|\cdot \|_{L^2(d\nu)}.
\]
Consequently by item 2) it follows that $L_R = L^2(d\n)$.
\newline 4) If $\mu$ is non-negative, then $L_R = \tilde{L}_R$ since $C^1_0(\R)$ is dense in $L^2(d\n)$, see Lemma 3.8 of \cite{KRT}.
\end{rem}
\textbf{Proof} (of Proposition \ref{r3.73}): This follows by Remark \ref{r3.6}, point 2). Indeed, any bounded function belongs to $L^2(d\nu)$ because $\nu$ is finite. \qed
\newline \textbf{Proof} (of Proposition \ref{p3.5}): It is enough to show that $\Lr$ is complete. We set $\chi(ds) = |R|(ds,\infty)$. Let $(f_n)$ be a Cauchy sequence in $\Lr$. We recall that
\[
\|f_n-f_m\|^2_R = \iin(f_n - f_m)^2(s)\chi(ds) + \fr \int_{\RR}|\mu|(ds_s,ds_2)(g_n-g_m)^2(s_1,s_2),
\]
where $g_n(s_1,s_2) = f_n(s_1) - f_n(s_2)$. Since both integrals are non-negative, the Cauchy sequence $(f_n)$ is Cauchy in $L^2(d\chi)$ and $(g_n)$ is Cauchy in $L^2(d|\mu|)$. Since $L^2(d\chi)$ is complete, there is $f \in L^2(d\chi)$ being the limit of $f_n$ when $n \lra \infty$. On the other hand, since $L^2(\RR, |\mu|)$ is complete, $g_n$ converges to some $g\in L^2(\RR;d|\m|)$.
\par It remains to show that
\[
g(s_1,s_2) = f(s_1)-f(s_2) \ \mu \ a.e., \ \ s_1,s_2 >0.
\]
Since $f_n \lra f$ in $L^2(d\chi)$, there is a subsequence $(n_k)$ such that $f_{n_k}(s) \lra f(s)$, for $s \notin N$, where $\chi(N)=0$. Consequently for $(s_1,s_2)\in N^c\times N^c$
\be \label{3.10}
f_n(s_1)-f_n(s_2) \lra f(s_1) - f(s_2).
\ee
Moreover, obviously if $s_1 = s_2$, (\ref{3.10}) holds. Hence for $(s_1,s_2) \in (N^c \times N^c) \cup D$,
\[
g(s_1,s_2) = f(s_1) - f(s_2),
\]
where we recall that $D = \left\{ (s,s)| s \in \R_+\right\}$ is the diagonal.
This concludes the proof if we show that $((N^c \times N^c) \cup D)^c$ is $|\mu|$-null.
\par This set equals
\[
(N^c \times N^c)^c \cap D^c
\]
and it is included in $\left((N \times \R_+) \cap D^c\right) \cup \left((\R_+ \times N) \cap D^c\right)$. Since $|\mu|$ and $|\mb|$ are equivalent outside $D^c$, it is enough to show that
\[
|\mb|((N \times \R_+) \cap D^c) = 0 \textrm{ and } |\mb|((\R_+ \times N) \cap D^c) = 0.
\]
Previous quantities are bounded by
\[
|\mb|(N \times \R_+), \textrm{ and } |\mb|(\R_+ \times N),
\]
which coincide with the marginal measures of $|\mb|$ evaluated on $N$. Assumption (C) allows to conclude. \qed
\begin{cor} \label{c3.7} 
If Assumptions  (C), (D) are in force, then $L_R$ is the closure of $C^1_0(\R_+)$ under $\|\cdot \|_{\cH}$; in particular $L_R$ is a "self-reproducing kernel space".
\end{cor}
\textbf{Proof}: According to Assumption (D), we have $\|\cdot\|_R = \|\cdot\|_{\cH}$ for $\psi \in C_0^1(\R_+)$. According to Proposition \ref{p3.5} $L_R$ is a Hilbert space equipped with $\| \cdot \|_{\cH}$, which is by definition the closure of $C_0^1(\R_+)$.   \qed
\begin{rem} \label{r3.8}
1)  We will see that the Paley-Wiener integral can be naturally defined on space $L_R$.
\newline 2) Corollary \ref{c3.7} is interesting because it shows that 
a natural space where the Wiener integral will be defined is complete under Assumptions (A), (B), (C), (D).
\newline 3) This result is of the same nature as the one of \cite{PT}, which shows that the space, where its Wiener integral is defined, is also complete with respect to the norm $\|\cdot \|_{\cH}$ when $X$ is a fractional Brownian motion, with parameter $H \leq \fr$. In Section \ref{s4.2} we have proved that Assumptions (A), (B), (C), (D) are indeed fulfilled in that case.
\end{rem}
\begin{prop} \label{p3.71} We suppose that Assumption (D) is fulfilled. Then any bounded variation function with compact support belongs to $L_R$ and
\be
\label{F371} \|\f \|^2_R = E\left(-\iin X d\f \right)^2.
\ee
\end{prop}
That property does not seem easy to prove in the general case.
\begin{cor}\label{c3.72} We suppose Assumption (D) to be fulfilled.
 Then every step function belongs to $L_R$. In particular, if $t > 0$,
 $1_{[0,t]} \in L_R$  

 and 
\[
\|1_{[0,t]} \|^2_{R} = E(X^2_t).
\]
\end{cor}
\begin{cor} \label{c3.721} Under Assumption (D), if $f: \R_+ \lra \R$ is a bounded variation function with compact support, then
\[
\int_{\RR}(f(t_1)-f(t_2))^2d|\m|(t_1,t_2) <\infty.
\]
\end{cor}
\textbf{Proof} (of Corollary \ref{c3.721}): This follows from Proposition \ref{p3.71} and the fact that $L_R \subset \Lr$. \qed
\newline \textbf{Proof:} (of Proposition \ref{p3.71}): Let $\f$ be a bounded variation function with compact support, defined on $\R_+$. Let $(\r_n)$ be a sequence of mollifiers converging to the Dirac delta function. We set
\[
\f_n = \r_n * \f.
\]
Since $\f_n$ is smooth with compact support, it belongs to $L_R$. Now $d\f_n \lra d\f$ and
\be
d\f_n \otimes d\f_m \xrightarrow[n,m \to \infty] {} d\f \otimes d\f \textrm{ weakly.} \nn
\ee
Since $R$ is continuous it follows that
\be
\int_{\RR}R(s_1,s_2) d\f_n(s_1) d\f_m(s_2) \lra \int_{\RR}R(s_1,s_2) d\f(s_1) d\f(s_2). \label{E3.71}
\ee
By Remark \ref{r3.1}, we have
\bea
\| \f_n\|_R^2 &=& E \left(-\iin X d\f_n \right)^2 \nn \\
\label{E3.711} \\
&=& \int_{\RR} R(s_1,s_2)d\f_n(s_1)d\f_n(s_2). \nn
\eea
Using \eqref{E3.71}, the limit of the right-hand side of \eqref{E3.711} gives
\[
E\left( - \iin X d\f \right)^2.
\]
Again using \eqref{E3.71}, it follows that
\[
\lim_{n,m \to \infty} \|\f_n-\f_m \|^2_R = 0,
\]
so $(\f_n)$ is Cauchy in $L_R$. Since $L_R$ is complete, there is $\psi \in L_R$ such that 
\[
\|\f_n-\psi \|_{R} \xrightarrow[n\to \infty] {} 0.
\]
 On the other hand $\f_n \lra \f$ $R(ds,\infty)$ a.e. since $R(ds,\infty)$ is a non-atomic measure. By Lebesgue's dominated convergence theorem,
\[
\lim_{n \to \infty}  \| \f_n -\f\|_{R(ds,\infty)} = 0.
\]
But we also have
\[
\| \f_n -\psi\|_{R(ds,\infty)} \leq \|\f_n -\psi \|_R \xrightarrow[n \to \infty]{} 0.
\]
By uniqueness of the limit, $\f = \psi$ $R(ds,\infty)$ a.e. and so $\psi = \f$ $ d\n$ a.e. this shows that $\f \in L_R$. The limit of the left-hand side in \eqref{E3.711} gives $\|\f \|^2_{R}$, which finally shows \eqref{F371}. \qed
\par A natural question concerns whether the constant function $1$ belongs to $L_R$. The answer is of course well-known if $X$ has a covariance measure structure with compact support because of Proposition \ref{r3.73}. Again it will be also the case if Assumptions (A), (B), (C) and (D) are fulfilled.
\begin{prop} \label{p3.74} If Assumptions  (C) and (D)  are fulfilled, then $1 \in L_R$ and $\| 1\|_{\cH}^2 = \| 1\|_{R}^2 = R(\infty,\infty) = E(X^2_\infty)$.
\end{prop}
\textbf{Proof}: For $n \in \mathbb{N}^*$, we consider a smooth, decreasing function $\hn : \R_+ \to [0,1]$  which equals 1 on $[0,n]$ and zero on $[n+1,\infty[$ and it is bounded by 1. Clearly $\hn \in L_R$ for any $n$ and
\bea
\| \hn -\hm\|^2_{\cH} &=& E \left( \iin (\hn -\hm)dX \right)^2 \nn\\
&=& E \left( -\iin Xd(\hn -\hm) \right)^2 = I(n,n)+I(m,m)-2I(n,m), \nn
\eea
where
\bea
I(n,m) &=& \int_{\RR}R(s_1,s_2) d\hn(s_1)d\hm(s_2), \nn \\
 &=& \int_{[n,n+1]\times[m,m+1]}R(s_1,s_2) d(1-\hn)(s_1)d(1-\hm)(s_2). \nn
\eea
We have
\[
\inf_{\xi \in [n,n+1]\times[m,m+1]}R(\xi) \leq I(n,m) \leq \sup_{\xi \in [n,n+1]\times[m,m+1]}R(\xi).
\]
Since $\lim_{s_1,s_2 \to \infty}R(s_1,s_2) = R(\infty,\infty)$, if follows that
\[
I(n,m) \xrightarrow[n,m \to \infty] {} R(\infty,\infty)
\]
and $\hn$ is a Cauchy sequence in $\| \cdot\|_{\cH}$. On the other hand $\hn \lra 1$ pointwise when $n \lra \infty$ and in particular a.e. with respect to the measure $|R|(dt,\infty)$. \qed
%\newline We already know by Corollary \ref{c3.72} that under Assumption (D), $1_{[0,t]} \in L_R$ for every $t \in [0,T]$. We evaluate now some quantities related to $\| 1_{[0,t]}\|_R$ which will be used in the sequel.
%\begin{prop} \label{p3.75}
%We assume Assumptions (D) is fulfilled
%\be
%\label{E375} \int_{\RR} d|\m|(s_1,s_2)(1_{[0,t]}(s_1)-1_{[0,t]}(s_2))^2 = 2 Var(X_t)-2 R(t, \infty).
%\ee
%In particular the left-hand side of (\ref{E375}) is bounded in t.
%\end{prop}
%OU EST UTILISEE CETTE PROPOSITION??????
%\textbf{Proof}: Since Assumption (D) is supposed, the left hand side of (\ref{E375}) equals
%\bea
%&&-\int_{\RR} d\m(s_1,s_2)(1_{[0,t]}(s_1)-1_{[0,t]}(s_2))^2 \nn \\
%&& \ \ \ \ = 2 \|1_{[0,t]} \|^2_{\cH} - 2 \int_0^t R(ds,\infty) = 2 Var(X_t) - 2 R(t, \infty) \nn 
%\eea
%because of Corollary \ref{c3.72}. \qed
\newline An important question concerns the separability of the Hilbert space $L_R$.
\begin{prop} \label{p376} Suppose the validity of Assumptions (C) and (D).
%(A), (B), (C($\nu$)), (D).
 Then the Hilbert space $L_R$ is separable. Moreover there is an orthonormal basis $(e_n)$ in $C_0^1$ of $L_R$.
\end{prop}
\textbf{Proof}: We denote by $\mathcal{S}$ the closed linear span $\{1_{[0,t]}, t\geq 0 \}$ into $\Lr$. We first prove
\be
\label{EIdLR} \mathcal{S} = L_R
\ee
a) $1_{[0,t]}$, $\forall t \geq 0$ belongs to $L_R$ because of Corollary \ref{c3.72}, so it follows that $\mathcal{S} \subset L_R$.
\newline b) We prove the converse inclusion. It is enough to show that $C^1_0(\R_+) \subset \mathcal{S}$. Let $\f \in C^1_0(\R_+)$ and consider a sequence of step functions of the type 
\[
\f_n(t)= \sum_{l} 1_{[t_l,t_{l+1}[} (t)\f(t_l),
\]
which converges pointwise to $\f$. Since the total variation of $\f_n$ is bounded by the total variation of $\f$, then $\f_n \lra \f$ weakly. 
Let $(Y_t)$ be a Gaussian process with the same covariance as $X$. A consequence of previous observations shows that
\be 
\label{YnPhi} \iin \f_n dY:= -\iin Y d\f_n \xrightarrow[n \to \infty]{} -\iin Y d\f \textrm{ a.s.}
\ee 
Since $Y$ is a Gaussian process, the sequence in the left-hand side of \eqref{YnPhi} is Cauchy in $\LOm$, so $(\f_n)$ is Cauchy in $L_R$ by Proposition \ref{p3.71}; the result follows because $\f_n \in \mathcal{S}$ for every $n$. This concludes b) and \eqref{EIdLR}.
\par Since $\|\cdot \|_R = \|\cdot \|_{\cH}$, taking into account the consideration preceding Remark \ref{r3.3} and \eqref{EIdLR}, $H_X$ is the closure in $L^2(\Om)$ of $-\iin X d\f$, $\f$ of the type $1_{[0,t]},$ $t \geq 0$. Since $X$ is continuous, $H_X$ (and therefore $L_R$) is separable. The existence of an orthonormal basis in $C^1_0(\R_+)$ follows by Gram-Schmidt orthogonalization procedure. \qed

\subsection{Path properties of some processes with stationary increments} \label{sJ}

In this subsection we are interested in expressing necessary and sufficient conditions under which the paths of Gaussian continuous processes with stationary increments restricted to any compact intervals, belong to $L_R$. We have some relatively complete elements of answer. 
\par We reconsider the example treated in Section \ref{s4.4}. Let $\tilde{X}$ be a process with weak, stationary increments, continuous in $L^2$ such that $\tilde{X}_0= 0$. We denote by $Q(t)= Var(\tilde{X}_t)$, and we consider again $X$ defined by $X_t = \tilde{X}_{t\wedge T}$. We recall that without restrictions to generality, we can suppose $Q(t)= Q(T)$, $t \geq T$.
\par We recall that in Proposition \ref{p1.7} we provided conditions so that Assumptions (A) and (B) are verified, i.e.
\begin{hyp} \label{hAB}
\begin{description}
	\item[i)] $Q$ is absolutely continuous with derivative $Q'$,
	\item[ii)] $F_Q(s):= s Q'(s)$, $s>0$ prolongates to zero by continuity to a bounded variation function. 
\end{description}
\end{hyp}
In Corollary \ref{c1.9} we provided conditions so that Assumption (D) is verified. This gave the following
\begin{hyp} \label{hD}
\begin{description}
	\item[i)] $Q$ is non-decreasing,
	\item[ii)] $Q''$ non-positive $\s$-finite measure. 
\end{description}
\end{hyp} 
\begin{prop} \label{pDom1} We suppose the validity of Hypothesis \ref{hAB}.
If 
\be
\int_{0+} Q(y)|Q''|(dy)< \infty, \label{hAB1}
\ee
then almost all paths of $X$ belong to $\tilde{L}_R$.
\end{prop}
\textbf{Proof}: We recall that Hypothesis \ref{hAB} implies that $Q''$ is a finite Radon measure on $]\d,\infty[$ for every $\d>0$. Hence \eqref{hAB1} implies that
\be
\iin Q(y)|Q''|(dy) < \infty. \label{EQQ2}
\ee
Since $\int_0^{\infty} |R|(ds,\infty)X_s^2 < \infty$ a.s. being $|R|(ds,\infty)$ a finite measure, it remains to prove that
\be
\label{EDom1} \int_{[0,T]^2} (X_{s_1}-X_{s_2})^2 |Q''|(ds_2-s_1)ds_1 < \infty. \ \textrm{ a.s.}
\ee
To prove \eqref{EDom1}, it is enough to evaluate the expectation of its left-hand side. We get
\[
\int_{[0,T]^2}|Q'' |(ds_2-s_1)ds_1Q(s_1-s_2) = 2 \int_{0}^T ds_1 \int_0^{s_1} Q(s_2)|Q''|(ds_2).
\]
This concludes the proof. \qed 
\begin{rem} \label{rDom2} If $X$ has a covariance measure structure, then Assumption \eqref{hAB1} is trivially verified.
\end{rem}
\begin{rem} 1) \label{RpDom2} Assumption \eqref{hAB1} implies \eqref{EQQ2} which ensures \eqref{EPD15} in Proposition \ref{PRD17}. Suppose \eqref{EQQ2},
if Assumptions (C), (D) are fulfilled, then Proposition \ref{PRD17} says that a.s. $X \in L_R$. 
\newline 2) In the sequel we will express necessary conditions.
\end{rem}
\begin{prop} \label{pDom2}
We suppose $X$ Gaussian  and continuous. We suppose again the validity of Hypotheses \ref{hAB}, \ref{hD} and the following technical conditions. There are $c_1,c_2 >0$, $\alpha_1 < 1$, $\a_2 >0$ such that
\be
c_1 t^{\a_1} \leq Q(t) \leq c_2 t^{\a_2}, \label{tech}
\ee
\be
-\int_{0+} Q(y)Q''(dy) = \infty \label{EInfinite}
\ee
then $X \notin \Lr$ a.s.
\end{prop}
\begin{cor} \label{crDom3}
Let $\tilde{X}$ be a continuous mean-zero Gaussian process with stationary increments such that $\tilde{X}_0 = 0$ a.s. $Q(t) = Var \tilde{X}_t$. Set $X_t = \tilde{X}_{t \wedge T}$, $t \geq 0$. We suppose Hypotheses \ref{hAB} and \ref{hD} together with \eqref{tech}. 
\newline Then $X \in L_R$ a.s. if and only if
\[
\int_{0+} Q(y) |Q''|(dy) < \infty.
\]
\textbf{Proof}: It follows from Remark \ref{RpDom2} and Proposition \ref{pDom2}, and the fact that $L_R \subset \Lr$.
\end{cor}
\begin{rem} \label{rDom3}
1) The importance of Corollary \ref{crDom3} is related to the problem of finding sufficient and necessary conditions on the paths of a continuous Gaussian process $X$ to belong to its "self-reproducing kernel space".
\par When it is the case, $X$ belongs to the natural domain of the divergence operator in Malliavin calculus (Skorohod integral); in the other cases $X$ will be shown to belong  the extended domain $Dom \d^*$, see Definition \ref{d431} introduced in the spirit of \cite{CN, MV}.
\newline 2) We conjecture that assumption \eqref{tech} and Hypothesis \ref{hD} can be omitted, but this would have considerably complicated the proof.
\end{rem}
\textbf{Proof} (of Proposition \ref{pDom2}): Since $X$ is continuous, therefore locally bounded, we observe that
\[
\int_0^T X_s^2 |R|(ds,\infty) < \infty \textrm{ a.s.}
\] 
To prove that $X \notin \Lr$ a.s., it will be enough to prove that
\be
\int_{\RR} (X_{s_1}-X_{s_2})^2 |Q''|(ds_2-s_1) = \infty \textrm{ a.s.} \label{E4Dom}
\ee
The left-hand side of \eqref{E4Dom} gives
\bea
&&2 \int_0^T ds_1 \int_0^{s_1} \left(X_{s_1}-X_{s_2} \right)^2(-Q'')(ds_2-s_1) = 2 \int_0^{T} ds_1 \int_{0}^{s_1} \left( X_{s_1} -X_{s_1-s_2}  \right)^2(-Q'')(ds_2) \nn \\
&& = 2 \int_0^T (-Q'')(ds_2) Q(s_2) \Phi(s_2), \label{E5Dom}
\eea
where
\[
\Phi(s_2) = \int_{s_2}^T ds_1 \frac{(X_{s_1}-X_{s_1-s_2})^2}{Q(s_2)}.
\]
In Lemma \ref{LL1} below we will show that
\[
\Phi(s_2) \xrightarrow[s_2 \to 0+]{} T \textrm{ a.s.}
\]
so a.s. $s_2 \longmapsto \Phi(s_2)$ can be extended by continuity to $[0,T]$. If \eqref{EInfinite} holds,
then \eqref{E5Dom} is also infinite and so \eqref{E4Dom} is established. It remains to establish the following lemma.
\begin{lem} \label{LL1} Under the hypotheses of Proposition \ref{pDom2}, we have
\be
Z_{\eps}:=	\frac{1}{Q(\eps)} \int_{\eps}^{T} ds \left(X_s - X_{s-\eps}\right)^2 \xrightarrow[\eps \to 0]{} T \textrm{ a.s.}
\label{ELL1}
\ee
\end{lem}
\textbf{Proof}: 1) We have $E(Z_{\eps})= T-\eps$ and this obviously converges to $T$ when $\eps \to 0$. In order to prove that the convergence in \eqref{ELL1} holds in $L^2(\Om)$, it would be enough to show that
\[
Var(Z_{\eps}) \xrightarrow[\eps \to 0]{} 0.
\]
2) In order to prove the a.s. convergence we will implement the program of \cite{GN}, see in particular Lemma 3.1. This will only be possible because of technical assumption \eqref{tech}. We will show that
\be
\label{ELL2} Var(Z_{\eps}) = O\left(\frac{\eps}{Q(\eps)}\right).
\ee
Consequently 
\[
Var(Z_{\eps})=O(t^{\alpha}),
\]
$\a = 1-\a_1$ and (3.1) in \cite{GN} is verified. The upper bound of \eqref{tech} allows to show that $X$ is H\"older continuous, by use of Kolmogorov lemma.
\newline 3) We prove finally \eqref{ELL2}. We remark that $Q(\eps) \neq 0$ for $\eps$ in a neighbourhood of zero, otherwise  \eqref{E4Dom} cannot be true. We have
\[
Var(Z_{\eps}) = \frac{1}{Q(\eps)^2} \int_{\eps}^T ds_1 \int_{\eps}^T ds_2 Cov\left((X_{s_1}-X_{s_1-\eps})^2,(X_{s_2}-X_{s_2-\eps})^2 \right).
\]
It is well-known that given two mean-zero Gaussian random variables $\xi$ and $\eta$
\[
Cov(\xi^2, \eta^2) = 3 Cov(\xi,\eta)^2.
\]
This, together with the stationary increments property, implies that
\[
Var(Z_{\eps}) = \frac{6}{Q^2(\eps)} \int_{\eps}^T ds_1 \int_0^{s_1-\eps} ds_2 \left[ Cov\left(X_{s_2+\eps}-X_{s_2},X_{\eps} \right)\right]^2.
\]
Since, by Cauchy-Schwarz
\[
Cov \left( X_{s_2+ \eps}-X_{s_2},X_{\eps}\right) \leq Q(\eps)
\]
then
\[
Var(Z_{\eps}) = \frac{3}{Q^2(\eps)} I(\eps) + O(\eps),
\]
where
\[
I(\eps) = \int_0^T ds_1 \int_{\eps}^{s_1} ds_2\left( -Cov\left(X_{s_2+\eps} -X_{s_2}, X_{\eps} \right) \right)^2.
\]
Since Hypothesis \ref{hD} holds, Assumption (D) is verified and 
\[
-Cov \left(X_{s_2+\eps}-X_{s_2}, X_{\eps} \right) \geq 0.
\]
Hence
\bea
I(\eps) &\leq& Q(\eps) \int_0^T ds_1 \int_{\eps}^{s_1} ds_2 \left(2 Q(s_2) - Q(s_2+\eps) -Q(s_2-\eps) \right) \nn \\
 &\leq& Q(\eps) \int_0^T ds_1 \int_{\eps}^{s_1} ds_2 \left(\int_{s_2-\eps}^{s_2} Q'(y)dy - \int_{s_2}^{s_2+\eps}Q'(y)dy \right). \nn
\eea
Using Fubini's theorem, we obtain
\bea
&&\frac{I(\eps)}{Q^2(\eps)} = \frac{1}{Q(\eps)} \int_0^T ds_1 \left( \int_0^{s_1} dy Q'(y) \int_{\eps \vee y}^{(y+\eps)\wedge s_1} ds_2 - \int_{\eps}^{s_1+\eps} dy Q'(y) \int_{(y-\eps)\vee \eps}^{y\wedge s_1} ds_2 \right) \nn \\
&&= \frac{1}{Q(\eps)} \left( \int_0^T ds_1 \left\{ \int_0^{\eps} dy Q'(y)y + \eps \int_{\eps}^{s_1-\eps} dy Q'(y) + \int_{s_1-\eps}^{s_1} dy Q'(y) (s_1-y)\right\} \right. \nn \\
&&- \frac{1}{Q(\eps)} \left( \int_0^T ds_1\left\{ \int_{\eps}^{2 \eps} dy Q'(y)(y-\eps)+ \eps \int_{2 \eps}^{s_1}dy Q'(y) + \int_{s_1}^{s_1+\eps} dy Q'(y)(s_1-y +\eps) \right\} \right). \nn
\eea
Performing carefully the calculations, in particular commuting $ds_1$ and $dy$ through Fubini's, it is possible to show that
\[
\frac{I(\eps)}{Q^2(\eps)} \leq O(\eps) + O\left(  \frac{\eps}{Q(\eps)} \right).
\]
Assumption \eqref{tech} allows to conclude. \qed

\subsection{Paley-Wiener integral and integrals via regularization} \label{s5.2}

%%  By default %Let us suppose that Assumptions (A), (B) are verified. 
%The map $I:C_0^1(\R) \subset L_R \lra L^2(\Omega)$ defined as $g \lra I(g)$ admits a linear continuous extension to $L_R$ with respect to $\| \cdot\|_R$. It will still be denoted by $I$.

We start introducing the definition of Paley-Wiener integral. 
\begin{prop} \label{cd3.9} Let $g \in L_R$, then
\be
E \left( \iin g dX\right)^2 = \|g \|^2_{\cH}. \label{E3633}
\ee
Therefore the map $g \lra \iin g dX$ is continuous with respect to $\| \cdot\|_{\cH}$.
\end{prop}
\textbf{Proof}: \eqref{E3633} follows from Remark \ref{r3.1}. The second part of the statement follows because 
\[
\|g \|_{\cH} \leq \|g \|_R.
\]  \qed
\newline At this point the map $I: C^1_0(\R_+) \subset L_R \lra \LOm$ defined as $g \lra I(g)$ admits a linear continuous extension to $L_R$. It will still be denoted by $I$. 
\begin{df}\label{d3.9}
Let $g \in L_R$. We define \textbf{ the Paley-Wiener integral} of $g$ with respect to $X$ denoted by $\iin gdX$ the random variable $I(g)$.
\end{df}
\begin{prop} \label{p3.91} Under Assumptions (C), (D), if $\f$ has bounded variation with compact support, then
\be
\label{E3.91} \iin \f dX = -\iin X d\f.
\ee
In particular
\[
\iin 1_{[0,t]} dX = X_t.
\]
\end{prop}
\textbf{Proof}: By definition, \eqref{E3.91} holds for $g \in C^1_0$. We introduce the same sequence $(\f_n)$ as in the proof of Proposition \ref{p3.71}. By \eqref{E3633}
\[
E\left( \iin (\f_n -\f)dX\right)^2 = \|\f_n -\f \|^2_R.
\]
This converges to zero when $n \lra \infty$ as it was shown in the proof of Proposition \ref{p3.71}. In the same proof it was established that $\iin X d\f_n \lra \iin X d\f$ in $\LOm$.  \qed 
\newline We recall briefly the notion of integrals via regularization in the spirit of \cite{RV1993} or \cite{RVsem}. We propose here a definite type integral.
\begin{df} \label{DEF} Let $Y$ be a process with paths in $L^1_{loc}(\R)$. We say that the \textbf{forward} (resp. \textbf{backward, symmetric}) \textbf{integral} of $Y$ with respect to $X$
exists, if the following conditions hold.
\begin{description}
	\item[a)] For $\eps > 0$ small enough the following Lebesgue integral
\bea
 I(\eps,Y,dX) = \iin&&Y_s\frac{X_{s+\eps}-X_s}{\eps} ds\nn \\
\textrm{(resp. }\iin&&Y_s\frac{X_{s}-X_{s-\eps}}{\eps} ds\nn \\
\iin&&Y_s\frac{X_{s+\eps}-X_{s-\eps}}{\eps} ds)\nn
\eea
with the usual condition $X_s=0$, $s \leq 0$ exists.
	\item[b)] $\lim_{\eps \to 0} I(\eps,Y, dX)$ exists in probability.
\end{description}
The limit above will be denoted by
\[
\iin Yd^-X \textrm{ (resp. } \iin Yd^+X, \iin Yd^oX).
\]
\end{df}

\begin{prop} \label{p3.11}
Let $f: \R_+ \lra \R$ cadlag bounded. Suppose the existence of $V_f: \R_+ \lra \R_+$ such that
\begin{description}
	\item[i)] 
	\[
	|f(s_2)-f(s_1)| \leq V_f(s_2-s_1), \ s_1,s_2 \geq 0,
	\]
	\item[ii)]
\be \label{3.11}
\int_{\RR} V_f^2(s_2-s_1) d|\mu|(s_1,s_2) < \infty
\ee
\end{description}
Then
\[
\iin f d^{*}X = \iin fdX, \ \ \ * \in \{-,+,0\}.
\]
In particular $\iin f d^* X$ exists.
\end{prop}
\textbf{Proof}: We consider the case $* = -$, the other cases being similar. The quantity
\[
\iin f(s) \frac{X_{s + \eps}-X_s}{\eps} ds
\]
equals
\[
\iin f_{\eps}(u)dX_u,
\]
where
\[
f_{\eps}(u) = \frac{1}{\eps} \int_{u-\eps}^u f(s)ds = \frac{1}{\eps} \int_{-\eps}^0 f(s+u)ds,
\]
with the convention that $f$ is prolongated by zero on $\R_{-}$.
\par It remains to show that $f_{\eps} \lra f$ in $L_R$. By Lebesgue's dominated convergence theorem and the fact that $f$ is bounded, cadlag and $|R|(ds,\infty)$ is non-atomic, we have
\[
\iin(f_{\eps}-f)^2(s) d|R|(ds,\infty) \xrightarrow[\eps \to 0] {} 0.
\]
It remains to show that
\[
\lim_{\eps \to 0} \int_{\RR} d|\mu|(s_1,s_2) ((f_{\eps}-f)(s_1)-(f_{\eps}-f)(s_2))^2 = 0.
\]
Indeed
\bea
&&|(f_{\eps}-f)(s_1)-(f_{\eps}-f)(s_2)| = \frac{1}{\eps} \left| \int_{-\eps}^{0} \left[(f(y + s_1)-f(s_1))-(f(y+s_2)-f(s_2)) \right]dy \right| \nn \\
 &&\ \ \ = \frac{1}{\eps} \left| \int_{-\eps}^{0} \left[(f(y + s_1)-f(y+s_2))-(f(s_1)-f(s_2)) \right]dy \right| \nn\\
 && \ \ \ \ \ \leq  \frac{1}{\eps}  \int_{-\eps}^{0} \left|(f(y + s_1)-f(y+s_2)\right|dy-|f(s_1)-f(s_2)| dy \leq 2 V(s_2-s_1).\nn
\eea
Since $f_{\eps} \lra f$, (\ref{3.11}) and Lebesgue's dominated convergence theorem allow to conclude. \qed

\subsection{About some second order Paley-Wiener integral} \label{s8}

We introduce now a second order Wiener integral of the type:
\[
I_2(g):=\int_{\RR} g(s_1,s_2)dX_{s_1}^1 dX_{s_2}^2,
\]
where $g:\RR \lra \R$ is a suitable function and $X^1,X^2$ are two independent copies of X.

In fact all the considerations can be extended to Wiener integrals with respect to $n$ copies $X^1, \ldots, X^n$ of $X$, but in order not to introduce technical complications we only consider the case $n=2$. This case will be helpful in section \ref{s8bis} in order to topologize the tensor product $L_R \otimes L_R$.

We will make use of tensor product spaces in the Hilbert framework. For a complete information about tensor product spaces and topologies the reader can consult \cite{Ryan}. We suppose here the validity of Assumption (C). We denote $\n=|R|(dt,\infty)$ as before.
If $g = g_1\otimes g_2$, $g_1, g_2 \in L_R$ then we set
\be
\label{EI4.1} I_2(g) = \iin g_1 dX_1 \iin g_2dX_2.
\ee
We remark that $g(s_1,s_2)= g_1(s_1)g_2(s_2)$.
We denote by $L_R \otimes L_R$ the algebraic tensor product space of linear combinations of functions of the type $g_1 \otimes g_2$, $g_1,g_2 \in L_R$.

We define $\tilde{L}_{2,R}$ as the space of Borel functions $g: \RR \lra \R$ such that
\be
\| g\|^2_{2,R} = \iin \n(dt)\| g(t,\cdot)\|^2_R + \frac{1}{2} \int_{\RR} d|\m|(s_1,s_2) \|g(s_1,\cdot)-g(s_2,\cdot) \|^2_R < \infty. \label{EI4.111}
\ee
An easy property which can be established by inspection is given below. 
\begin{lem} \label{lI4.1}
For $g: \RR \lra \R$ we have 
\bea
\| g\|^2_{2,R} &=& \int_{\RR} \n(ds_1)\n(ds_2) g^2(s_1,s_2)\nn \\
 &&+ \frac{1}{4} \int_{\RR \times \RR}(g(s_1,t_1)-g(s_2,t_1)-g(s_1,t_2)+g(s_2,t_2))^2 d|\m|(t_1,t_2)d|\m|(s_1,s_2) \nn\\
&& + \fr \iin d\n(s)\left\{\int_{\RR}(g(s,t_1)-g(s,t_2))^2 d|\m|(t_1,t_2) + \int_{\RR}(g(s_1,s)-g(s_2,s))^2 d|\m|(s_1,s_2) \right\} \nn\\
&=& \iin \n(dt)\| g(\cdot,t)\|^2_R + \fr \int_{\RR} d|\m|(t_1,t_2) \| g(\cdot,t_1)-g(\cdot,t_2)\|^2_{\RR}.\nn
\eea
\end{lem}

\begin{rem} \label{rI71}
We observe that second term of the right-hand side equals
\[
\frac{1}{4} \int_{\R^2_+ \times \R^2_+} \left(\Delta_{]s_1,s_2] \times ]t_1,t_2]}g \right)^2 d|\m|(t_1,t_2)d|\m|(s_1,s_2),
\]
where $\Delta_{]s_1,s_2] \times ]t_1,t_2]}g$ is the planar increment introduced in 
Section \ref{s0}.
\end{rem}
\begin{rem} \label{rI4.2}
\begin{enumerate}
	\item The (semi)-norm $\| \cdot\|_{2,R}$ derives from an inner product. We have
	\bea
	\left<f,g\right>_{2,R} =&& \iin \n(ds) \left\langle f(s,\cdot), g(s,\cdot)\right\rangle_{R} \nn \\ 
	&&+ \fr \int_{\RR} d|\m|(s_1,s_2) \left\langle f(s_1,\cdot)-f(s_2,\cdot),g(s_1,\cdot)-g(s_2,\cdot)  \right\rangle_{R} \nn
	\eea
	\item An analogous expression to Lemma's \ref{lI4.1} statement, can be written for $\left<\cdot,\cdot \right>_{2,R}$ instead of $\|\cdot\|_{2,R}$.
	\item If $f = f_1 \otimes f_2$, $f_1, f_2 \in L_R$ then $f \in \tilde{L}_{2,R}$. If $g= g_1 \otimes g_2$, $g_1,g_2 \in L_R$
	\[
	\left<f,g\right>_{2,R} = \left<f_1,g_1\right>_R \left<f_2,g_2\right>_R.
	\]
	\item $L_R \otimes L_R$ is included in $\tilde{L}_{2,R}$. In particular any linear combination of the type $\phi \otimes \phi$ belong to $L_R \otimes L_R$.
	\item Similarly to the proof of Proposition \ref{p3.5}, taking into account Assumption (C), it is possible to show that $\tilde{L}_{2,R}$ is complete and it is therefore a Hilbert space.
	\item Since $\left<\cdot,\cdot\right>_{2,R}$ is a scalar product and because of 3., it follows that $\|\cdot \|_{2,R}$ is the Hilbert tensor norm of $L_R \otimes L_R$. For more information about tensor topologies, see e.g. \cite{Ryan}. The closure of $L_R \otimes L_R$ with respect to $\|\cdot \|_{2,R}$ can be identified with the Hilbert tensor product space $L_R \otimes^h L_R$; it will be denoted by $L_{2,R}$. 
	\item $\tilde{L}_{2,R}$ can also be equipped with the scalar product $\left<\cdot, \cdot\right>_{2,\cH}$
	\bea
	\left<f,g\right>_{2,\cH} =&& \iin R(ds,\infty) \left<f(s,\cdot),g(s,\cdot)\right>_{\cH}\nn \\
	&& -\fr \int_{\RR} \m(ds_1,ds_2)\left<f(s_1,\cdot)-f(s_2,\cdot),g(s_1,\cdot)-g(s_2,\cdot)\right>_{\cH}.\nn
	\eea
\end{enumerate}
\end{rem}
\begin{rem} \label{rI4.3}
Similar considerations as in Remark \ref{rI4.2} can be made for the inner product $\left<\cdot, \cdot\right>_{2,\cH}$.
\begin{enumerate}
	\item If $f,g$ are as in 3. of Remark \ref{rI4.2}, then
	\[
	\left<f,g\right>_{2,\cH} = \left<f_1,g_1\right>_{\cH}\left<f_2,g_2\right>_{\cH}.
	\]
	\item We denote by $\| \cdot\|_{2,\cH}$ the associated norm. Analogous expressions as for Lemma \ref{lI4.1} can be found for $\|\cdot \|_{2,\cH}$.
	\item If Assumption (D) is fulfilled then $L_R$ can be identified with $\cH$ and $\left<\cdot,\cdot\right>_{2,\cH}$ coincides with $\left<\cdot,\cdot\right>_{2,R}$. The Hilbert tensor product $\cH \otimes ^h \cH$ can be identified with $L_{2,R}$.
	\item If $f \in \tilde{L}_{2,R}$ then
	\[
	\| f\|_{2,\cH} \leq \| f\|_{2,R}.
	\]
\end{enumerate}
\end{rem}
The double integral application $g \longmapsto I_2(g)$ extends by linearity through (\ref{EI4.1}) to the algebraic tensor product $L_R \otimes L_R$.
\begin{prop} \label{pI4.4}
$I_2:L_R \otimes L_R \lra L^2(\Om, \mathcal{F}, P)$ extends continuously to $L_{2,R}$. In particular for every $g\in L_{2,R}$ we have
\[
E\left(I_2(g)^2\right) = \| g\|^2_{2,R}.
\]
\end{prop}
\textbf{Proof}: Let $g \in L_{2,R}$, so $g = \sum_{i=1}^n g_{i1} \otimes g_{i2}$, $\ g_{i1},g_{i2}\in L_R$, then
\bea
E(I_2(g)^2) &=& \sum_{i,j=1}^n E(I_2(g_{i1}\otimes g_{i2})I_2(g_{j1}\otimes g_{j2}))\nn \\
&=&\sum_{i,j=1}^n E\left( \iin g_{i1}dX^1 \iin g_{i2} dX^2 \iin g_{j1}dX^1 \iin g_{j2} dX^2\right)\nn \\
&=& \sum_{i,j=1}^n E\left( \iin g_{i1}dX^1 \iin g_{j1} dX^1 \right)E \left( \iin g_{i2}dX^2 \iin g_{j2} dX^2\right) \nn
\eea
using the independence of $X^1$ and $X^2$. Therefore by Remark \ref{rI4.3} 1. and bilinearity of the inner product, it follows
\bea
E(I_2(g))^2  &=& \sum_{i,j=1}^n \left<g_{i1},g_{j1}\right>_{\cH}\left<g_{i2},g_{j2}\right>_{\cH}\nn \\
&=&\sum_{i,j=1}^n \left<g_{i1}\otimes g_{i2}, g_{j1} \otimes g_{j2}\right>_{2,\cH} = \| g\|^2_{2,\cH} \leq \| g\|^2_{2,R}. \nn
\eea
This allows to conclude the proof of the proposition. \qed
\begin{rem} \label{rI4.5} The proof of Proposition \ref{pI4.4} allows to establish (as a by product) that for $g \in L_{2,R}$ the double integral is unambiguously  defined.
\end{rem}
Given $g \in L_{2,R}$, we denote
\[
I_2(g)= \int_{\RR} g(s_1,s_2) dX^1_{s_1}dX^2_{s_2}.
\]
This quantity is called \textbf{double (Paley-)Wiener integral} of $g$ with respect to $X^1$ and $X^2$.

We will characterize now some significant functions which belongs to $L_{2,R}$.
\begin{lem} \label{PD14} Suppose that $1_{[0,t]} \in L_R$ for every $t > 0$ and Assumption  (C). Then for any $t_1,t_2 > 0$, $y_1,y_2>0$, $h = 1_{]t_1,t_2]\times ]y_1,y_2]}$ belongs to $L_{2,R}$.
\end{lem}
\textbf{Proof:} Since $1_{]t_1,t_2]}, 1_{]y_1,y_2]} \in L_R$, clearly $h \in L_R \otimes L_R \subset L_{2,R}$.\qed

%Before giving the proof of Proposition \ref{PD15}, we provide a useful, intermediate result.
%\begin{lem} \label{PD16} Let $a < b \leq c < d$, $\m(]a,b]\times]c,d]) = E((X_b-X_a)(X_d-X_c))$.
%\end{lem}
%\textbf{Proof:} Let $\eps > 0$. Since $]a,b]\times ]c+\eps,d]$ does not intersect the diagonal, it is clear that
%\[
%E((X_b-X_a)(X_d-X_{c+\eps})) = \m(]a,b]\times]c+\eps,d])
%\]
%The left-hand side converges to
%\[
%E((X_b-X_a)(X_d-X_c))
%\]
%since $X$ is continuous in $L^2$. The right-hand side goes to $\m(]a,b]\times]c,d])$ because $(-\m)$ is a non-negative measure.\qed
%We introduce below the notion of bounded variation functions $g: \RR \lra \R$ which vanish outside the first quater of the plane.
%\begin{df} \label{DPBV} Let $g: \RR \lra \R$. We say that it fulfills the \textbf{PBV property} if
%\begin{enumerate}
%	\item $g$ has bounded planar variation,
%	\item $g(t_1,0)$, $g(0,t_2)$ have bounded variation for any $t_1, t_2 >0$.
%\end{enumerate}
%\end{df}
%\begin{rem} \label{RD150} A consequence of \eqref{EPBV} is the following. $g: \RR \lra \R$ fulfills the PBV property, if and only if  there is a finite signed measure  $\chi$ on $\RR$ 
%\[
%g(t_1,t_2) = \chi(]0,t_1]\times]0,t_2]).
%\] 
%\end{rem}
\begin{rem} \label{RDPBV}
1) If $g$ is a sum of functions of the type $g^1\otimes g^2$, where $g^1,g^2: \R_+\lra \R$ are bounded variation functions with compact support, then $g$ has bounded planar variation.
\newline 2) If $g$ is as in item 1) and $X^1,X^2$ are independent copies of $X$, then
\bea
\int_{\RR} g(t_1,t_2) dX^1_{t_1} dX^2_{t_2} &=& \int_{]0,\infty[^2} X^1_{t_1}X^2_{t_2}dg(t_1,t_2),\nn 
\eea
where the right-hand side is a Lebesgue integral with respect to the signed measure $\chi$ such that
\[
g(t_1,t_2) = \chi(]0,t_1]\times]0,t_2]).
\]
This follows because of the following reasons.
\par If $g = g_1 \otimes g_2$, $g_1$, $g_2$ have bounded variation with compact support then
\begin{description}
	\item[i)] $I_2(g)= \iin g_1 dX^1 \iin g_2 dX^2$,
	\item[ii)] $\int_{]0,\infty[^2} \f dg = \int_{]0,\infty[^2} \f(s_1,s_2) dg_1(s_1) dg_2(s_2)$,
	\item[iii)] $\iin g_1 dX^1 = - \iin X_s^1 dg_1(s)$, because of Proposition \ref{p3.91}.
\end{description}
\end{rem}
A significant proposition characterizing elements of $L_{2,R}$ under Assumption (D) is the following. 
\begin{prop} \label{PD150} We suppose, that Assumptions (C) and (D) are verified. Moreover we suppose that $h: \RR \lra \R$ has bounded planar variation. Then $h \in L_{2,R}$.
\end{prop}
\begin{rem} \label{RD151} If $X_t = X_T$, $t \geq T$ then the statement of Proposition \ref{PD150} holds if $h|_{[0,T]^2}$ has bounded planar variation.
\par Indeed, by definition of $L_{2,R}$, if $h$ is prolongated by zero outside $[0,T]^2$ denoted by $\bar{h}$, then $\|h - \bar{h} \|_{2,R} =0$. In particular $h = \bar{h}  \  \nu_\infty \otimes \nu_\infty$ a.e. since
\[
\|h \|_{L^2(\nu_\infty) \otimes L^2(\nu_\infty)} \leq \|h - \bar{h} \|_{2,R},
\]
where $\nu_\infty = R(ds,\infty).$
\end{rem}
%\begin{rem} \label{RLD151} $h \in C^2_0(\RR)$ fulfills the assumptions of Proposition \ref{PD150} with planar quadratic variation
%\be
%\label{E152} \|h \|_{pv} := \int_{]0,\infty[^2} \left| \frac{\pa^2 h}{\pa t_1 \pa t_2}(t_1,t_2)\right| dt_1 dt_2
%\ee
%and
%\bea
%\|h(t_1,\cdot) \|_{var}&=& \iin \left|\frac{\pa h}{\pa t_2} (t_1,0)\right| dt_1, \nn \\
%\|h(\cdot, t_2)\|_{var}& = & \iin \left|\frac{\pa h}{\pa t_1}(0,t_2) \right| dt_2. \nn
%\eea
%\end{rem}
\textbf{Proof} of Proposition \ref{PD150}: 
Let $N > 0$ and $t_i:=t_i^N := \frac{i}{N}$, $0 \leq i \leq N^2$. According to Corollary \ref{c3.72} $1_{]t_i, \times t_{i+1}]}, 1_{]t_j,t_{j+1}]}$ belongs to $L_R$ for any $0 \leq i,j \leq N^2$.
We denote 
\bea
h^N (s_1,s_2) &=& \sum_{i,j=0}^N h(t_i,t_j)1_{]t_i,t_{i+1}]}(s_1) 1_{]t_j,t_{j+1}]} (s_2)\nn.
\eea
Of course $h^N$ belongs to $L_R \otimes L_R \subset L_{2,R}$.
\par On the other hand $h^N \lra h$ for every continuity point. The total variation of $dh^N$ is bounded by
\bea
&&\sum_{i,j =0}^{N^2} 1_{]t_i,t_{j+1}]} \otimes 1_{]t_j,t_{j+1}]} \left|\Delta h_{]t_{i},t_{i+1}] \times ]t_j, t_{j+1}]} \right| \nn
 %+ \sum_{i=0}^{N^2} 1_{]t_i,t_{i+1}]} |h(t_{i+1},0) - h(t_i,0)| \nn \\
%&& + \sum_{j=0}^{N^2} 1_{]t_j,t_{j+1}]} |h(0,t_{j+1})-h(0,t_j)| + |h|(0,0). \nn
\eea
Previous quantity is bounded by $\|h \|_{pv}$.
%\[
%\|h \|_{pv} + \|h(0,\cdot) \|_{var} + \|h(\cdot,0) \|_{var} + |h(0,0)|.
%\]
Finally $h^N$ converges weakly to $h$, by the theory of two-parameter distribution functions of measures. Therefore, if $X^1$ and $X^2$ are two independent copies of $X$, then 
\be
\label{F111} \int_{]0,\infty[^2} X_{t_1}X_{t_2} dh^N(t_1,t_2) \xrightarrow[N \to \infty] {} \int_{]0,\infty[^2} X_{t_1}X_{t_2} dh(t_1,t_2) \  a.s.
\ee
%In order to reduce complexity of notations, but without restriction of generality, we suppose the sequence $h^N$ and $h$ vanishing on the axes.
\par By Remark \ref{RDPBV} 2), we have
\be \label{E300} \int_{\RR} \left( h^N -h^M\right) (t_1,t_2) dX_{t_1}^1 dX_{t_2}^2 = \int_{]0,\infty[^2} X^1_{t_1}X^2_{t_2} d\left(h^N -h^M\right) (t_1,t_2).
\ee
%Previous considerations imply that this converges to zero in $\LOm$. 
By Fubini's, the fact that $X^1$ and $X^2$ are independent and \eqref{E300}, it follows
\bea
\|h^N -h^M \|^2_{2,R} &=& E \left( \int_{]0,\infty[^2} (h^N-h^M)(t_1,t_2) dX^1_{t_1} dX^2_{t_2} \right)^2 \nn \\
&=& \int_{]0,\infty[^4} R(t_1,s_1)R(t_2,s_2) d(h^N-h^M)(t_1,t_2) d(h^N-h^M)(s_1,s_2). \nn
%&=& \left(\int_{\RR} R(t_1,t_2) d(h^N-h^M) (t_1,t_2) \right) \xrightarrow[N,M \to \infty]{} 0, \nn
\eea
This converges to zero because $dh^N\otimes dh^M$ weakly converges when $N,M \to \infty$ and $(t_1,s_1,t_2,s_2) \mapsto R(t_1,s_1)R(t_2,s_2)$ is a continuous function. Consequently the sequence $(h^N)$ is Cauchy in $L_{2,R}$.
\par Since $L_{2,R}$ is complete, there is $\psi: \RR \lra \R \in L_{2,R}$ such that 
\[
\|h^N -\psi \|_{2,R} \xrightarrow[N\to \infty]{}0.
\]
By definition of $\|\cdot \|_{2,R}$, we have
\be
\label{F113} \|h^N - \psi \|^2_{L^2(d\n_{\infty})^{\otimes 2}} \leq \|h^N- \psi \|^2_{2,R} \lra 0,
\ee
where again
\[
\n_{\infty} = R(ds,\infty).
\]
So there is a subsequence $(N_k)$ such that $\|h^{N_k} - \psi \|\lra 0$ $ \n_{\infty}\otimes \n_{\infty}$ a.e. 
\par Since $h^N \lra h$ excepted on a countable quantity of points and $\n_{\infty} \otimes \n_{\infty}$ is non-atomic, then $h^N \lra h$ $\n_{\infty} \otimes \n_{\infty}$ a.e. Finally $h= \psi$ $\n_{\infty} \otimes \n_{\infty}$ a.e. and therefore
\[
\| h -h^N\|_{2,R} = 0
\]
and so $h \in L_{2,R}$. \qed
\par A side-effect of the proof of Proposition \ref{PD150} is the following.
\begin{prop} \label{R1000} If $h: \RR \lra \R$ has bounded planar variation, then
\bea
\int_{\RR} h(s_1,s_2) dX^1_{s_1} dX^2_{s_2} &=& \int_{]0,\infty[^2} X^1_{s_1} X^2_{s_2} d\chi(s_1,s_2), \nn
\eea
where as usual $\chi(]0,s_1] \times]0,s_2]) = \Delta_{]0,s_1] \times ]0,s_2]}h$.
\end{prop}
\begin{rem} \label{R3648}
From Proposition \ref{pI4.4} and Proposition \ref{R1000} we obtain
\bea
\| h\|^2_{2,R} &=&E \left(\int_{\RR} h(s_1,s_2)dX^1_{s_1}dX^2_{s_2} \right)^2 \nn\\ 
 \label{ER3648} \\
 &=& \int_{]0,\infty[^4} R(t_1,s_1)R(t_2,s_2) dh(t_1,t_2)dh(s_1,s_2). \nn
\eea
\end{rem}
Another consequence of Proposition \ref{PD150} is the following.
\begin{prop} \label{PD160} We suppose the following.
\begin{description}
	\item[a)] Assumptions (C) and (D);
	\item[b)] there is $r_0 > 0$ such that 
	\be
	\sup_{r \in ]0,r_0]} \int_{\RR} d|\m| (t_1,t_2) Var(X_{t_1+r}-X_{t_2+r}) < \infty \label{ELTR},
	\ee
	\item[c)] $X_t = X_T$, $t \geq T$.
\end{description}
Then, for $a \in \R$ small enough, $h(s,t )=1_{[0,(t+a)_+]}(s) \in L_{2,R}$. Moreover
\be
 \|h \|^2_{2,R} = \int_{\RR} d(-\m)(t_1,t_2) Var\left(X_{(t_2+a)_+}-X_{(t_1+a)_+} \right)+ \iin R(ds,\infty)Var(X_{(s+a)+}). \nn
\ee
\end{prop}
\textbf{Proof}: We regularize the function $h$ in $t_1$. Let $\r$ be a smooth function with compact support on $\R_+$, $\r_{\eps}(x) = \frac{1}{\eps} \r \left(\frac{x}{\eps} \right)$, $x \in \R$. We set
\bea
F^{\eps}(s,t) &=& \int_{\R} \r_{\eps}(s-s_1) 1_{[0,(t+a)_+]} (s_1)ds_1 \nn \\
\label{EEPD160}
\\
&=& \int_{\frac{s-(t+a)_+}{\eps}}^{\frac{s}{\eps}} \r(\tilde{s}_1)d\tilde{s}_1 + \iin R(dt,\infty) Var(X_t). \nn
\eea
Of course we have 
\be
\label{EEPD161} \frac{\pa F^{\eps}}{\pa s}(s,t) = \r \left(\frac{s}{\eps}\right) - \r \left(\frac{s-(t+a)_+}{\eps} \right).
\ee
Since $F^{\eps}$ is smooth, by Remark \ref{RD151} and Proposition \ref{PD150}, it follows that $F^{\eps} \in L_{2,R}$. 
\par It remains to show that $F^{\eps} \lra h$ in $L_{2,R}$. First of all, we observe that $F^{\eps} \lra h$ pointwise. We need to show that $\|F^{\eps}-h \|_{2,R} \lra 0$ when $\eps \to 0$. We have 
\[
\|F^{\eps}-h \|_{2,R} = I_1(\eps) + I_2(\eps),
\]
where
\bea
I_1(\eps) &=& \iin R(dt,\infty) \|F^{\eps}(\cdot,t)-1_{[0,(t+a)_+]} \|^2_R \nn \\
I_2(\eps) &=& \int_{\RR} d|\m|(t_1,t_2) \left\|F^{\eps}(\cdot,t_1)-F^{\eps}(\cdot,t_2) -\left(1_{[0,(t_1+a)_+]}-1_{[0,(t_2+a)_+]} \right)\right\|^2_{R}. \nn
\eea
We will first evaluate
\be
\label{Lepst} \|F^{\eps}(\cdot, t) \|^2_R,
\ee
\be
\label{Lepst12} \|F^{\eps}(\cdot,t_1) - F^{\eps}(\cdot,t_2) \|^2_R.
\ee
Now
\bea
\|F^{\eps}(\cdot,t) \|^2_R &=& E \left(\iin F^{\eps}(s,t) dX_s \right)^2 = E \left(\iin X_s \frac{\pa F^{\eps}}{\pa s}(s,t) ds \right)^2 \nn \\
&=& \frac{1}{\eps^2} \int_{\RR} ds_1 ds_2 R(s_1,s_2) \left[\r\left(\frac{s_1}{\eps}\right) - \r\left(\frac{s_1 -(t+a)_+}{\eps} \right) \right] \nn \\
&& \ \ \ \ \ \ \ \ \ \ \ \ \ \ \ \ \left[ \r \left(\frac{s_2}{\eps}\right) - \r\left( \frac{s_2 - (t+a)_+}{\eps} \right) \right] \nn \\
&=& I_{++}(\eps) - I_{+-} (\eps,t) -I_{-+}(\eps,t) + I_{--}(\eps,t), \nn
\eea
where after an easy change of variable, 
one can easily see that
\[
\sup_{t \geq 0} |I_{++}(\eps)+ I_{+-}(\eps,t)+ I_{-+}(\eps,t) | \xrightarrow[\eps \to 0]{} 0
\]
because $R(0,s_2)= R(s_1,0) \equiv 0$, $\forall s_1,s_2 >0$. On the other hand
\bea
I_{--}(\eps,t) &=& \frac{1}{\eps^2} \int_{\RR} ds_1 ds_2 R(s_1,s_2) \r \left(\frac{s_1-(t+a)_+}{\eps} \right) \r\left(\frac{s_2-(t+a)_+}{\eps} \right) \nn \\
&=& \int_{\RR} d\tilde{s}_1 d\tilde{s}_2 R\left((t+a)_+ +\eps \tilde{s}_1, (t+\eps)_+ + \eps \tilde{s}_2 \right) \r(\tilde{s}_1) \r(\tilde{s}_2). \nn
\eea
By Lebesgue's dominated convergence theorem,
\bea
I_1(\eps) &=& \int_{\RR} d\r(s_1)d\r(s_2) \left[R\left((t+a)_+ +\eps s_1, (t+a)_+ + \eps s_2\right)-R((t+a)_+,(t+a)_+) \right]\frac{R(dt,\infty)}{2} \nn\\
&+& J(\eps), \nn
\eea
where $\lim_{\eps \to 0} J(\eps) = 0$,
so
\[
\iin \|F^{\eps}(\cdot,t) \|^2_R R(dt,\infty) \xrightarrow[\eps \to 0] {} \iin \|1_{[0,(t+a)_+]} \|^2_R R(dt,\infty) = \iin Var(X_{(t+a)+})R(dt,\infty).
\]
Similarly, we can show that
\[
\iin \left<F^{\eps}(\cdot,t), 1_{[0,(t+a)_+]} \right>_R R(dt, \infty) \xrightarrow[\eps \to 0]{} \iin \|1_{[0,(t+a)_+]} \|^2_R R(dt,\infty).
\]
This implies that $\lim_{\eps \to 0} I_1(\eps)= 0$. Concerning $I_2(\eps)$, we need to evaluate \eqref{Lepst12}. We observe that
\bea
\|F^{\eps}(\cdot,t_1) - F^{\eps}(\cdot,t_2) \|^2_R &=& E\left( \iin X_s \left(\frac{\pa F^{\eps}}{\pa s}(s,t_1) - \frac{\pa F^{\eps}}{ \pa s} (s,t_2)\right)ds \right)^2 \nn \\
&=& E\left( \iin X_s \left[ \r\left(\frac{s-(t_2+a)_+}{\eps} \right)- \r \left( \frac{s-(t_1+a)_+}{\eps} \right)\right] ds\right)^2 \nn \\
&=& K_{++}(\eps,t_1,t_2) -K_{-+}(\eps,t_1,t_2) -K_{+-}(\eps,t_1,t_2) + K_{--}(\eps,t_1,t_2), \nn
\eea
where
\bea
K_{++}(\eps,t_1,t_2) &=& \int_{\RR} ds_1 ds_2 R(s_1,s_2) \r \left(\frac{s_1 -(t_2+a)_+}{\eps} \right) \r\left( \frac{s_2 -(t_2+a)_+}{\eps} \right), \nn \\
K_{+-}(\eps,t_1,t_2) &=& \int_{\RR} ds_1 ds_2 R(s_1,s_2) \r \left(\frac{s_1 -(t_2+a)_+}{\eps} \right) \r \left( \frac{s_2 -(t_1+a)_+}{\eps} \right), \nn \\
K_{-+}(\eps,t_1,t_2) &=& \int_{\RR} ds_1 ds_2 R(s_1,s_2) \r \left(\frac{s_1 -(t_1+a)_+}{\eps} \right) \r \left( \frac{s_2 -(t_2+a)_+}{\eps} \right), \nn \\
K_{--}(\eps,t_1,t_2) &=& \int_{\RR} ds_1 ds_2 R(s_1,s_2) \r \left(\frac{s_1 -(t_1+a)_+}{\eps} \right) \r \left( \frac{s_2 -(t_1+a)_+}{\eps} \right). \nn 
\eea
Consequently
\bea
K_{++}(\eps,t_1,t_2) &=& \int_{\RR} ds_1 ds_2 \r(s_1)\r(s_2)R((t_2+a)_+ + \eps s_1, (t_2+a)_+ + \eps s_2), \nn \\
K_{+-}(\eps,t_1,t_2) &=& \int_{\RR} ds_1 ds_2 \r(s_1)\r(s_2)R((t_2+a)_+ + \eps s_1, (t_1+a)_+ + \eps s_2), \nn \\
K_{-+}(\eps,t_1,t_2) &=& \int_{\RR} ds_1 ds_2 \r(s_1)\r(s_2)R((t_1+a)_+ + \eps s_1, (t_2+a)_+ + \eps s_2), \nn \\
K_{--}(\eps,t_1,t_2) &=& \int_{\RR} ds_1 ds_2 \r(s_1)\r(s_2)R((t_1+a)_+ + \eps s_1, (t_1+a)_+ + \eps s_2). \nn
\eea
Hence
\bea
&& \int_{\RR} d|\m| (t_1,t_2) \|F^{\eps}(\cdot,t_1)-F^{\eps}(\cdot,t_2) \|^2_R \\
&& = \int_{\RR} d\r(s_1)d\r(s_2) \int_{\RR} d|\m|(t_1,t_2) \nn \\
&& \ \ \ \ \ Cov\left( X_{(t_2+a)_+ + \eps s_1}-X_{(t_1 +a)_+ + \eps s_1}, X_{(t_2+a)_+ + \eps s_2} - X_{(t_1 + a)_+ + \eps s_2} \right).
\eea
By Fubini's, Cauchy-Schwarz,  choosing the support of $\r$ small enough, and taking into account hypothesis b) of the statement, previous expression equals
\bea
&&\int_{\RR} d|\m|(t_1,t_2) \iin d\r(s) Var \left(X_{(t_2+a)_+ + \eps s}- X_{(t_1+a)_+ + \eps s} \right) \nn \\
&& = \iin d\r(s) \int_{\RR} d|\m|(t_1,t_2) Var\left (X_{(t_2+a)_+ + \eps s}-X_{(t_1+a)_+ + \eps s} \right). \nn
\eea
Lebesgue's dominated convergence theorem says that previous expression goes to
\bea
&&\int_{\RR} d|\m|(t_1,t_2) Var\left(X_{(t_2+a)_+}-X_{(t_1+a)_+} \right) \nn \\
\label{H190} \\
&& = \int_{\RR} d|\m|(t_1,t_2) \left\|1_{[0,(t_2+a)_+]}-1_{[(t_1+a)_+]} \right\|^2_R.\nn
\eea
This shows that
\[
\lim_{\eps \to 0} \int_{\RR} d|\m|(t_1,t_2) \|F^{\eps}(\cdot,t_1) -F^{\eps}(\cdot, t_2) \|^2_R
\]
equals the right-hand side of \eqref{H190}.
By similar arguments we can show that
\[
\int_{\RR} d|\m|(t_1,t_2) \left<F^{\eps}(\cdot,t_1)-F^{\eps}(\cdot,t_2), 1_{[0,(t_1+a)_+]}-1_{[0,(t_2+a)_+]} \right>_R
\]
converges again to the right-hand side of \eqref{H190}. This finally shows $\lim_{\eps \to 0} I_2(\eps) = 0$ and the final result. \qed
\par An interesting consequence is the following.
\begin{prop} \label{PRD17}
We suppose Assumptions (C), (D). Let $g \in L_{2,R}$. Suppose that $X$ fulfills the following assumption
\be
\label{EPD15} \int_{\RR} Var(X_{t_1}-X_{t_2})d|\m|(t_1,t_2) < \infty.
\ee
Then $g(s, \cdot) \in L_R$, $ R(ds,\infty)$ a.e., $s \longmapsto \iin g(s,t)dX_t \in L_R$ a.s. and it belongs to $L^2(\Om;L_R)$. Moreover
\be \int_{\RR} g(s,t) dX^1_s dX_t = \iin \left(\iin g(s,t)dX_t\right) dX^1_s. \label{EM11}
\ee
if $X^1$ is an independent copy distributed as $X$.
\end{prop}
\begin{cor} \label{CRD18} We suppose Assumptions (C), (D), \eqref{ELTR} and $X_t = X_T$ for $t \geq T$. 
\newline 1) We have $s \longmapsto X_s \in L_R$ a.s. and it belongs to $L^2(\Om;L_R)$.
\newline 2) Let $X^1$ be an independent copy of $X$. If $h(s_1,s_2) = 1_{[0,s_1 \wedge T]}(s_2)$, then
\begin{description}
	\item[i)] 
	\[
	\int_{\RR} h(s_1,s_2) dX_{s_1}^1 dX_{s_2} = \iin X_s dX_s^1.
	\]
	\item[ii)] 
	\be 
	\label{F1003} E \left(\int_{\RR} h(s_1,s_2)dX_{s_1}^1 dX_{s_2} \right)^2 = E\left(\|X \|^2_R \right).
	\ee
\end{description}
\end{cor}
\textbf{Proof}(of Corollary \ref{CRD18}): It is a consequence of Proposition \ref{PD160} setting $g=h$, with $a = 0$ and Proposition \ref{PRD17}. \qed
\newline \textbf{Proof}(of Proposition \ref{PRD17}): The fact that $g(s,\cdot)\in L_R$ for $R(ds,\infty)$ a.e. comes from the definition of $L_{2,R}$. Let $g^N \in L_R \otimes L_R$ of the type $g^N(s,t)= \sum_{i=1}^N f_i(s)h_i(t)$, $f_i,h_i \in L_R$ and 
\be
\|g^N - g \|^2_{2,R}  \xrightarrow[N \to \infty] {} 0. \label{EPDR17}
\ee
We denote by $Z(s)= \iin g(s,t)dX_t$, $Z^N(s)= \iin g^N(s,t) dX_s$. We observe that $Z^N(s)= \sum_{i=1}^N\left( \iin h_i(t)dX_t\right)f_i(s)$. Clearly $Z^N \in L_R$ a.s. The result would follow if we show the existence of a subsequence $(N_k)$ such that $Z^{N_k} \lra Z$ a.s. in $L_R$. For this it will be enough to show that
\[
E \left(\|Z^N -Z \|^2_{R}\right) \xrightarrow[N \to \infty] {}0.
\]
We have indeed
\bea 
E \left(\|Z^N-Z \|^2_{R} \right) &=& E \left(\iin (Z^N-Z)^2(s)R(ds,\infty) \right) \nn \\
	&-& \fr E\left(\int_{\RR} \left((Z^N-Z)(s_1) -(Z^N-Z)(s_2)\right)^2 d\m(s_1,s_2) \right) \nn \\
&=& \iin R(ds,\infty)E \left(\iin (g^N-g)(s,t)dX_t\right)^2\nn \\ 
&-& \fr \int_{\RR} d\m(s_1,s_2)) E\left( \iin \left[(g^N-g)(s_1,t) -(g^N-g)(s_2,t)\right] dX_t \right)^2. \nn
\eea
By Assumption (D) and Corollary \ref{cd3.9}, previous expression equals
\bea
&&\iin R(ds,\infty) \left\|(g^N-g)(s,\cdot) \right\|^2_R \nn \\
&& + \fr \int_{\RR} d\m(s_1,s_2) \left\|(g^N-g)(s_1,\cdot)-(g^N-g)(s_2,\cdot) \right\|_R^2\nn \\
&& = \|g^N-g \|^2_{2,R}. 
\eea
The result follows by \eqref{EPDR17}. \qed
%\section{Elements of Malliavin Calculus}

\section{Basic considerations on Malliavin calculus} \label{s6}

\setcounter{equation}{0}

The aim of this paper is to implement Wiener analysis in the case when our basic process $X$ fulfills Assumptions (A), (B) and (C($\n$)). In the sequel we will also often suppose the validity of Assumptions (C), (D). The spirit is still the one of \cite{KRT} in which the process $X$ was supposed to have a covariance measure structure but in a much more singular context. The target of this is the study of a suitable framework of Skorohod calculus with It\^o formulae including the case when the covariance is singular. We also explore the connection with calculus via regularization. 
\par Let $X = (X_t)_{t \in [0,\infty]}$ be an $L^2$-continuous process with continuous paths. For simplicity we suppose $X_0 = 0$. We denote by $\Co$ the set of continuous functions defined on $\R_+$ vanishing at zero with a limit at infinity. As in \cite{KRT}, we will also suppose that the law $\Xi$ of process $X$ on $\Co$ has full support, i.e. the probability that $X$ belongs to any non-empty, open subset of $\Co$ is strictly positive. This allows to state the following result.
\begin{prop}\label{p4.1} We set $\Omega_0 = \Co$, equipped with  its Borel $\sigma$-algebra and probability $\Xi$. We denote by $\Fc$ the linear space of $f(l_1, \ldots, l_m)$, $m\in \mathbb N^*$, $f \in C^\infty_b(\R^m)$, $l_1,\ldots, l_m \in \Omega_0^*$. Then $\Fc$ is dense into $L^2(\Omega_0,\Xi)$.
\end{prop}
\begin{rem} \label{r4.2}
\begin{description}
	\item[i)] A reference for this result is \cite{MR}, Section II.3.
	\item[ii)] We apply Proposition \ref{p4.1} on the canonical probability space related to a continuous square integrable process $X$. 
\end{description}
\end{rem}
We introduce a technical assumption, which will be verified in the most examples.
\be \label{4.0} 1_ {[0,t]} \in L_R, \ \ \forall t \geq 0, 1 \in L_R.
\ee
\newline For instance it is fulfilled if Assumptions (C) and (D) hold or if $X$ has a covariance measure structure, see Corollary \ref{c3.72} and Proposition \ref{r3.73}. 
\begin{rem} \label{r4.4} Taking into account \eqref{4.0}, we denote by $\bL$ the linear space of functions $\bar{f}: \R_+ \lra \R$ such that there is $f \in L_R$ with 
\be 
\bar{f}(t) =\left<f,1_{[0,t]}\right>_{\cH}. \label{4.2}
\ee
$\bL$ is the classical self-reproducing kernel space appearing in the literature.
We equip $\bL$ with the Hilbert norm inherited from $L_R$ i.e. $\| \bar{f}\|_{\bL}= \| f\|_R$. Therefore $\bar{f}_n \lra \bar{f}$ in $\bL$ if and only if $f_n \lra f$ in $L_R$. We set $\gamma_{\infty}= \sup_{t \geq 0} \sqrt{Var(X_t)}$.
Since for $0 \leq s < t,$
\[
	|\bar{f}(t)-\bar{f}(s)|= \left|E\left( (X_t -X_s) \iin f dX \right) \right| \leq \left\{ E(X_t -X_s)^2\right\}^{\fr}\|f \|_{\cH}
\]
and $X$ is continuous in $L^2(\Om)$. We have the following.
\begin{enumerate}	
	\item If $f\in L_R$, then $\sup_{t \geq 0}|\bar{f}(t)|\leq \gamma_{\infty}\| f\|_{\cH} \leq \gamma_{\infty} \| f\|_R = \gamma_{\infty} \| \bar{f}\|_{\bar{L}_R}.$
	\item $\bar{L}_R \subset C_b(\R_+)$.
\end{enumerate}
\end{rem}

\par We denote by $Cyl$ the set of smooth and cylindrical random variables of the form
\be \label{4.1}
F=f\left( \iin \phi_1 dX, \ldots, \iin \phi_m dX\right),
\ee
where $f \in C_b^\infty(\R^m)$, $\phi_1, \ldots, \phi_m \in C_0^1(\R)$ and $\iin \phi_i dX, 1\leq i\leq m$ still denotes the Paley-Wiener integral developed in Section \ref{s5}.

\par An important basic consequence of Proposition \ref{p4.1} for developing Malliavin calculus is the following.
\begin{thm}\label{t4.3} $Cyl$ is dense into $L^2(\Om)$.
\end{thm}
Before entering the proof we make some preliminary considerations. We first suppose that $\Om$ coincides with the canonical space $\Om_0$ and $X_t(\o)=\o(t)$, $t \geq 0$, so $P= \Xi$. In this case, if $f \in C_0^1(\R_+)$ (which is an element of $L_R$), the following Wiener integral
\[
\iin f dX(\o)= -\iin X_s(\o)df(s)= -\iin  \o(s)df(s)
\] 
is pathwise defined.
\begin{lem} \label{l4.5}
Let $l: \Om_0 \lra \R$ be linear and continuous. There is a sequence $(g_n)$ in $C_0^1(\R)$, $a_n \in \R$ with $\left(
\iin g_n dX \right)(\o)+a_n X_\infty \lra l(\o)$, $\forall \o \in \Om_0$ and so in particular $\Xi$ a.s.
\end{lem}
\textbf{Proof}: Since $l:\Om_0 \lra \R$ is linear and continuous, there is a finite signed, Borel measure $\ell$ on $\overline{\R}_+$ such that for every $\tilde{h} \in \Om_0$
\[
  l(\tilde{h}) = -\int_{]0,\infty[}\tilde{h} d\ell + \tilde{h}(+\infty)\ell(\{ +\infty\})
\]
Now
\[
  l(X)= -\iin X d\ell + X_\infty \ell(\{ +\infty\})
\]
We set
\bea
&&g_n(x)=\int_0^n \rho_n(x-y)g(y)dy, \nn\\
&&g(x)=\ell([x,\infty[). \nn
\eea
where $(\rho_n)$ is the usual sequence of mollifiers with compact support approaching Dirac delta function. In particular $dg_n \lra \ell|_{\R_+}$ weakly. We set
\[
l_n(\tilde{h})= -\iin \tilde{h}dg_n, \ \ a_n= \ell(\{+\infty\})
\]
so that
\[
l_n(X)=\iin g_n dX.
\]
Since $l_n(\tilde{h}) \lra -\iin \tilde{h} d
\ell$ pointwise, the result follows. \qed

\begin{lem} \label{r4.6} The statement of Theorem \ref{t4.3} holds whenever  $\Om=\Om_0$, $P=\Xi$.
%Lemma \ref{l4.5} implies the result of the theorem whenever. F (so a.s.) by a sequence of $f(\iin \phi_0 dX,\ldots, \iin \phi_m dX)$ for some $\phi_0, \ldots, \phi_n \in C^1_0(\R_+)$ and $\phi_0 =1$. Since $f$ is bounded, the convergence also holds in $L^2(\Om, \Xi)$.
\end{lem}
\textbf{Proof}: By Proposition \ref{p4.1} it is enough to show that any element of $\mathcal{F}C^{\infty}_b$ can be approached by a sequence of random variables in $Cyl$. Let $F \in \mathcal{F}C^{\infty}_b$ given by $f(l_1, \ldots, l_m)$ as in Proposition \ref{p4.1}. By truncation it is clear that we can reduce to the case, when $f$ is bounded. Lemma \ref{l4.5} implies that it can be pointwise approximated (so a.s.) by a sequence of random variables of the type $f(\iin \phi_0 dX,\ldots, \iin \phi_m dX)$ for $\phi_0= 1$, and $\phi_1 \ldots, \phi_n \in C^1_0(\R_+)$. Since $f$ is bounded, the convergence also holds in $L^2(\Om;\Xi)$. \qed
\newline
\textbf{Proof} (of Theorem \ref{t4.3}): 
Any r.v. $h \in \LOm$ can be represented through $F(X)$, where $F \in L^2(\Om_0,\Xi)$. According to Lemma \ref{r4.6} there is a sequence of elements of the type $f(\iin \phi_0 di, \ldots , \iin \phi_m di)$, $i(s)=s$, $\phi_0 =1,\phi_1, \ldots, \phi_n \in C_0^1(\R_+)$, $f \in C_b^\infty(\R)$ converging in $L^2(\Om_0, \Xi)$ to $F$. Since Wiener integrals $\iin \phi_j dX$, $1 \leq j \leq n$, can be pathwise represented, then 
\[
f\left(\iin \phi_0 dX, \ldots , \iin \phi_m dX\right) = \left.f\left(\iin \phi_0 di, \ldots , \iin \phi_m di\right)\right|_{i=X}
\]
and the general result follows. \qed

\section{Malliavin derivative and related properties} 
\label{s7}

\setcounter{equation}{0}

In this section we suppose again Assumptions (A), (B) and (C($\n$)). We will suppose from now on that $X$ is Gaussian.
We start with a technical lemma.
\begin{lem} \label{LVNonDeg} Let $\phi_1, \ldots \phi_n \in C^1_0(\R_+)$ orthogonal with respect to $\langle \cdot, \cdot \rangle_\cH$,
 not vanishing. Then then the law of the vector 
\[
V= \left(\iin \phi_1 dX, \ldots, \iin \phi_m dX \right)
\]
has full support in the sense that for any non empty open set $I$ of $\R^n$, $P\left\{V \in I \right\}>0$.
\end{lem}
\textbf{Proof}: Clearly we can reduce the question to the case $I = \prod_{j=1}^m ]a_j,b_j[$, $a_j<b_j$. Since the random variables $\iin \phi_1 dX, \dots, \iin \phi_m dX$ are independent, it is enough to write the proof in the case $m=1$, $\phi= \phi_1 \in C^1_0(\R)$, $\phi \neq 0$. 
\newline Let $\Xi$ be the law of $X$ on $\Om_0$. Then
\be \label{Ephi} 
P\left\{\iin \phi dX \in ]a_1,b_1[ \right\} = \Xi\left\{ \o \in \Om_0| l(\o)\in ]a_1,b_1[ \right\},
\ee
where $l(\o)= -\int \o d\phi$, which is clearly an element of the topological dual $\Om_0^*$. Since $l$ is continuous,
\be \label{EPhi1}
\left\{\o \in \Om_0 | l(\o) \in ]a_1,b_1[ \right\}
\ee
is an open subset of $\Om_0$.
Since $\Xi$ has full support, it remains to show that the set \eqref{EPhi1} is non empty.
\newline It is always possible to find $\o_0 \in \Om_0$ such that $l(\o_0) \neq 0$. Otherwise the derivative $\dot{\phi}$ would be orthogonal with respect to the $L^2(\R_+)$ norm to the linear space
\[
\{\o \in L^2(\R_+) \cap C(\R_+)| \o(0)=0 \}.
\]
This would  not be possible since that space is dense in $L^2(\R_+)$. Consequently, there exists $\l \in \R$ such that
\[
l(\l \o_0)= \l l(\o_0) \in ]a_1,b_1[.
\]
It is enough to choose $\l$ between $\frac{a_1}{l(\o_0)}$ and ${\frac{b_1}{l(\o_0)}}$. Finally $\l \o_0$ belongs to the set defined in \eqref{EPhi1}. \qed
\newline For $F \in Cyl$ of the form (\ref{4.1}), we define
\be
D_tF = \sum_{i=1}^n \pa_i f\left(\iin \phi_1 dX, \ldots, \iin \phi_m dX \right) \phi_i(t). \label{ExprD}
\ee
\begin{rem} \label{r421} Let $F \in Cyl$. Since $\phi_i \in L_R$ 
and $\pa_i f, 1 \le i \le n $ are bounded, then $t \lra D_tF \in L_R$ a.s. Moreover
\[
E \left( \|DF\|^2_R\right) < \infty.
\]
Consequently $DF\in L^2(\Om;L_R)$.
\end{rem} 

\begin{prop} \label{R4210} Expression (\ref{ExprD}) does not depend on the explicit form (\ref{4.1}).
\end{prop}
\textbf{Proof: } We can of course reduce the problem as follows. Let $f: \R^n \lra \R$, $f \in C^\infty_b(\R^n)$, $\phi_1,\ldots,\phi_n \in C_0^1$ such that
\[
f(Z_1, \ldots, Z_n) =0,
\]
where $Z_i = \iin \phi_i dX$. We need to prove that
\be
\label{ER4210} \sum_{k=1}^n \pa_k f(Z_1, \ldots, Z_n) \phi_k =0 \quad a.s.
\ee
By a classical orthogonalization procedure with respect to the inner product $\left<\cdot,\cdot\right>_{\cH}$, there is $m \leq n$, $A= (a_{ij})_{1 \leq i \leq n, 1\leq j\leq m}$ such that $\phi_i = \sum_{j=1}^m a_{ij} \psi_j$, $\psi_1, \ldots \psi_m $ being orthogonal. Writing $Y_j = \iin \psi_jdX$, we also have
\be
\tilde{f}(Y_1,\ldots, Y_m) = 0, \label{ER4209}
\ee 
with
\[
\tilde{f}(y_1,\ldots,y_m) = f\left(\sum_{j=1}^m a_{1j}y_j, \ldots, \sum_{j=1}^m a_{nj}y_j\right).
\]
By usual rules of calculus, (\ref{ER4210}) implies that
\be \sum_{l=1}^m \pa_l \tilde{f}(Y_1, \ldots, Y_m)\psi_l = \sum_{k=1}^m \pa_k f(Z_1, \ldots, Z_n)\phi_k. \label{ER4211}
\ee
\eqref{ER4209} implies that
\be \label{ER4500} 
\int_{\R^m} \tilde{f}^2(y_1, \ldots, y_m) d\m_{V}(y_1,\ldots,y_n) = 0,
\ee
where $\m_{{V}}$ is the law of $(Y_1, \ldots, Y_n)$. Since $\tilde{f}$ is continuous and because of Lemma \ref{LVNonDeg}, \eqref{ER4500} implies that $\tilde{f}\equiv 0$. This finally allows to conclude (\ref{ER4210}). \qed

Before going on, we need to show that $D: Cyl \lra L^2(\Om)$ is closable. We observe first the following property.
\begin{prop} \label{PR4211} Let $F \in Cyl$, $h\in L_R$. Then 
\[
E(\left<DF,h\right>_\cH)= E\left( F\iin h dX\right).
\]
\end{prop}
\textbf{Proof:} It is very similar to Lemma 6.7 of \cite{KRT} or Lemma 1.1 of \cite{Nualart}. \qed
\begin{rem} \label{r424} $Cyl$ is a vector algebra. Moreover, if $F,G \in Cyl$, then 
\be \label{E424} D(FG) = GDF + FDG. 
\ee
\end{rem}
A consequence of Proposition \ref{PR4211} and Remark \ref{r424} is the following.
\begin{cor} \label{CR4212} Let $F,G \in Cyl$, $h \in L_R$. Then
\bea
&&E( G\left<DF,h\right>_{\cH})\nn\\ 
&& \ \ \ \ \ = E(-F\left<DG,h\right>_{\cH}) + E(FG\iin h dX).\nn
\eea
\end{cor}
Finally we can state the following result.
\begin{prop} \label{PR4213} The map $D: Cyl\lra L^2(\Om;L_R)$ is closable.
\end{prop}
\textbf{Proof:} Let $F_n$ be a sequence in $Cyl$ such that $\lim_{n \to \infty} E(F_n^2)=0$ and there is $Z \in L^2(\Om;L_R)$ such that $\lim_{n \to \infty} E(\| DF_n -Z\|_R^2)=0$. We need to prove that $Z = 0$ a.s. It is enough to show that $\|Z \|^2_{\cH}=0$ a.s. Since $\cH$ is separable and $C_0^1$ is dense in $\cH$, it is enough to show that $\left<Z,h\right>_{\cH} = 0$ a.s. Since $Cyl$ is dense in $L^2(\Om)$, we only have to prove that
\[
E(\left<Z,h\right>_{\cH}G)=0 \ \forall G \in Cyl.
\]
By Corollary \ref{CR4212}, previous expectation equals
\bea && \lim_{n \to \infty} E\left( \left<DF_n,h\right>_{\cH}G\right) \nn \\
&& \ \ \ \ =\lim_{n\to \infty}\left(E\left( -F_n\left<DG,h\right>_{\cH}\right)+ E\left(F_n G \iin h dX \right)  \right) =0.
\eea
This concludes the proof of the proposition. \qed
\newline
We denote by $\DD$ the space constituted by $F\in \LOm$ such that there is a sequence $(F_n)$ of the form (\ref{4.1}) verifying the following conditions.
\begin{description}
	\item[i)] $F_n \lra F$ in $\LOm$,
	\item[ii)] $E\left(\|DF_n - Z\|_R^2 \right) \xrightarrow[n \to \infty] {} 0$,
\end{description}
for some $Z \in L^2(\Om;L_R)$. In agreement with Proposition \ref{PR4213}, we denote $DF = Z$.
\newline The set $\Dd$ will stand for the vector subspace of $\LOm$ constituted by functions $F$ such that there is a sequence $(F_n)$ of the form (\ref{4.1}) with
\begin{description}
	\item[i)] $F_n \lra F$ in $\LOm$,
	\item[ii)] $E\left(\|DF_n - DF_m\|_{\cH}^2\right) \xrightarrow[n,m \to \infty] {} 0$.
\end{description}
Note that $\DD \subset \Dd$. $\DD$, equipped with the scalar product
\[
\left<F,G\right>_{1,2} = E(FG) + E(\left<DF,DG\right>_R)
\]
is a Hilbert space. 
\newline From previous definitions we can easily prove the following. 
\begin{prop} \label{AfterdefD} Let $(F_n)$ be a sequence in $\DD$ (resp. $\Dd$), $F \in L^2(\Om)$, $\mathcal{Y} \in L^2(\Om;L_R)$ such that
\[
E\left((F_n-F)^2+\|DF_n - \mathcal{Y} \|^2_R\right) \xrightarrow[n \to \infty] {} 0.
\]
(resp.
\[
 E\left((F_n-F)^2+\|D(F_n - F_m)\|^2_{\cH}\right) \xrightarrow[m,n \to \infty] {} 0. )
\]
Then $F \in \DD$ and $\mathcal{Y}=DF$ (resp. $F \in \Dd$).
\end{prop}

\begin{rem} \label{r422} If Assumption (D) is fulfilled, then $\DD = \Dd$ and
\[
\left<F,G\right>_{1,2} = E(FG) + E(\left<DF,DG\right>_{\cH}).
\] 
\end{rem}
\begin{rem} \label{r423} The notation $\DD$ does not have the same meaning as in \cite{KRT}. Indeed $\|\cdot \|_{|\cH|}$ introduced there is not exactly a norm.
\end{rem}
\begin{rem} \label{rr4212} By definition of $\Dd$ the statement of Corollary \ref{CR4212} extends to $F, G \in \Dd$. We have therefore the following
\[
E \left( G \left<DF,h\right>_{\cH}\right) = E \left( -F\left<DG,h\right>_{\cH}\right)+ E\left(FG \iin h dX \right)
\]
for every $F,G \in \Dd$, $h \in L_R$.
\end{rem}
\begin{prop} \label{PPol} Let $f:\R \lra \R$, be absolutely continuous. Let $\phi \in L_R$. 
%We suppose the validity of Assumption (D) (EST-CETTE HYPOTHESE UTILE???). 
We suppose that $f'$ is subexponential. 
\newline Then $f\left(\iin \phi dX \right) \in \Dd$ and 
\[
D_r f\left( \iin \phi dX\right) = f' \left( \iin \phi dX\right) \phi(r).
\]
\end{prop}
\begin{rem} \label{RPol} 
\begin{enumerate}
	\item A function $g: \R^n  \lra \R$ is said to be \textbf{subexponential} if there is $\g >0$, $c > 0$ with $|f(x)| \leq c e^{\g|x|}$, 
$\forall x \in \R^n$.
	\item In particular if $f$ is a polynomial, previous result holds.
\end{enumerate}
\end{rem}
\textbf{Proof} (of Proposition \ref{PPol}) 
i) We first suppose $f \in C_b^{\infty}(\R)$. There is a sequence $\phi_n$ in $C_0^1$ such that $\| \phi-\phi_n\| \xrightarrow[n \to \infty] {} 0$. Clearly
\[
E\left( f\left( \iin \phi dX \right)- f\left(\iin \phi_n dX \right)\right)^2 \xrightarrow[n \to \infty] {} 0
\]
since
\[
\iin \phi_n dX \lra \iin \phi dX
\]
in $L^2(\Om)$ and by Lebesgue dominated convergence theorem. On the other hand 
\[
D_t f\left( \iin \phi_n dX\right) = f'\left( \iin \phi_n dX \right) \phi_n(t), \ \ t \geq 0,
\]
so 
\bea
&&E\left(\left\| Df\left(\iin \phi_n dX \right)-f'\left( \iin \phi dX\right) \phi\right\|^2_R \right) \nn \\
&& \leq \| \phi-\phi_n\|^2_R E\left( f'\left( \iin \phi dX\right)^2\right) + E\left(f'\left(\iin \phi_n dX \right) - f'\left(\iin \phi dX \right) \right)^2\|\phi_n \|^2_R. \nn
\eea
This converges to zero by usual integration theory arguments. The result
 for $f \in C^\infty_b(\R)$ follows by Proposition \ref{AfterdefD}.
\newline
ii) We suppose now that $f'$ is subexponential nad let $\phi \in L_R$. $\iin \phi dX$ is Gaussian zero-mean variable with covariance $\s^2= \| \phi\|^2_R$. In fact it is the limit in $\LOm$ of r.v. of the type $\iin \phi_n dX$, $\phi_n \in C^1_0$.
We proceed setting $\tilde{f}'_M = (f'\wedge M)\vee(-M)$ for $M >0$ and
\[
\tilde{f}_M:= f(0) + \int_0^x \tilde{f}'_M(y)dy.
\]
We also set
\[
f_M(x) = \int_{\R} \r_{\frac{1}{M}}(x-y)\tilde{f}_M(y)dy,
\]
where $\r_{\eps}$ is a sequence of Gaussian mollifiers converging to the Dirac delta function. It is easy to show that
\be
\int_{\R}(f_M-f)^2(x)p_{\sigma}(x)dx \xrightarrow[M \to \infty] {} 0, \label{FDD12}
\ee
\be
\int_{\R}(f'_M-f')^2(x)p_{\sigma}(x)dx \xrightarrow[M \to \infty] {} 0, \label{FDD12der}
\ee
where $p_{\sigma}$ is the density related to the Gaussian law $N(0,\sigma^2)$. \eqref{FDD12} implies that
\[
(f_M-f)\left( \iin \phi dX\right)\lra 0
\]
in $\LOm$. By point i) of the running proof we have
\[
D_r f^M\left(\iin \phi dX \right)= \left(f^M \right)'\left(\iin \phi dX \right) \phi(r).
\]
\eqref{FDD12der} implies
\[
E\left(\left\|D_{\cdot} f_M\left( \iin \phi dX\right) - f'\left( \iin \phi dX\right)\phi \right\|^2 \right) \xrightarrow[M \to \infty]{} 0.
\]
which together with Proposition \ref{AfterdefD} clearly gives the result.
%\newline iii) We suppose $f \in C^1_b (\R)$. For $\eps >0$, we define
%\[
%f_{\eps}(x) = \int_{\R} \r_{\eps} (x-y)f(y)dy,
%\] 
%where $(\r_{\eps})$ is as at point ii). It follows that
%\[
%f'_{\eps}(x) = \int_{\R} \r_{\eps}(x-y) f'(y) dy.
%\]
%There is a constant $C > 0$ such that,
%\bea
%|f_{\eps}(x)| &\leq & C(1 + |x|) \nn \\
%|f_{\eps}(x)| &\leq & C, \ \ \forall x \in \R \nn \\ 
%\eea
%for every $\eps > 0$. Since $f_{\eps} \to f$, $f'_{\eps} \to f'$ pointwise, by use of Lebesgue dominated convergence theorem and again by point i) of the present proof, the result follows again via Proposition \ref{AfterdefD}.
\qed

Proposition \ref{PPol} extends to the case, where $f$ depends on more than one variable. The proof is a bit more complicated, but it follows the same idea. Therefore we omit it.
\begin{prop} \label{PPolN}
Let $f: \R^n \lra \R$ of class $C^1$, with subexponential partial derivatives. Let $\phi_1, \ldots, \phi_n \in L_R$.
%ASSUMPTION D (VRAIMENT???) 
Then
\[
f\left( \iin \phi_1 dX, \ldots, \iin \phi_n dX\right) \in \Dd
\]
and
\bea
&&D_r f\left(\iin \phi_1 dX, \ldots, \iin \phi_n dX \right) \nn \\ 
&& = \sum_{j=1}^n \pa_j f\left( \phi_1 dX, \ldots, \iin \phi_n dX\right) \phi_j(r).
\eea
\end{prop}
We establish some immediate properties of the Malliavin derivative.
\begin{lem} \label{l425} Let $F \in Cyl$, $G \in \DD$. Then $F \cdot G \in \DD$ and (\ref{E424}) still holds.
\end{lem}
\textbf{Proof}: According to the definition of $\DD$, let $(G_n)$ be a sequence in $Cyl$ with the following properties.
\begin{description}
	\item[i)] $E(G_n -G)^2 \xrightarrow[n \to \infty]{} 0,$
	\item[ii)] $E \left( \iin (D_r(G_n -G))^2|R|(dr,\infty)\right) \xrightarrow[n \to \infty]{} 0,$ 
	\item[iii)] $E \left( \int_{\RR}d|\m|(r_1,r_2)(D_{r_1}(G_n -G)-D_{r_2}(G_n-G))^2\right) \xrightarrow[n \to \infty]{} 0.$
\end{description}
Since $F \in L^\infty (\Om)$ then $FG_n \lra FG$ in $\LOm$. Remark \ref{r424} implies that
\[
D(F G_n) = G_n DF + FDG_n.
\]
It remains to show ii) and iii) for $G_n$ (resp. $G$) replaced with $FG_n$. We only check ii), because iii) follows similarly.
If $F$ is of the type (\ref{4.1}) then
\[
DF = \sum_{i=1}^m Z_i \phi_i,
\]
where $\phi_i \in L_R$, $Z_i \in L^\infty(\Om)$. This implies, by subadditivity, that
\[
\iin |R|(dr,\infty)(D_rF)^2 \leq 2^m\left(\sum_{i=1}^m \| Z_i\|^2_\infty \iin \phi_i^2(r)|R|(dr,\infty)\right).
\]
Hence
\bea
&& E\left( \iin |R|(dr,\infty)(G_n-G)^2(D_r F)^2\right) \nn \\
&&\leq E(G_n -G)^2 \max_{i \in \{1, \ldots , m \}}\| Z_i\|_{\infty}^2 \left( 2^m \sum_{i=1}^m \iin \phi_i^2(r)|R|(dr,\infty)\right) \xrightarrow[n \to \infty]{} 0.\nn
\eea
Moreover, since $F \in L^{\infty}$ and taking into account ii)
\[
E\left( \iin |R|(dr,\infty)(F D_r(G_n-G))^2\right) \xrightarrow[n \to \infty]{} 0.
\]
Hence ii) is proven for $F(G_n - G)$ instead of $G_n-G$. \qed

A natural question is the following. Does $X_t$ belong to $\Dd$ for fixed $t$? The proposition and corollary below partially answers the question.
\begin{prop} \label{PTX}
If $\psi \in L_R$, then $\iin \psi dX \in \DD$ and $D_t\left( \iin \psi dX \right) = \psi(t)$.
\end{prop}
\textbf{Proof}: We consider a sequence $(\psi_n)$ in $C_0^1(\R)$ such that $\|\psi -\psi_n \|_R\xrightarrow[n \to \infty]{}0.$ We know that $\iin \psi_n dX \in Cyl$. Obviously
 \[
 E\left( \iin (\psi_n-\psi)dX \right)^2 = \| \psi -\psi_n\|^2_{\cH} \leq \|\psi -\psi_n \|^2_R \xrightarrow[n \to \infty]{} 0.
 \]
On the other hand
\[
D_r\iin \psi_n dX = \psi_n(r), \  r>0,
\]
so
\[
E \left( \left\| D\iin (\psi_n -\psi_m) dX\right\|^2_R \right) = \| \psi_n -\psi_m \|^2_R \xrightarrow[n,m \to \infty]{}0. 
\] \qed
\begin{cor} \label{CXI} If $1_{[0,t]} \in L_R$, then $X_t \in \DD$ and $DX_t = 1_{[0,t]}$.
\end{cor}
\begin{rem} \label{RXII} The conclusion of Corollary \ref{CXI} holds if Assumptions (C) and (D) hold, see Corollary \ref{c3.72}.
\end{rem}

\section{About vector valued Malliavin-Sobolev spaces} \label{s8bis}

\setcounter{equation}{0}

We suppose here the validity of Assumption (C) and we use the notations introduced in Section \ref{s8}. We denote again $\n(dt) = |R|(dt,\infty)$. 

\par We will first define $Cyl(L_R)$ as the set of smooth cylindrical random elements of the form
\[
u(t)  = \sum_{l=1}^n \psi_l(t)G_l, \ t \in \R_+,
\]
$G_l \in Cyl$, $\psi_l \in C_0^1(\R_+)$. If $u \in Cyl(L_R)$, we define
\[
\tilde{D}_s u(t) = \sum_{l=1}^n \psi_l(t)D_sG_l, s,t \geq 0.
\]
Clearly $\tilde{D}u = (\tilde{D}_s u(t))$ belongs to $L_{2,R}$ for each underlying $\omega \in \Omega$.
\begin{rem} \label{R4520} 
\begin{enumerate}
	\item If $u \in Cyl(L_R)$, it is easy to see that a.s. the paths of $\tilde{D}u$ belong to $L_R \otimes L_R$.
	\item   Taking into account Assumption (C), if $u \in Cyl(L_R)$, then $u(t) \in \DD$.
	\item By analogous arguments as in Proposition \ref{PR4213}, is is possible to show that $\tilde{D}: Cyl(L_R) 		\lra L^2(L_{2,R})$ is well-defined and closable. This allows to set $Z = \tilde{D}u$, called the \textbf{Malliavin derivative} of process $u$.
\end{enumerate}
\end{rem}
Similarly to $\DD$ we will define $\DDd$. We denote $\DDd$ the vector space of random elements $u: \Om \lra L_R$ such that there is a sequence $(u_n)$ in $Cyl(L_R)$ and 
	\begin{description}
		\item[i)] $\| u -u_n\|^2_R \xrightarrow[n \to \infty]{}0$ in $\LOm$.
		\item[ii)] There is $Z:\Om \lra L_{2,R}$ with $\|\tilde{D}u_n -Z \|_{2,R} \lra 0$ in $\LOm$.  
	\end{description}
We denote $Z$ again by $\tilde{D}u$.
%\begin{cor} \label{PCYR} Under Assumption (A), (B), (C), (D) and (\ref{ELTR}). We have 
%\[
%f(X) \in \Dd(L_R)
%\]
%and
%\[
%D_r f(X_t) = f'(X_t)1_{[0,t]}(r).
%\]
%\end{cor}
%\textbf{Proof:} This is a consequence of Proposition \ref{PYR} PLUS LOIN p. 4.55B!!!!!!!! and Proposition \ref{PD12}.\par We observe that $u_n(t) \in Cyl$ for every $t \geq 0$. By definition of $\tilde{D}$ on $Cyl(L_R)$ we have
%Finally the result follows. 
\begin{rem}\label{R4530}
\begin{description}
	\item[i)] Let $(u_n)$ be a sequence in $\DDd$, $u \in L^2(\Om; L_R)$, $Z \in L^2(\Om; L_{2,R})$. If 
	\[
	\lim_{n\to \infty} E \left(\|u-u_n \|_R^2 + \|\tilde{D}u_n -Z \|^2_{2,R} \right) =0,
	\] 
	it is not difficult to show that $u \in \DDd$ and
	\[
	\tilde{D}u = Z.
	\] 
	\item[ii)] Let $u_t = \psi(t)G$, $t \geq 0$; $G\in \Dd$, $\psi \in L_R$. Then $u \in \DDd$. Moreover $\tilde{D}_r u(t )= \psi(t) D_r G$, $r \geq 0$ . This follows by point i) and the fact that $u$ can be approximated by $u_t^n = \psi_n(t)G_n$, where $\psi_n \in C_0^1$ and $G_n \in Cyl$.
\end{description}
\end{rem}
\begin{rem} 
\label{r453}
\begin{enumerate}
	\item $\DDd$ is a Hilbert space if equipped with the norm $\| \cdot\|$ associated with the inner product
	\[
	\left< u,v \right> = E\left(\left<u ,v\right>_R + \left<\tilde{D}u,\tilde{D}v\right>_{2,R}\right).
	\]
	Moreover $Cyl(L_R)$ is dense in $\DDd$.
	\item We convene here that
	\[
	\tilde{D}u: (s,t) \longmapsto \tilde{D}_s u(t).
	\]
	\item If Assumption (D) is fulfilled, it is possible to show that $\DDd = \Ddd$, where $\Ddd$ is constituted by the vector space of random elements $u:\Om \lra L_R$ such that there is a sequence $(u_n)$ is $Cyl(L_R)$ with the following properties
	\begin{description}
		\item[i)] $\| u-u_n\|^2_{\cH} \xrightarrow[n \to \infty]{}0 $ in $\LOm$.
		\item[ii)] 	There is $Z:\Om \lra L_R \otimes^h L_R = L_{2,R}$ with
		\[
		\| \tilde{D}u_n-Z\|_{2,R} \xrightarrow[n \to \infty]{}0
		\]
		in $\LOm$.
	\end{description}
	\item If there is a sequence $u_n$ verifying points i), ii), then it is not difficult to show that $u \in \Ddd$. Of course $Z = \tilde{D}u$.
\end{enumerate}
\end{rem}
We focus the attention on some technical point. The derivative $\tilde{D}$ of process $(u(t))$ may theoretically not be compatible with the family of derivatives of random variables $u(t)$.
\begin{prop} \label{P4521} Let $u \in \DDd$. Then $\nu(dt)$ a.e. $u(t)\in \DD$ and
\[
D_r u(t) = \tilde{D}_r u(t),\  \n\otimes \n \otimes P. \ a.e.
\]
\end{prop}
\textbf{Proof:} Since $u \in \DDd$, there is a sequence $u_n \in Cyl(L_R)$ such that $u_n \lra u$ in $\DDd$. According to \eqref{EI4.111} and Point ii), it follows that
\[ 
E \left(\iin \n(dt) \| \tilde{D}_{\cdot}u_n(t)- \tilde{D}_{\cdot}u(t)\|_R^2) \right) \xrightarrow[n \to \infty]{}0.
\]
Consequently $\n(dt)$ a.e. we have
\[ 
E \left( \| \tilde{D}_{\cdot}u_n(t)- \tilde{D}_{\cdot}u(t)\|_R^2 \right) \xrightarrow[n \to \infty]{}0.
\]
By a similar argument, it follows that
\[
\lim_{n \to \infty} \iin \n(dt) E(u_n(t)-u(t))^2 = 0.
\]
We observe that $u_n(t) \in Cyl$ for every $t \geq 0$. By definition of $\tilde{D}$ on $Cyl(L_R)$, we have 
\[
\tilde{D}u_n(t) = Du_n(t).
\]
Finally the result follows. \qed
\par From now on we will not distinguish between $D$ and $\tilde{D}$.
\par A delicate point consists in proving that the process $X \in \Ddd$. First we state a lemma.
\begin{lem} \label{l4.55quater} We suppose Assumption (D). Let $g \in C^1$ such that there is $T> 0$ with $g(t) = g(T)$, $t \geq T$. Then $g \in L_R$ and for every $f \in C^1_0$
\be
\label{E4551} \left< f,g\right>_R = \int_{\RR} f'(s_1)g'(s_2) R(s_1,s_2) ds_1 ds_2.
\ee
\end{lem}
\textbf{Proof}: We consider a family of functions $\chi^n$ in $C^{\infty}_b(\R_+)$ such that $\chi^n =1$ on $[0,n]$ and $\chi^n =0$ on $[n+1,\infty]$. We define $g_n = g \chi^n$. For $n,m$, $n>m$, we have
\bea
\|g_n-g_m \|^2_R &=& E \left(\iin (g_n-g_m) dX\right)^2 \nn \\
&=& E \left( \int_{]0,\infty[^2} X_{s_1}X_{s_2} d(g_n-g_m)(s_1)d(g_n-g_m)(s_2)\right) \xrightarrow[n,m\to \infty] {}0. \nn
\eea
This shows that $g_n$ is Cauchy; $g_n$ is also Cauchy in $L^2(\n)$. Consequently, there is a subsequence $(n_k)$ such that $g_{n_k} \lra g$ in $L^2(d\n)$. Since $g_n \lra g$ pointwise, then $g \in L_R$. By Remark \ref{r3.1}, we recall that \eqref{E4551} holds for every $f,g \in C^1_0$. Therefore it holds for $f$ and $g_n$. Letting $n \lra \infty$ on both sides, the result follows. \qed
\newline We operate now a restriction on $X$, supposing the existence of $T>0$ with $X_t= X_T$ if $t \geq T$. 
\begin{prop} \label{l4.55ter} We suppose Assumption (D), \eqref{ELTR} and $X_t = X_T$, $t \geq T$. Let $f \in L_R$; there is $\f = \f_f \in L_R$ such that
 \be
 \left<f,X\right>_{R} = \iin \f dX \textrm{ a.s.} \label{EEE10}
 \ee
\end{prop}
\textbf{Proof}: 
By Lemma \ref{l4.55quater}, we observe for every $f \in C^1_0(\R)$, $g \in C^1,$ constant after some $T>0$, we observe 
\be
\label{E4552} \left< f,g\right>_R = -\iin d\f_f(s)g(s),
\ee
where
\[
\f_f(s) = \iin R(s_1,s)f'(s_1)ds_1.
\]
Taking into account Assumption (A), $\f_f$ has bounded variation. 
\par The next step will be to prove that
\be
\left< f,X\right>_R = -\iin d\f_f (s) X_s, \forall f\in C^1_0(\R). \label{E4556}
\ee
We will set $g = X$.
\newline 1) We denote $h(s_1,s_2)= 1_{[0,s_1 \wedge T]}(s_2)$ and we consider again the approximating sequence $(F^{\eps})$ as in the proof of Proposition \ref{PD160}. We recall that  each $F^{\eps}$ verify has bounded planar variation and therefore, belongs to $L_{2,R}$. We also had
\[
F^{\eps} \lra h
\]
in $L_{2,R}$. By construction it also converges pointwise.
\newline  2) Let $X^1$ be an independent copy of $X$. By isometry of the double Wiener integral, it follows that
\be \label{EEG2} E \left( \int_{\RR} \left(F^{\eps}(t_1,t_2)-h(t_1,t_2) \right)dX^1_{t_1}dX_{t_2}\right)^2 \xrightarrow[\eps \to 0]{} 0.
\ee
3) Taking into account Remark \ref{RDPBV} and Proposition \ref{p3.91}, we can easily show that
\[
\int_{\RR} F^{\eps} (t_1,t_2) dX^1_{t_1} dX_{t_2} = - \iin dX^1_{t_1} \iin X_{t_2} \frac{\pa F^{\eps}}{\pa t_2} (t_1,t_2) dt_2,
\]
where $F^\eps$ is given in \eqref{EEPD160}.
\newline4)  By Proposition \ref{PRD17} and item 3), we have
\[
\int_{\RR} \left(F^{\eps}(t_1,t_2) -h(t_1,t_2) \right) dX^1_{t_1} dX_{t_2} = - \iin dX^1_{t_1} \Phi^{\eps}(t_1,X),
\]
where
\bea
\Phi^\eps(t_1,x) &=& \iin dt_1 x(t_2) \frac{\pa F^{\eps}}{\pa t_1}(t_1,t_2) - x(t_1) \nn \\
&=& \frac{1}{\eps} \iin dt_2 x(t_2) \r\left( \frac{t_1 - t_2}{\eps}\right) - x(t_1). \label{EEG4}
\eea
\eqref{EEG2} gives
\be \label{EEG3} E\left(R^{\eps}(X)^2 \right) \xrightarrow[\eps \to 0] {} 0,
\ee
where
\[
R^{\eps}(x) = \iin dX^1_{t_1} \Phi^{\eps}(t_1,x).
\]
Taking the conditional expectation with respect to $X$, we get
\bea
 &&E \left(R^{\eps}(X)^2 \right) = E \left(\tilde{R}^{\eps}(X) \right), \nn \\
 &&\tilde{R}^{\eps}(x) = E\left(R^{\eps}(x) \right)^2 = \|\Phi^{\eps}(\cdot,x) \|^2_R.  \nn
\eea
Therefore there is a sequence $(\eps_n)$ such that
\[
\left\|\Phi^{\eps_n} (\cdot,X)\right\|^2_R \xrightarrow[n \to \infty]{}0 \ a.s.
\]
Setting
\[
X^{\eps}_t = \Phi^{\eps}(\cdot,X),
\]
we have shown that $\|X^{\eps}-X \|_R \xrightarrow[\eps \to 0] {} 0$.
\newline 5) By \eqref{EEG4}, obviously $X^{\eps} \lra X$ pointwise a.s.
\newline 6) By \eqref{E4552}, we have
\be
\left<f,X^{\eps} \right>_R = - \iin X^{\eps}_s d \f_f(s), \forall f \in C^1_0(\R). \label{E4557}
\ee
Since $X^{\eps} \lra X$ a.s. in $L_R$, $X^{\eps} \lra X$ pointwise. Lebesgue's dominated convergence theorem allows to take the limit, when $\eps \lra 0$ in \eqref{E4557}. This establishes \eqref{E4556}.
\newline In order to conclude the validity of \eqref{EEE10}, taking into the isometry property of stochastic integral, we need to show that the linear operator $f \longmapsto \iin \f_f dX$ from $C^1_0$ to $\LOm$ is continuous with respect to $\|\cdot \|_{R}$.
\par Let $(f_n)$ be a sequence in $C^1_0$ converging to $0$ according to the $L_R$-norm. Corollary \ref{CRD18} implies that $X \in L^2(\Om;L_R)$. Cauchy-Schwarz implies that
\bea
E \left(\iin \f_f dX \right)^2 &=& E\left( \left<f_n,X \right>_{R}^2\right) \nn \\
&\leq& \| f_n\|^2_R E\left(\|X \|^2_{R} \right) \xrightarrow[n \to \infty]{}0. \nn
\eea
This concludes the proof of \eqref{E4556}.
 \qed
\begin{prop} \label{R4550} Under Assumptions (C) and (D) and again \eqref{ELTR} together with $X_t = X_T$, $t \geq T$, we have
\[
X \in \Dd(L_R)
\]
and
\[
D_{t_2}X_{t_1} = 1_{[0,t_1 \wedge T]}(t_2).
\]
\end{prop}
\textbf{Proof:} Let $(e_n)$ be an orthonormal basis of $L_R= \cH$ which is separable by Proposition \ref{p376}. By Corollary \ref{CRD18}, $X \in L_R$, so
\[
X = \sum_{i=1}^{\infty} F_i e_i \textrm{ in } \cH \ a.s.,
\]
where
\[
F_i= \left<X,e_i\right>_{\cH}.
\]
We recall that 
\bea
E(\|X \|^2_{\cH})&=& E\left(\iin R(ds,\infty)X^2_s-\fr \int_{\RR}d\m(s_1,s_2) (X_{s_1}-X_{s_2})^2 \right) \nn \\ 
&=& \iin R(ds,\infty)E(X^2_s)+ \fr \int_{\R^2}d(-\m)(s_1,s_2)Var(X_{s_1}-X_{s_2}) \nn
\eea
which is finite by assumption. Since
\[
\| X\|^2_{\cH} = \sum_{i=0}^{\infty} |F_i|^2, 
\]
taking the expectation we get
\[
\sum_{i=0}^{\infty} E(F_i^2) < \infty.
\]
This shows
\be \label{ERAO} \lim_{n\to \infty} \|X-X^n \|^2_R = 0.
\ee
It remains to show that the sequence $(D(\sum_{i=0}^{n}F_i e_i))_{n \geq 0}$ is Cauchy in $L_{2,R}$. It is enough to show
\be 
\label{E455} E\left( \left\| \sum_{i=n}^{\infty} DF_i \otimes e_i\right\|^2_{\cH \otimes \cH}\right) \xrightarrow[n \to \infty]{} 0.
\ee
According to Proposition \ref{l4.55ter} there is $\phi_{i} \in L_R$ such that $F_i= \iin  \phi_i dX$. Proposition \ref{PTX} says that $DF_i= \phi_i$. Consequently the left-hand side of (\ref{E455}) equals
\bea
E\left(  \sum_{i=n}^{\infty}\|DF_i \|^2_{\cH}\right) &=& E\left( \sum_{i=n}^{\infty}\|\phi_i\|^2_{\cH}\right) \nn \\
&=& \sum_{i=1}^n E(F_i^2) \xrightarrow[n \to \infty]{} 0, 
\eea
where the last equality is explained by Proposition \ref{cd3.9}.
\newline It remains to show that
\[
D_{t_1}F_{t_2} = h(t_1,t_2),
\]
with $h(t_1,t_2) = 1_{[0,t_1 \wedge T]}(t_2)$. We observe that
\[
DX^n = \sum_{i=1}^n e_i \otimes DF_i = \sum_{i=1}^n e_i \otimes \phi_i,
\]
so that
\bea
\|DX^n -h \|^2_R &=& \iin R(dt,\infty) \left\|\sum_{i=1}^n e_i(t)\phi_i-1_{[0,t\wedge T]} \right\|^2_R \nn \\
\label{ERA1}\\
&+& \int_{\RR} (-d\m)(t_1,t_2) \left\|\sum_{i=1}^n e_i(t_1)\phi_i - 1_{[0,t_1\wedge T]}- e_i(t_2)\phi_i + 1_{[0,t_2\wedge T]}\right\|^2_R. \nn
\eea
We have
\bea
&&\left\|\sum_{i=1}^n e_i(t)\phi_i - 1_{[0,t\wedge T]} \right\|^2_R = E \left\{\sum_{i=1}^n \iin(e_i(t)\phi_i - 1_{[0,t\wedge T]})dX \right\}^2 \nn \\
&&= \sum_{i,j=1}^n E \left(e_i(t) e_j(t) \iin \phi_i dX \iin \phi_j dX \right) \nn\\
\label{ERA2} \\
&-& 2 \sum_{i=1}^n e_i(t) E \left( \iin \phi_i dX X_t \right)  + E(X_t^2) \nn\\
&&=  E \left(\sum_{i=1}^n e_i(t) \iin \phi_i dX -X_t \right)^2= E\left(X^n_t-X_t \right)^2. \nn
\eea
By a similar reasoning, it follows that
\bea
&&\left\|\sum_{i=1}^n (e_i(t_1)-e_i(t_2)) \phi_i- 1_{[0,t_1 \wedge T]}+1_{[0,t_2 \wedge T]} \right\|^2_R \nn \\
\label{ERA3}\\
&& = E\left((X_{t_1}^n -X_{t_1})-(X^n_{t_2}-X_{t_2}) \right)^2. \nn
\eea
Therefore coming back to \eqref{ERA1} and taking into account \eqref{ERA2} and \eqref{ERA3}, we have
\[
\|DX^n-h \|^2_R = \|X-X^n \|^2_R \xrightarrow[n\to \infty]{}0
\] 
because of \eqref{ERAO}. \qed
\begin{rem}\label{RLYR} Adapting slightly the proof of Proposition \ref{R4550}, under the same assumptions, we have $X_{\cdot +r} \in \DDd$, for $r \in \R$ small enough.
\end{rem}
%\par We state first a technical lemma.
%\begin{lem} VERIFIER SI ON UTILISE CE LEMME \label{LYR} We suppose Assumptions (A), (B), (C($\n$)). Let $Y \in \Dd(L_R)$ and 
%$\psi \in L_R$, then $\left<Y, \psi\right>_R \in \Dd$ and $D_t\left<Y,\psi\right>= \left<D_tY,\psi\right>_R$.
%\end{lem}
%\textbf{Proof}: If $Y \in Cyl(L_R)$, then the result is obvious. The general case follows approximating $Y$ by a sequence $(Y^N)$ in $Cyl(L_R)$ and applying Cauchy-Schwarz in $L_R$ and $L_{2,R}$. \qed

\begin{prop} \label{PYR} Let $f \in C^2_b$ and $Y \in \Dd(L_R)$ such that
\be
\sup_{t \leq T}\|DY_t\| \in L^{\infty} \label{E5000}.
\ee
Then $f(Y) \in \Dd(L_R)$ and
\be
\label{EPYR0} Df(Y)=f'(Y)DY
\ee
in the sense that
\[
D_{t_2}f(Y_{t_1}) \equiv f'(Y_{t_1})D_{t_2}Y_{t_1}.
\]
\end{prop}
\begin{cor} \label{PCYR} Under Assumptions (C), (D), (\ref{ELTR}), $X_t=X_T$ if $t \geq T$, we have 
\[
f(X) \in \Dd(L_R)
\]
and
\[
D_r f(X_t) = f'(X_t)1_{[0,t]}(r).
\]
\end{cor}
\textbf{Proof:} This is a consequence of Proposition \ref{PYR} and Proposition \ref{R4550}.
\begin{rem} \label{RYR} If $Y \in Cyl(L_R)$ of the form $\sum_{i=1}^m F^i \psi_i$, $\psi_i \in C^1_0$, $F^i \in Cyl$, $f \in C^2_b$ we have 
\[
D_{t_2}f(Y_{t_1}) = \sum_{i=1}^n f'(Y_{t_1})\psi_i(t_1)D_{t_2}F^i.
\]
It is obviously a.s. an element of $L_{2,R}$ since $DF^i \in L_R$ a.s. and $f'(Y)\psi_i \in L_R$ by Propositions \ref{p3.4000} and \ref{p3.401}.
\end{rem}
\textbf{Proof} (of Proposition \ref{PYR}): We proceed in five steps.
\par a) We suppose that $Y\in Cyl(L_R)$, $f \in C^{\infty}_b(\R)$. Complications come from the fact that $f(Y)$ does not necessarily belong to $Cyl(L_R)$. Let $\psi \in L_R$. We show that 
\[
\left<f(Y),\cdot\right> \in Cyl
\]
and
\be
D(\left<f(Y),\psi\right>) = \left<f'(Y)DY, \psi\right>. \label{EPYR1}
\ee
\par b) We make some general considerations about approximations.
\par c) We suppose that $f \in C^2_b(\R)$, $Y \in Cyl$. For $\psi \in L_R$, we show that $\left<f(Y),\psi\right> \in \Dd$ and (\ref{EPYR1}) holds.
\par d) We suppose that $Y \in \Dd(L_R)$, $f \in C^2_b(\R)$. For $\psi \in L_R$ we show that
\[
\left<f(Y),\psi\right> \in \Dd
\]
and (\ref{EPYR1}) holds.
\par e) We conclude the proof.
\par We will proceed now in details step by step.
%rozwiniecie stepow ktore z gubsza zostaly opisane powyzej
\par a) Let $F^1, \ldots, F^m \in Cyl$, $\psi_1, \ldots, \psi_m \in C^1_0$ such that $Y = \sum_{i=1}^m F^i\psi_i$. Since $f \in C_b^{\infty}$, using the definition of inner product on $L_R$ and the definition of Malliavin derivative on $Cyl$, it follows that $\left<f(Y),\psi\right>_R \in Cyl$ and
\[
D\left( \left<f(Y), \psi\right>\right) = \sum_{i=1}^m DF^i\left<f'(Y)\psi_i,\psi\right>_R.
\]
This coincides with
\[
\left< f'(Y_{\cdot})DY_{\cdot}, \psi\right>_R
\]
taking into account Remark \ref{RYR}.
\par b) Consider the case $f \in C^2_b(\R)$. We regularize setting
\[
f_{\eps}(y)= \int_{\R}dz f(y + \eps z) \r(z),
\]
where $\r(y)= \frac{1}{\sqrt{2 \pi}}e^{-\frac{y^2}{2}}$. We denote $C_f= \|f' \|_{\infty}$, so 
\[
\sup_{y}|f_{\eps}(y)-f(y)| \leq \eps C_f \int_{\R}|z|\r(z)dz = \eps C_f \sqrt{\frac{2}{\pi}}.
\]
We observe that for every $\eps > 0$
\be
|f_{\eps}(y_1)-f_{\eps}(y_2)| \leq C_f|y_1 -y_2|.  \label{E4010}
\ee
Let $Y \in \Dd(L_R)$ fulfilling (\ref{E5000}). In this framework, we prove the following results
\be
\label{E3998} E \left( \|f(Y) \|^2_{R} \right) <\infty,
\ee
\be
\label{E3999} E\left(\|f'(Y)DY\|^2_{2,R} \right)<\infty.
\ee
\be
\label{E4001} E \left( \| (f-f_{\eps})(Y)\|^2_{R}\right) \xrightarrow[\eps \to 0]{} 0,
\ee
\be
\label{E4002} E\left( \left\| \left[ f'_{\eps}(Y) - f'(Y)\right]DY \right\|^2_{2,R} \right) \xrightarrow[\eps \to 0]{} 0. 
\ee
Indeed (\ref{E3998}) and (\ref{E3999}) follow by similar arguments as for (\ref{E4001}) and (\ref{E4002}). We only prove the two latter formulae.
\bea
E\left( \| (f_{\eps}-f)(Y)\|^2_R\right) &=& E\left(\iin |R|(dt,\infty)(f_{\eps}-f)^2(Y_t) \right) \nn \\
&& \ \ \ \ +E \left(\int_{\RR}d|\m|(s_1,s_2)\left[ (f_{\eps}-f)(Y_{s_1})-(f_{\eps}-f)(Y_{s_2})\right]^2 \right). \nn
\eea
(\ref{E4010}) implies that
\be
\label{E4011}|(f_{\eps}-f)(Y_{s_1})-(f_{\eps}-f)(Y_{s_2})| \leq 2 C_f|Y_{s_1}-Y_{s_2}|.
\ee
Since $f_{\eps} \lra f$ pointwise when $\eps \to 0$ and using Lebesgue's dominated convergence theorem, (\ref{E4001}) follows. 
\par Concerning (\ref{E4002}), using similar arguments and the fact that $f''$ is bounded, we obtain
\[
E \left( \| f'_{\eps}(Y)-f'(Y)\|_{R}^2\right) \xrightarrow[n \to \infty]{} 0.
\] 
So
\bea
E\left( \|(f_{\eps}'(Y)-f'(Y))DY \|^2_{{2,R}}\right)&\leq& E\left(\|(f_{\eps}-f)'(Y)\|D_{\cdot}Y \|_R \|^2_{R} \right) \nn \\
&& \ \ \ \ \ \ = I_1(\eps)+ I_2(\eps)+I_3(\eps) , \nn
\eea
where
\bea
I_1(\eps) &=& E \left( \iin |R|(dt,\infty)[(f_{\eps}-f)'(Y_t) ]^2 \| DY_t\|^2\right), \nn \\
I_2(\eps) & = & E \left(\int_{\RR} d|\m|(t_1,t_2) [(f_{\eps}-f)'(Y_{t_1})-(f_{\eps}-f)'(Y_{t_2}) ]^2 \| D Y_{t_1}\|^2_R \right), \nn \\
I_3(\eps) &=& E \left( \int_{\RR}d|\m|(t_1,t_2)(f_{\eps}-f)'(Y_{t_2})\left( \|DY_{t_1} \|_{R}-\|DY_{t_2} \|_R\right)^2\right). \nn
\eea
All the integrands converge a.s. and for any $(t_1,t_2)$ when $\eps \to 0$. We apply (\ref{E4011}) replacing $f_{\eps}, f$ with $f_{\eps}', f'$. The fact that $\sup_{t \leq T}\| DY_t\|_R \in L^2$, Lebesgue's dominated convergence theorem and Cauchy-Schwarz show that $I_i(\eps) \lra 0$, $i=1,2,3$.
\par c) We go on with the proof. If $Y \in Cyl$ clearly $Y \in \Dd(L_R)$ and (\ref{E5000}) is verified. Using (\ref{EPYR1}), it remains to show
\be
\label{E5001} E \left(\left<(f -f_{\eps})(Y), \psi\right>_R^2 \right) \xrightarrow[\eps \to 0]{} 0,
\ee
\be
\label{E5002} E \left( \left<(f'_{\eps}(Y)-f'_{\d}(Y))DY, \psi\right>^2_{{2,R}}\right) \xrightarrow[\eps,\d \to 0]{} 0.
\ee
The left-hand side of (\ref{E5001}) is bounded by
\[
\|\psi \|^2_R E(\| (f -f_{\eps})(Y)\|^2_R).
\]
because of Cauchy-Schwarz.
This together with (\ref{E4001}) implies \eqref{E5001}. (\ref{E5002}) holds again because of Cauchy-Schwarz and (\ref{E4002}). 
\par d) 
We first observe that $f(Y) \in L_R$ a.s. by Proposition \ref{p3.4000}. Let $Y \in \Ddd$ and a sequence $(Y^n)$ in $Cyl(L_R)$ such that
\[
\lim_{n \to \infty} E \left( \left\|Y-Y^n \right\|^2_{2,R}\right).
\]
We have
\be
E \left( \| f(Y^{n})-f(Y)\|^2_R\right) \leq \| f'\|^2_{\infty} E(\| Y^{n}-Y\|^2_R) \xrightarrow[n \to \infty]{} 0 \label{E6000}
\ee
and
\be
E \left(  \| f'(Y^{n})DY^{n}-f'(Y)DY\|^2_{2,R} \right) \xrightarrow[n \to \infty]{} 0. \label{E6001}
\ee
Then
\bea
\| f'(Y^{n})DY^{n})-f'(Y)DY\|^2_{2,R} &\leq&  \| f'(Y^{n})(DY^{n}-DY)\|^2_{2,R} + \|(f'(Y^{n})-f(Y))DY \|^2_{{2,R}} \nn\\
&\leq& \|f' \|^2_{\infty}\|DY^{n}-DY \|^2_{{2,R}} + \|(f'(Y^{n})-f'(Y))DY \|^2_{2,R}. \nn
\eea
The first term goes to zero since $Y^{n} \lra Y$ in $\Dd(L_R)$. The second one converges because $f''$ is bounded, using Lebesgue's dominated convergence theorem. This shows the validity of (\ref{E6000}) and (\ref{E6001}).
\par The next difficulty consists in showing that $U:=\left<f(Y), \psi\right> \in \Dd$ if $\psi \in L_R$. This will be the case approximating it through $U^{n}$, where
\[
U^{n} = \left<f(Y^{n}), \psi\right>_{R}. 
\]
Indeed, by item c) we have $U^n \in \Dd$ and taking into account Proposition \ref{AfterdefD} it remains to show that
\begin{description}
	\item[i)] $E((U^{n}-U)^2) \xrightarrow[n \to \infty]{} 0$,
	\item[ii)] $E\left(\|DU^{n}-DU^{m} \|^2_{R} \right)\xrightarrow[n,m \to \infty]{} 0.$ 
\end{description}
Concerning i) we can easily obtain
\[
E(U^{n}-U)^2 \leq \| \psi\|^2_R E( \| f(Y^{n})-f(Y)\|^2_R).
\]
This converges to zero because of (\ref{E6000}). 
As far as ii) is concerned, we can prove that
\be \label{FA55} \lim_{n \to \infty} E\|DU^{n}-\left<f'(Y)DY,\psi\right>_R \|^2_R=0.
\ee
Indeed, by item c) and \eqref{EPYR1}
\[
DU^{n} = \left<f'(Y^{n})D_{\cdot}Y, \psi\right>,
\]
so the left-hand side of (\ref{FA55}) gives
\bea
&&E \left(\left<f'(Y^{n})DY^{n}-f'(Y)DY, \psi\right> \right)_R^2 \nn \\
&& \ \ \ \ \leq \| \psi\|^2_R E \left(\| f'(Y^{n})DY^{n}-f'(Y)DY)\|^2_{{2,R}} \right) \xrightarrow[n \to \infty]{} 0  \nn
\eea
because of (\ref{E6001}). This concludes the proof of d).
\par e) Since $L_R$ is a separable Hilbert space, we consider an orthonormal basis $(e_n)^{\infty}_{n=0}$. We can expand a.s.
\[
f(Y) = \lim_{N \to \infty} I_N(f(Y)),
\]
where
\[
I_N(f(Y))= \sum_{n=0}^N\left<f(Y),e_n\right>e_n
\]
and the convergence holds in $L_R$. According to d) $\left<f(Y),e_n\right>\in \Dd$ and so $I_N(f(Y)) \in \Dd(L_R)$. It remains to show
\be
E\left( \| f(Y)-I_N(f(Y))\|^2_R \right) \xrightarrow[N \to \infty]{} 0, \label{E7000}
\ee 
\be \label{E7001} \|DI_N(f(Y))-DI_{M}(f(Y)) \|^2_{2,R} \xrightarrow[N,M \to \infty]{} 0. 
\ee
(\ref{E7000}) follows using Parseval's and Lebesgue's dominated convergence. Indeed
\bea
&&\| f(Y) - I_N(f(Y))\|^2_R = \sum_{n=N+1}^{\infty}  \left<f(Y),e_n\right>^2 \nn \\
&& \ \ \ \ \leq \sum_{n=0}^{\infty} \left<f(Y),e_n\right>^2 = \| f(Y)\|^2_R.\nn
\eea
$\| f(Y)\|^2_R$ is integrable because of (\ref{E3998}). Concerning (\ref{E7001}), taking $M>N$, we observe that
\[
DI_N(f(Y))-DI_{M}(f(Y)) = \sum_{n=N+1}^M \sum_{m=0}^{\infty}\left<f'(Y) DY, e_n\otimes e_m\right> e_n \otimes e_m,
\]
so by Parseval's in $L_{2,R}$ we have
\be
\label{E7010} \| DI_N(f(Y))-DI_M(f(Y))\|^2_{2,R} = \sum_{n=N+1}^M \sum_{m=0}^{\infty} \left<f'(Y)DY, e_n\otimes e_m\right>^2.
\ee
Now previous quantity converges a.s. to zero when $N,M \lra \infty$. Moreover (\ref{E7010}) is bounded by
\[
\| f'(Y)DY\|^2_{2,R}.
\]
Lebesgue's dominated convergence theorem finally implies (\ref{E7001}).  \qed

%\begin{rem} \label{RPYR} EN PRINCIPE FAUX. EN TENIR COMPTE POUR LES APPLICATIONS. Assumption (\ref{E5000}) and $f \in C^2_b$ in the statement of Proposition \ref{PYR} is restrictive only because we assume $Y \in \Dd(L_R)$. Following the schema of the previous proof, it is possible to obtain the conclusion (\ref{EPYR0}), if $Y=X$ and $f \in C^1$ with polynomial growth.
%\end{rem}
An easier but similar result to Proposition \ref{PYR} is the following
\begin{prop} \label{PYRRA} Let $Z$ be a random variable in $\Dd$, $f \in C^2_b$, $DZ\in L^{\infty}$. Then $f(Z) \in \Dd$ and
\[
Df(Z)=f'(Z)DZ.
\]
\end{prop}
\textbf{Proof:} It follows by similar, but simpler arguments than those of Proposition \ref{PYR}. \qed
\par Let $Y$ be a stochastic process such that $Y_t \in \Dd$ $\forall t \in \R_+$. Let $a: \RR \lra \R$ be Borel integrable function. We look for conditions on $a$ so that the process
\[
Z_t = \iin a(s,t) Y_s ds
\]
belongs to $\Dd(L_R)$. A partial answer is given below. We first proceed formally. If it exists, the Malliavin derivative is given by $D_rZ_t = Z_1(r,t)$ where
\[
Z_1(r,t) = \iin a(s,t)D_rY_s ds. 
\]
We need now another technical lemma.
\begin{lem} \label{LD101} Let $(d\r_t)$ be a $\sigma-$finite signed Borel measure on $\R_+$. Let $(Y_t)$ be a stochastic process fulfilling the following properties
\begin{description}
	\item[i)] For every $t \geq 0$ $Y_t \in \Dd$.
	\item[ii)] $t \longmapsto Y_t$ is continuous and bounded on $supp \ d\r_t$ in $\Dd$. In particular $t \longmapsto Y_t$ is continuous and bounded on $supp \ d\r_t$ in $L^2$ and  $t \longmapsto D_{\cdot}Y_t$ is continuous and bounded on $supp \  d\r_t$ in $L^2(\Om; L_R)$.
\end{description}
Let $g \in L^2(d\r_t)$. Then
\be
\label{ELD101} \iin g(t) Y_t d\r_t \in \Dd
\ee
and 
\be
\label{ELD102} D_r \left(\iin g(t) Y_td\r_t \right) = \iin g(t)D_rY_t d\r_t.
\ee
\end{lem}
\textbf{Proof:} We denote $t_i^n = i2^{-n}$, $i= 0, \ldots, n2^n$. We set
\[
\zeta^n = \sum_{i=1}^{n2^n} \int_{t_{i-1}^n}^{t_i^n}  Y_{t_i^n} g(s)d\r_s = \iin g(s) Y_s^n d\r_s,
\]
where
\[
	Y_s^n = \left\{
\begin{array}{clrr}
Y_{t_i^n} &\ \ \textrm{ if } s\in ]t_{i-1}^n,t_i^n], s\leq n, \\
0 &\ \ \textrm{ if } s > n.
\end{array} \right.
\]
It follows
\be
\label{ELD103} E\left( \iin (Y_s -Y_s^n)^2d\r_s\right) \xrightarrow[n \to \infty]{} 0,
\ee
since $t \longmapsto Y_t$ is continuous in $\LOm$. By Cauchy-Schwarz, it follows that $\zeta^n \lra \iin g(s) Y_s d\r_s$ belongs to $\LOm$. Since $\zeta^n$ is a linear combination of random variables issued from process $Y$, then $\zeta^n \in \Dd$ and 
\be
\label{ELD104} D_r\zeta^n = \iin g(t)D_r Y_t^n d\r_t, \ r \geq 0.
\ee
Since $t \longmapsto DY_t$ is continuous in $L^2(\Om;L_R)$, it follows
\be
\label{ELD105} E\left( \iin \|D_{\cdot}Y_t^n-D_{\cdot}Y_t \|^2_R d\r_t\right) \xrightarrow[n \to \infty]{} 0.
\ee
Again by Cauchy-Schwarz inequality it follows that
\[
D_r\zeta^n \lra \iin D_rY_t g(t) d\r_t, \ r \geq 0
\]
in $L^2(\Om;L_R)$. By Proposition \ref{AfterdefD}, the conclusions (\ref{ELD101}) and (\ref{ELD102}) hold. \qed
\begin{prop} \label{TPD101} Let $(d\r_t)$ be a finite Borel measure on $\R_+$, $a: \RR \lra \R$ be a Borel function, $(Y_t)$ be a stochastic process. We suppose the following.
\begin{description}
	\item[i)] $a(s,\cdot) \in L_R$ for $d\r_s$ a.e.
	\item[ii)] $\iin \|a(s,\cdot) \|^2_R d \r_s <\infty$ .
	\item[iii)] Assumptions i), ii) on $Y$ stated in Lemma \ref{LD101} hold.
\end{description}
We suppose that Assumptions (C), (D) hold. Then the process
\[
Z_t= \iin a(s,t)Y_s \ d\r_s, \quad \ t \geq 0,
\]
belongs to
$\Dd(L_R)$ and 
\be \label{EPD102} D_r Z_t = \iin a(s,t)D_r Y_s\  d\r_s, \quad \ r \geq 0.
\ee
\end{prop}
\textbf{Proof:} According to Proposition \ref{p376} there is an orthonormal basis $(e_n)$ of $L_R$ included in $C_0^1(\R)$. For a.e. $d\r_s$, we have
\[
a(s,\cdot) = \lim_{n \to \infty} a^n (s,\cdot) \  d\r_s  \ a.e.,
\]
where for $n \geq 1$
\[
a^n(s,\cdot)= \sum_{m=0}^n\left<a(s,\cdot),e_m\right>_R e_m.
\]
Moreover by Parseval's, for a.e. $d\r_s$
\be
\| a(s,\cdot)\|^2_R = \sum_{m=0}^{\infty} \left<a(s,\cdot),e_m\right>^2_R. \label{EEPSV} 
\ee
Let $m\geq 0$. According to hypothesis ii) and Cauchy-Schwarz it follows that $g_m(t):=\left<a(s,\cdot),e_m\right>_R$ belongs to $L^2(d\r_s)$. By Lemma \ref{LD101}, we obtain 
\[
\iin \left<a(s,\cdot),e_m\right>_R Y_s d\r_s \in \Dd
\]
and
\[
D_r\left( \iin \left<a(s,\cdot),e_m)\right>_R Y_s d\r_s\right)=\iin \left<a(s,\cdot),e_m)\right>_R D_r Y_s d\r_s.
\]
We denote
\[
Z_t^m = \iin a^m(s,t)Y_s d\r_s.
\] 
By linearity and Remark \ref{R4530} ii) $Z^m \in \Dd(L_R)$ and
\[
DZ_t^m = \iin a^m(s,t)D Y_s d\r_s, \ \ t \geq 0.
\]
Using Remark \ref{R4530} i), it will be enough to show
\begin{description}
	\item[a)] $\lim_{m \to \infty} E(\|Z^m -Z \|^2_R)=0$,
	\item[b)] $\lim_{m\to \infty} E(\|DZ^m -\mathcal{Y} \|^2_{2,R})=0$, where
	\[
	\mathcal{Y}(r,t) = \iin a(s,t)D_r Y_s d\r_s.
	\]
\end{description}
a) First, we observe that $Z \in \tilde{L}_R$ a.s. because, avoiding some technical details, we have
\bea
&&\left\| \iin a(s,\cdot) Y_s d\r_s \right\|_R \leq const. \iin \|a(s,\cdot)Y_s \|_Rd\r_s \nn\\
&&= const. \iin \| a(s,\cdot)\|_R|Y_s| d\r_s \leq const. \sqrt{\iin \|a(s,\cdot) \|_R^2d\r_s \iin |Y_s|^2d\r_s}. \nn
\eea
Again Cauchy-Schwarz and condition ii), imply that
\[
E\left(\left\|\iin a(s,\cdot)Y_s d\r_s \right\|^2_R \right)<\infty.
\]
Taking into account (\ref{EEPSV}) and Lebesgue's dominated convergence theorem we can show that
\bea
&&E(\|Z-Z^m \|^2) \nn\\
&&\label{E417septe}\\
&& \ \ \ \ \leq const. \left\{ E\left( \iin |Y_s|^2 d\r_s\right)\iin \| a-a^m(s,\cdot)\|^2_R d\r_s\right\} \xrightarrow[m \to \infty]{} 0.\nn
\eea
b) By similar arguments, we can show that
\[
\mathcal{Y} \in \tilde{L}_{2,R} \  a.s.
\]
and 
\[
E\left( \| \mathcal{Y}\|^2_{2,R}\right)<\infty.
\]
Moreover
\[
DZ^m-\mathcal{Y}= \iin(a^m-a)(s,t)D_rY_s d\r_s.
\]
Consequently, by similar arguments as in (\ref{E417septe}) it follows that
\[
E\left(\|DZ^m-\mathcal{Y} \|^2_{2,R} \right) \xrightarrow[m \to \infty]{} 0.
\]
This concludes the proof of Proposition \ref{TPD101}. \qed

An application of previous proposition is the following. It holds under Assumptions (C), (D).
\begin{prop} \label{p4.1517A} Let $Y$ be a process, continuous in $L^2$, such that $Y_t \in \Dd$ for every $t \geq 0$ and $t \longmapsto D_{\cdot}Y_t$ is continuous in $L^2(\Om;L_R)$. Let $\eps >0$ and denote
\[
Y_t^{\eps} = \int_{(t-\eps)^+}^{(t+\eps) \wedge T} Y_s ds.
\]
Then $Y^{\eps} \in \Dd(L_R)$.
\end{prop}
\textbf{Proof:} In view of applying Proposition \ref{TPD101}, we set $\r(t)=t \wedge T,$
\[
a(s,t)= 1_{]t-\eps,t+\eps]\cap ]0,T]}(s)1_{[0,T]}(t),
\]
which also gives
\[
a(s,t)= 1_{[s-\eps,s+\eps[\cap [0,T+\eps[}(t)1_{[0,T]}(s).
\]
We have $Y_t^\eps = \iin a(s,t)Y_s d\r_s$. According to Assumption (D) and Corollary \ref{c3.72} $a(s,\cdot) \in L_R$, for every $s \geq 0$, Assumption i) of Proposition \ref{TPD101} is verified. Again by Corollary \ref{c3.72}
\[
\| a(s,\cdot)\|^2_R = Var(X_{(T\wedge s)+\eps}-X_{(s-\eps)^+}) \leq 2 Var(X_{(T\wedge s)+\eps})+2Var(X_{(s-\eps)^+}).
\]
Since $X$ is continuous in $L^2$, $s \longmapsto \|a(s,\cdot) \|_R$ is bounded and Assumption ii) of Proposition \ref{TPD101} is verified.
\par Point iii) of Proposition \ref{TPD101} follows by the continuity assumption on $Y$ and $DY$ and because $\r$ has compact support. \qed
\par In the sequel, we will apply Proposition \ref{p4.1517A} to $Y=g(X)$ with $g$ having polynomial growth. 
\par The lemma below allows to improve slightly the statement of Proposition \ref{TPD101}.
\begin{lem} \label{PD100} Let $(Y_t)$ (resp. $(Y^n_t))$ be a process (resp. a sequence of processes) such that $Y_t, Y_t^n \in \Dd$, $\forall t \in \R_+$ and
\bea
\label{ED100} &&E\left(\iin Y_t^2 d\r_t + \iin \|D_{\cdot}Y_t \|^2_R d\r_t \right) < \infty \\
\label{ED101} && E \left( \iin (Y_t-Y_t^n)^2d\r_t+ \iin \|D_{\cdot}(Y_t-Y_t^n) \|^2_R d\r_t\right) \xrightarrow[n \to \infty]{} 0.
\eea
Let $a: \RR \lra \R$ be a Borel function such that
\be
\label{ED102} \iin d \r_s\| a(s,\cdot)\|^2_R < \infty.
\ee
We set
\bea
Z_t &=& \iin a(s,t) Y_s d\r_s, \nn \\
Z_1(r,t) &= &\iin a(s,t)D_rY_s d\r_s, \nn \\
Z_t^n& = &\iin a(s,t)Y_s^n d\r_s, \nn \\
Z_1^n(r,t) &=& \iin a(s,t) D_r Y_s^n d\r_s. \nn
\eea
We have the following properties:
\begin{description}
	\item[i)] $Z \in \tilde{L}_R$, $Z_1 \in \tilde{L}_{2,R}$ and 
	\[
	E(\| Z\|_R^2 + \| Z_1\|^2_{2,R}) <\infty.
	\]
	\item[ii)] 
	\[ \lim_{n \to \infty} E\left( \|Z^n-Z \|^2_R + \|Z_1^n -Z_1 \|^2_R\right) =0.
	\]
\end{description}
\end{lem}
\textbf{Proof:} We only prove point i) since the other follows similarly.
\bea
\| Z\|_R^2&=& I_1 +I_2,\nn\\
\| Z_1\|^2_{2,R} &=& I_3 +I_4, \nn
\eea
where
\bea
I_1&=& \iin |R|(dt,\infty)\left(\iin a(s,t)Y_s d\r_s \right)^2 ,\nn\\
I_2&=& \int_{\RR} d|\m|(t_1,t_2)\left(\iin(a(s,t_1)-a(s,t_2))Y_s d\r_s \right)^2,\nn \\
I_3&=& \iin |R|(dt,\infty) \left\|\iin a(s,t)D_{\cdot}Y_s d\r_s \right\|^2_R, \nn \\
I_4&=& \int_{\RR}d|\m|(t_1,t_2)\left\| \iin (a(s,t_1)-a(s,t_2))D_{\cdot}Y_s d\r_s\right\|^2. \nn
\eea
Cauchy-Schwarz implies that
\[
I_1 \leq \iin |R|(dt,\infty) \left( \iin a^2(s,t)d\r_s\right) \left( \iin Y_u^2 d\r_u \right),
\]
\[
I_2 \leq \int_{\RR} d|\m|(t_1,t_2) \iin(a(s,t_1)-a(s,t_2))^2 d\r_s \iin Y_s^2 d\r_s.
\]
Consequently
\[
E(I_1 + I_2) \leq E\left( \iin Y_s^2 d\r_s\right) \iin \| a(s,\cdot)\|_R^2 d\r_s <\infty.
\]
On the other hand, Bochner integration theory implies
\[
I_3+I_4 \leq \iin \|D_{\cdot}Y_s \|^2_R d\r_s \iin d\r_s \|a(s,\cdot) \|^2_R.
\]
Taking the expectation it follows
\[
E(I_3+I_4) <\infty. 
\] \qed
\newline Next result allows to relax the boundedness property on $\| Y\|_R$ and $\|D_{\cdot} Y \|_R$ required in Proposition \ref{TPD101}.
\begin{cor} \label{CPD102} Let $(d\r_t)$ be a finite measure on $\R_+$, $a: \RR \lra \R$ be a Borel function. Let $(Y_t)$ be a stochastic process continuous in $\Dd$ such that
\[
E \left( \iin Y_t^2 d\r_t + \iin \|D_{\cdot}Y_t\|^2_R d\r_t\right)<\infty.
\]
We only suppose i), ii) as in Proposition \ref{TPD101}. Then the same conclusion as therein holds.
\end{cor}
\textbf{Proof:} We use Lemma \ref{PD100} and we approximate $Y$ by $Y^m$, where $Y^m = \phi^m(Y)$ for $\phi=\phi^m:\R \lra \R$ smooth, with
\[
\phi(y) = \left\{
\begin{array}{clrr}
y &\ \ , |y|\leq m\\
0 &\ \ , |y|>m+1.
\end{array} \right.
\]
Clearly
\be \label{E41522} E\left( \iin (Y_t-Y^m_t)^2 d\r_t\right) \xrightarrow[m \to \infty]{} 0 
\ee
by Lebesgue's dominated convergence theorem. Moreover Proposition \ref{PYR} implies
\[
DY^m = \phi'(Y)DY,
\]
and of course $\phi$ is smooth, bounded such that
\[
\phi'(y) = \left\{
\begin{array}{clrr}
1 &\ \ , |y|\leq m\\
0 &\ \ , |y|>m+1.
\end{array} \right.
\]
Therefore we have again
\[
DY^m-DY = DY(1- \phi'(Y))
\]
and by similar arguments, we obtain
\be \label{E41523} E \left(\iin \|D_{\cdot}(Y^m_t -Y_t )\|^2_R d\r_t \right) \xrightarrow[m \to \infty]{}0.
\ee
Finally (\ref{E41522}), (\ref{E41523}) and Lemma \ref{PD100} allow to conclude. \qed

\section{Skorohod integrals} \label{s9}

\setcounter{equation}{0}

\subsection{Generalities} \label{s10.1}

We suppose again Assumptions (A), (B), (C) by default. Similarly as in \cite{CN} and \cite{MV}, we will define two natural domains for the divergence operators, i.e. Skorohod integral, which will be in some sense the dual map of the related Malliavin derivative. We denote by $L^2(\Om;L_R)$ the set of stochastic processes $(u_t)_{t \in [0,T]}$ verifying $E(\|u \|)^2_R < \infty$. We say that $u \in L^2(\Om; L_R)$ belongs to $Dom(\d)$ if there is a zero-mean square integrable random variable $Z$ such that
\be
\label{dual} E(FZ) = E(\left<DF,u\right>_{\cH})
\ee
for every $F \in Cyl$. In other words we have
\be
\label{E431} E(FZ) = E\left( \iin R(ds,\infty)D_sFu_s \right)- E\left( \int_{\RR} \m(ds_1,ds_2)(D_{s_1}F -D_{s_2}F)(u_{s_1}-u_{s_2})\right)
\ee
\newline for every $F \in Cyl$. Using the Riesz theorem, we can see that $u\in Dom(\d)$ if and only if the map
\[
F \lra E(\left<DF,u\right>_{\cH})
\]
is continuous linear form with respect to $\|\cdot \|_{\LOm}$. Since $Cyl$ is dense in $\LOm$, $Z$ is uniquely characterized. We will call $Z = \iin u \d X$ the \textbf{Skorohod integral} of $u$ towards $X$.
\par We continue investigating general properties of Skorohod integral.
\begin{df} \label{deef} If $u\jeden \in Dom(\d)$,  for any $t \geq 0$, then we define
\[
\int_0^t u_s \d X_s := \iin u_s \jeden \d X_s.
\]
\end{df}
A consequence of the duality formula defining Skorohod integral appears below.
\begin{rem} \label{r4.5.1}
If (\ref{dual}) holds, then it will be valid by density for every 
$F \in \DD$. 
\end{rem}
\begin{prop} \label{p4.5.1} Let $u \in Dom(\d)$, $F \in \DD$. Suppose $F\iin u_s \d X_s \in \LOm$. Then $Fu \in Dom(\d)$ and
\[
\iin Fu_s \d X_s = F\iin u_s \d X_s - <DF,u>_{\mathcal{H}}.
\]
\end{prop}
\textbf{Proof}: The proof is very similar to the one of Proposition 6.4 in \cite{KRT}. Let $F_0 \in Cyl$. We need to show 
\[
E\left( F_0 \left\{ F\iin u_s \d X_s- <DF,u>_{\cH}\right\} \right) = E\left( <DF_0,Fu>_{\cH}\right).
\]
We proceed using Lemma \ref{l425}, which says that $F_0F \in \DD$ and \eqref{E424} holds (with $G=F_0$), together with Remark \ref{r4.5.1}, which extends the duality relation. \qed
\par We state now Fubini's theorem, which allows to interchange Skorohod and measure theory integrals. When $X$ has a covariance measure structure, this was established in \cite{KRT}, Proposition 6.5. When $X$ is a Brownian motion the result is stated in \cite{Nualart}.

\begin{prop} \label{p4.5.2} Let $(G, \mathcal{G},\l)$ be a $\s$-finite measure space. Let $u: G\times \R_+\times \Om \lra \R$ be a measurable random field with the following properties.
\begin{description}
	\item[i)] For every $x \in G, u(x,\cdot)\in Dom(\d)$.
	\item[ii)]
	\[
	E\left( \int_G d\l (x)\|u(x,\cdot)\|^2_{R}\right) < \infty.
	\] 
	\item[iii)] There is a measurable version in $\Om \times G$ of the random field
	\[
	\left( \iin u(x,t)\d X_t\right)_{x \in G}.
	\]
	\item[iv)] It holds that
	\[
	\int_G d \l (x) E\left( \iin u(x,t) \d X_t \right)^2 <\infty.
	\]
	Then $\int_G d\l (x) u(x,\cdot) \in Dom(\d)$ and
	\[
	\iin \left( \int_G d\l(x)u(x,\cdot) \right) \d X_t = \int_G d\l (x)\left(\iin u(x,t) \d X_t \right).
	\]
\end{description}
\end{prop}
\textbf{Proof}: We will prove two following properties,
\begin{description}
	\item[a)] $\int_G d\l (x) u(x,\cdot) \in L^2(\Om,L_R)$
	\item[b)] For every $F \in Cyl$ we have 
	\be \label{F451}
	E\left( F\left(\int_G d\l(x) \iin u(x,\cdot) \d X_t \right) \right) = E\left( \left< DF, \int_G d\l (x)u(x, \cdot)\right>_{\cH} \right).
	\ee
\end{description}
Without restriction to the generality we can suppose $\l$ to be a finite measure. Concerning a), Jensen's inequality implies
\bea E\left( \left\| \int_G d\l(x) u(x, \cdot)\right\|_R^2 \right) \leq \l(G) E\left(\int_G d\l(x) \left\| u(x, \cdot)\right\|_R^2 \right) = \l(G)\int_G d\l(x)\|u(x,\cdot) \|_R^2 <\infty \nn
\eea
because of ii). For part b), by classical Fubini's theorem, the left-hand side of (\ref{F451}) gives 
\bea
&& \int_G d\l (x) E \left(F \iin u(x,t) \d X_t \right) = \int_G d\l (x) E \left(\left<DF, u(x,\cdot)\right>_{\cH} \right)\nn \\
&& \label{F452} \\
&& = E\left( \int_G d \l(x)\left<DF,u(x,\cdot)\right>_{\cH}\right). \nn
\eea
This is possible because
\[
|\left<DF,u(x,\cdot)\right>_{\cH}| \leq \|DF \|_{\cH} \|u(x,\cdot) \|_R.
\]
(\ref{F452}) equals the right-hand side of (\ref{F451}) because
\bea
&&\int_G d\l(x)\left<DF,u(x,\cdot)\right>_{\cH}\nn \\
&& = \int_G d\l(x) \left( \iin D_s F u(x,s) R(ds,\infty) \right. \nn\\
&& \ \ \ \ \ \ \ \ \ \ \ \ \ \ \  \left. - \fr \int_{\RR}(D_{s_1}F -D_{s_2}F)(u(x,s_1)-u(x,s_2))d\m(s_1,s_2) \right) \nn\\
&& = \iin R(ds,\infty) D_sF\int_G u(x,s) d\l(x) \nn \\
&& \ \ \ \ \ \ \ \ \ \ \ \ \ \ \ - \fr \int_{\RR}(D_{s_1}F-D_{s_2}F) \int_G (u(x,s_1)-u(x,s_2)) d\l(x) d\m(s_1,s_2). \nn
\eea
Last equality is possible by means of classical Fubini's theorem and assumption ii). This equals
\[
\left<DF,\int_G d\l(x) u(x,\cdot) \right>
\]
and the proof is concluded. \qed

\subsection{Malliavin calculus and Hermite polynomials} \label{s3.3bis}

We introduce here shortly Hermite polynomials. For more information, refer to \cite{Nualart}, Section 1.1.1. Those polynomials have the following properties. For every integer $n \geq 1$
\begin{description}
	\item[i)] $nH_n(x)=xH_{n-1}(x)-H_{n-2}(x)$,
  \item[ii)] $H'_n(x) = H_{n-1}(x)$,
	\item[iii)] $H_0(x)\equiv 1$, $H_{-1}(x)= 0$.
\end{description}
From Proposition \ref{PPol} the following result follows.

\begin{prop} \label{Pexph}
Let $h \in L_R$.
\begin{description}
	\item[i)] For any $n \geq 1$, $F:=H_n\left( \iin h dX\right) \in \Dd$ and
	\[
		D_tF = H_{n-1} \left(\iin h dX \right)h(t).
	\] 
	\item[ii)] $F = \exp\left( \int_0^{\infty} h dX\right)$. Then $F \in \Dd$ and
	\[
		D_tF = F h(t).
	\]
\end{description}
\end{prop}
We recall that $\{H_k, k\geq n \}$ constitute a basis of the linear span generated by $\{1, ,\ldots,x^n \}$. We denote by $\mathcal{E}_{Herm}$ the linear subspace of $\LOm$ constituted by all finite linear combinations of elements of the type $H_n\left(\iin \phi_n dX \right)$, $\phi_n \in C^1_0(\R_+)$, $n \in \mathbb{N}$.
\begin{prop} \label{C432}
	$\bar{\mathcal{E}}_{Herm}= L^2(\Om)$.
\end{prop}
\textbf{Proof:} We first observe that
\[
	\left\{ \exp\left( \iin h dX\right), h\in C^1_0(\R_+) \right\}
\]
is total in $\LOm$. The idea is to show that a random variable $F \in \LOm$, such that $E\left( F\exp \left(\iin h dX \right)\right)=0$ for every $h \in C^1_0(\R_+)$,  fulfills
\[
	E \left( F g\left(\iin h_1 dX, \ldots, \iin h_n dX \right)\right)=0,
\]
for every $h_1,\ldots, h_n \in C^1_0 (\R_+)$ and $g \in C_b^{\infty}(\R^n)$. This can be done adapting the proof of Lemma 1.1.2 of \cite{Nualart}.
\par Let us now consider $F \in L^2(\Om)$, $h\in C^1_0(\R_+)$ such that
\be
	\label{EHN1} E\left(FH_n\left(\iin h dX\right)\right)=0, \forall  \in \mathbb{N}.
\ee
It remains to show
\be
\label{EHN2} E\left(F\exp\left(\iin h dX \right) \right) =0.
\ee
By (\ref{EHN1}) it follows obviously that
\[
E\left(F \left(\iin h dX \right)^n \right) =0, \forall n \in \mathbb{N}
\]
and consequently (\ref{EHN2}) holds.\qed

We denote by $\mathcal{E}_n$ the linear span of $H_n\left( \iin \phi dX\right), \phi \in C^1_0(\R)$, $\|\phi \|_{\cH}=1$, and by $\mathcal{H}_n$ the closure of $\mathcal{E}_n$ in $\LOm$. We recall that all the considered Wiener integrals are Gaussian random variables. Adapting Theorem 1.1.1 and Lemma 1.1.1 of \cite{Nualart}, we obtain the following result.
\begin{prop} \label{p432}
\begin{enumerate}
	\item The spaces $(\cH_n)$ are orthogonal.
	\item $\LOm = \oplus_{n=0}^\infty \cH_n$. 
\end{enumerate}
\end{prop}
We discuss here some technical points related to Malliavin derivative and chaos spaces.
\begin{lem} \label{L4321}
	Let $n \geq 1$. The map $D: \mathcal{E}_n \subset \LOm \lra L^2(\Om; L_R)$ verifies the following. For any sequence $(F_k)$ in $\mathcal{E}_n$ converging to zero in $\LOm$, $(DF_k)$ is Cauchy in the sense that
	\[
		\lim_{k,l \to \infty} E\left(\| DF_{k} -DF_{l}\|^2_{\cH} \right) = 0.
	\]
\end{lem}
\textbf{Proof: } The result will follow if, for every $F \in \mathcal{E}_n$ we prove
\be
	E \left( \|DF \|^2_{\cH}\right) = n E(F^2).
\ee
Let $F = \sum_{k=1}^{m} H_{n} \left(\iin h_k dX \right)$, $h_k \in C^1_{0}$. Item i) of Proposition \ref{Pexph} and Lemma 1.1.1 of \cite{Nualart} imply that
\bea
E \left( \|DF \|^2_{\cH}\right) &=& \sum_{k,l =1}^m E \left(H_{n-1}\left( \iin h_k dx\right)H_{n-1} \left(\iin h_l dX\right)  \right) \left<h_l,h_k \right>_{\cH} \nn \\
&=& \sum_{k,l=1}^m \left<h_k,h_l \right>_{\cH} \frac{1}{(n-1)!} \left\{E \left(\iin h_k dX \iin h_l dX \right) \right\}^{n-1} \nn.
\eea
In fact $\iin h dX$ is a standard Gaussian random variable. This gives
\be
 \frac{1}{(n-1)!} \sum_{k,l =1}^{m} \left<h_k,h_l \right>_{\cH} \left<h_k,h_l \right>_{\cH}^{n-1} = n \frac{1}{n!} \sum_{k,l=1}^n \left<h_k,h_l \right>^n_{\cH} = n E(F^2) \nn
\ee
again by Lemma 1.1.1 of \cite{Nualart}.
\begin{cor} \label{R4321}
Let $n \geq 1$.
\begin{description}
	\item[i)] $\cH_n \subset \Dd$
	\item[ii)] If Assumption (D) is verified, then $D: \mathcal{H}_n \subset \LOm \lra L^2(\Om;L_R)$ is continuous.
	\item[iii)] Suppose that Assumption (D) is verified. For every $F \in \cH_n$, we have $\left< DF,h \right>_{\cH} \in \cH_{n-1}$, $\forall h \in L_R$.
\end{description}
\end{cor}
\textbf{Proof}: i) Let $F \in \cH_n$ and $(F_k)$ a sequence in $\mathcal{E}_n$ converging to $F$ in $L^2(\Om)$. By Lemma \ref{L4321} and Proposition \ref{AfterdefD}, $F \in \Dd$.
\newline ii) It is an obvious consequence of Lemma \ref{L4321}.
\newline iii) Let $h \in L_R$. By items i) and ii) $T_{h}: \cH_n \lra \LOm$ defined by $T_h(F) = \left<DF,h \right>_{\cH}$ is continuous. By Proposition \ref{Pexph} i) the image of $\mathcal{E}_n$ through $T_h$ is included in $\cH_{n-1}$. Since $\cH_{n-1}$ is a closed subspace of $\LOm$, the result follows. \qed
\begin{prop} \label{p4322}
We suppose the validity of Assumptions (C) and (D). Let $F \in \Dd$, $h \in L_R$. Then there is a sequence $(F_n)$ such that $F_n \in \cH_n$, $h \in L_R$,
\[
F= \sum_{n=0}^{\infty}F_n,
\]
\[
\left<DF,h\right>_{\cH} = \sum_{n=1}^{\infty}\left<DF_n,h\right>_{\cH}, \ \ \forall h \in L_R,
\]
where the convergence holds in $L^2(\Om)$.
\end{prop}
\textbf{Proof: } Let $h \in L_R$ and $F = \sum_{n=0}^{\infty}F_n$ according to Proposition \ref{p432}. By Corollary \ref{R4321}, iii), $\left<DF_{m+1},h \right>_{\cH}$ belongs to $\cH_m$; since $\left<DF,h\right>_{\cH}\in \LOm$, we need to show that for $m \geq 0$
\be
E \left(\left<DF, h\right>_{\cH}\mathcal{G}_m \right) = E \left( \left<DF_{m+1},h\right>_{\cH}\mathcal{G}_m\right) \label{EEE}
\ee
for every $\mathcal{G}_m \in \cH_m$. To prove (\ref{EEE}), taking into account Corollary \ref{R4321} i) and Remark \ref{rr4212}, we write
\bea 
E\left(\left<DF,h\right>_{\cH}\mathcal{G}_m \right) &=& E \left(F\left\{\mathcal{G}_m \iin h dX - \left<D \mathcal{G}_m,h\right>_{\cH}\right\} \right) \nn \\
&=& \sum_{n=0}^{\infty} E\left(F_n\left\{\mathcal{G}_m \iin h dX - \left<D \mathcal{G}_m,h\right>_{\cH}\right\} \right) \nn \\
&=& \sum_{n=0}^{\infty} E\left( \left<DF_n,h\right>_{\cH} \mathcal{G}_{m}\right) = E \left( \left<DF_{m+1},h\right>_{\cH} \mathcal{G}_m\right). \nn
\eea  \qed

\subsection{Generalized Skorohod integrals and Hermite polynomials} \label{s10.3}

\par We define now, implementing the idea of \cite{CN} and \cite{MV}, an extension of $Dom(\d)$ denoted by $Dom(\d)^*$. The idea is to use a similar duality relation to (\ref{E431}), but keeping in mind that an element $u$ of $\Dstar$ will not necessarily live in $L^2(\Om;L_R)$. We denote by $\mathcal{L}^2$ the space of processes $(u_t)_{t \geq 0 }$ with
\[
E\left( \iin u_s^2|R|(ds,\infty)+\iin |u_s|^2\bar{m}(ds)\right) < \infty,
\]
where $\bar{m}$ is the marginal measure of $|\bar{\m}|$.
\newline We will denote by $\mathcal{M}$ the linear space of Borel functions $f:\R_+ \lra \R$ such that
\be
\| f\|_{\mathcal{M}}:=\iin f^2(s)|R|(ds,\infty) + \int_{\RR} (f(s_1)-f(s_2))^2 |\bar{\mu}|(ds_1,ds_2) < \infty. \label{EMM}
\ee
\begin{rem} \label{r331}
\begin{enumerate}
	\item Obviously $\| \cdot\|_{\mathcal{M}} \leq \| \cdot \|_{R}$ and so $\Lr \subset \mathcal M$.
	\item $\mathcal M$ is complete (therefore it is a Hilbert space) because of Assumption (C). The argument is similar to the one used in proof of Proposition \ref{p3.5}.
	\item Expanding the second integral of \eqref{EMM} and applying Cauchy-Schwarz, we can show that $\mathcal{L}^2 \subset \mathcal{M}$.
\end{enumerate}
\end{rem}

\begin{df} \label{d431}
A process $u \in L^2(\Om;\mathcal M)$ is said to belong to $\Dstar$ if there is a square integrable r.v. $Z \in \LOm$ such that
\bea 
 E(FZ) &=& E \left( \iin R(ds,\infty)D_s F u_s \right.\nn \\
 &&\label{E434}\\ 
  &-& \left. \int_{\RR} \mu(ds_1,ds_2) ((D_{s_1}F-D_{s_2}F)(u_{s_1}-u_{s_2})) \right) \nn
\eea
for any $F = H_n\left( \iin h dX\right)$, $h \in C^1_0$, $n \geq 0$ .
\end{df} 
\begin{rem} \label{r431}
\begin{enumerate}
	\item If $u \in \Dstar$ then (\ref{E434}) holds by linearity for every $F \in \cE_{Herm}$.
	\item Since $\cE_{Herm}$ is dense in $\LOm$, $Z$ is uniquely determined. $Z$ will be called \textbf{the $*$-Skorohod integral} of $u$ with respect to $X$, it will be denoted $\iin u_s \d^*X_s$. 
	\item The right-hand side of (\ref{E434}) is well-defined for any $u \in L^2(\Om,\mathcal M)$ and 
	\newline $F = f\left(\iin \f_1dX, \ldots, \iin\f_n dX \right)$, $\f_1, \ldots, \f_n \in C_{0}^1(\R_{+})$, $f \in C^1_{pol}(\R^n)$, in particular for $F$ as in Definition \ref{d431}. We observe that in this case $F^i = \pa_i f(\iin \f_1 dX, \ldots, \iin \f_n dX)$, $1\leq i \leq n$, is a square integrable random variable.
	\par Indeed, we only discuss the second addend of the right-hand side of \eqref{E434}, since the first one is obviously meaningful. 
	\newline This is bounded by 
\bea
&&E \left(\int_{\RR}|\m|(ds_1,ds_2) |(D_{s_1}F -D_{s_2}F)(u_{s_1}-u_{s_2})| \right)\nn \\
&& = \sum_{^j=1}^n E \left(\int_{\RR}|\m|(ds_1,ds_2) |F^i(\f_j(s_1)-\f_j(s_2))(u_{s_1}-u_{s_2})| \right) \nn \\
&& \leq \sum_{j=1}^n \|\f'_j \|_{\infty} \left(\int_{\RR} |\bar{\m}|(ds_1,ds_2)|u_{s_1}-u_{s_2}|F^j\right) \nn \\
&& \leq \sum_{j=1}^n \|\f'_j \|_{\infty} \sqrt{ E (F^j)^2 E \left(\int_{\RR} |\bar{\m}|(ds_1,ds_2)|u_{s_1}-u_{s_2}| \right)^2} \nn \\
&& \leq \sum_{j=1}^n \|\f'_j \|_{\infty} \|F^j \|_{\LOm} \sqrt{|\bar{\m}|(\RR)} E\left(\int_{\RR}|\bar{\m}|(ds_1,ds_2) |u_{s_1} -u_{s_2}|^2 \right). \nn
\eea
If $u \in L^2(\Om; \mathcal{M})$, then previous quantity is finite.
\end{enumerate}
\end{rem}
\begin{prop} \label{pDDstar} $Dom(\delta) \subset \Dstar$. 
\end{prop}
\textbf{Proof: } Let $u \in Dom(\d)$, $Z = \iin u \d X$. First of all $u \in L^2(\Om; \mathcal{M})$ since $\|u\|_{\mathcal{M}}  \leq \|u\|_R$. For any $F \in Cyl$ we have
\be
\label{EDDstar} E(FZ) = E \left(\left<DF,u\right>_{\cH} \right).
\ee
Relation (\ref{EDDstar}) extends to the elements $F$ of the type $f\left( \iin h dX\right)$, $f \in C^1$ with subexponential derivative. For this, it is enough to provide the same type of approximation sequence as in the proof (item ii) ) of Proposition \ref{PPol}. If $(f^M)$ is such a sequence, setting $F^M = f^M \left(\iin h dX \right)$ and taking into account Proposition \ref{PPol}, we clearly obtain
\bea
&& \lim_{M\to \infty} E \left( \left<DF^M, u\right>_{\cH}\right) = E\left(\left<DF,u\right>_{\cH} \right),\nn\\
&& \lim_{M \to \infty} E(F^MZ)= E(FZ).
\eea
This implies in particular that (\ref{EDDstar}) holds for $F$ of the type $H_n\left(\iin h dX \right)$, where $n \geq 0$, $h \in C^1_0.$ For such an element, the right-hand side of (\ref{EDDstar}) coincides with the right-hand side of (\ref{E434}) and the result follows. \qed

\section{Itô formula in the very singular case} \label{s10}

\setcounter{equation}{0}

We suppose again Assumptions (A), (B), (C) by default. We start with a technical observation.
\begin{lem} \label{lI1} Let $(G_1,G_2)$ be a Gaussian vector such that $Var G_2 = 1$. Let $f \in C^1(\R)$ such that $f'$ is subexponential. Then 
\[
nE(f(G_1)H_n(G_2))=E(f'(G_1)H_{n-1}(G_2))Cov(G_1,G_2).
\]
\end{lem}
\begin{rem} \label{RI1} It follows in particular that $E(H_n(G_2))=0$, $\forall n\geq 1$.
\end{rem}
\textbf{Proof}: According to relation i) about Hermite polynomials, the left-hand side equals $I_1-I_2$,
where
\bea
I_1 &=& E(f(G_1)G_2H_{n-1}(G_2)\nn\\
I_2 &=& E(f(G_1)H_{n-2}(G_2)\nn.
\eea
According to Wick theorem, recalled briefly in Lemma \ref{Wick} below, $I_1$ gives.
\[
E(f'(G_1)H_{n-1}(G_2))Cov(G_1,G_2) + E(f(G_1)H'_{n-1}(G_2)).
\]
Using relation ii) about Hermite polynomials, we have
\[
E(f(G_1)H_{n-1}'(G_2) = E(f(G_1)H_{n-2}(G_2)) = I_2
\]
and the result follows. \qed
\par The lemma below was recalled and for instance proved in \cite{FGR}.
\begin{lem} \label{Wick}(Wick) Let $\underline{Z} = (Z_1, \ldots, Z_N)$ be a zero-mean Gaussian vector, $\phi \in C^1(\R^N,\R)$ such that the derivatives are subexponential. Then for $1 \leq l \leq N$, we have
\[
E(Z_l \phi(\underline{Z})) = \sum_{j=1}^N Cov(Z_l,Z_j)E(\pa_j \phi(\underline{Z})).
\]
\end{lem}
Applying Lemma \ref{lI1} iteratively, it is possible to show the following.
\begin{prop} \label{pI2} Let $f \in C^{n+2}(\R)$, such that $f^{(n+2)}$ is subexponential. Let $(G_1,G_2)$ be a Gaussian vector such that $Var(G_2)=1$. We have the following.
\begin{description}
	\item[a)] $n! E(f(G_1)H_n(G_2)) = E(f^{(n)}(G_1))Cov(G_1,G_2)^n,$ 
	\item[b)] $(n-1)! E(f'(G_1)H_{n-1}(G_2)) = E(f^{(n)}(G_1))Cov(G_1,G_2)^{n-1}, $
	\item[c)] $n! E(f''(G_1)H_n(G_2)) = E(f^{(n+2)}(G_1))Cov(G_1,G_2)^n)$. 
\end{description}
 \end{prop}
Let $(X_t)$ be a process such that $1_{[0,t]}\in L_R$ and $X_t = \iin 1_{[0,t]}dX$ for every $t \geq 0$. We recall that under Assumption (D) this is always verified.
\par We denote $\gamma(t)=Var(X_t)$.
\begin{rem} \label{remark}
\begin{description}
	\item[a)] $R(t,\infty) = Var(X_t) -Cov(X_t,X_{\infty}-X_t)$ 
	\item[b)] $t \lra Var(X_t)$ has locally bounded variation if and only if $t \lra Cov(X_t,X_\infty - X_t)$ has locally bounded variation.
\end{description}
\end{rem}
\begin{rem} \label{r43}
It is not easy to find general conditions so that $t \lra \g(t)$ has locally bounded variation even though this condition is often realized. We give for the moment some examples.
\begin{enumerate} 
	\item If $X$ has a convergence measure structure $\g$ has always bounded variation, see \cite{KRT}, Lemma 8.12.
	\item If $X_t = \tilde{X}_{t \wedge T}$, and $(\tilde X_t)$ is a process with weak stationary increments, then $\g(t)=Q(t)$ has always bounded variation, under for instance the assumptions of Proposition \ref{p1.7}.
	\item \[
X_t = \int_0^t G(t-s) dW_s, \ G \in L^2_{loc}(\R).
\]
In this case $\g(t)=\int_0^t G^2(u)du$, which is increasing and therefore locally of bounded variation.
	\item In all explicit examples considered until now, e.g. fractional Brownian motion, bi-fractional Brownian motion, then $\gamma$ has locally bounded variation.
\end{enumerate}
\end{rem}
We can now state the It\^o's formula in the singular case. We recall that from the beginning we suppose $E(X_t)=0$, $\forall t \geq 0$.
\begin{prop} \label{p34} We denote $\g(t)= Var(X_t)$, which is supposed to have locally bounded variation. Let $f:\R \lra \R$ of class $C^\infty_b$. Then $f'(X)1_{[0,t]}$ belongs to $\Dstar$ and
\[
\int_0^t f'(X_s) 1_{[0,t]}(s)\d^* X_s = f(X_t)-f(X_0) -\fr \int_0^t f''(X_s)d\g_s.
\]
\end{prop}
\textbf{Proof} (of Proposition \ref{p34}): We proceed similarly to \cite{CN} and \cite{MV}. We first observe that
\[
g(X_s)1_{[0,t]}(s)
\]
belongs to $\mathcal{L}^2 \subset \mathcal M$, for every bounded function $g$. In particular this is true for $g=f'$. Moreover 
\[
Z_f:= f(X_t)-f(X_0)-\fr \int_0^t f''(X_s)d\g_s
\]
belongs to $\LOm$, since $f$ has linear growth, $f''$ is bounded and $X$ is square integrable. In agreement with (\ref{E434}) for $u_s=f'(X_s)\jeden(s)$ we have to prove that for any $F$ of the type $H_n(\iin \phi dX)$, $\phi \in C^1_0$, $n \geq 0$.
\bea
E(FZ_f) &=& \int_0^t R(ds,\infty)E(D_s F f'(X_s))\nn \\
\label{E440}\\
&-& \int_{\RR} \m(ds_1,ds_2)E((f'(X_{s_1})\jeden(s_1)-f'(X_{s_2})\jeden(s_2))(D_{s_1}F-D_{s_2}F)).\nn
\eea
Without restriction to generality we suppose $\|\phi \|_{\cH}=1$. In this case by Proposition \ref{PPol},  we have 
\[
D_sF=H_{n-1}\left(\iin \phi dX \right) \phi(s).
\]
The right-hand side of (\ref{E440}) becomes
\bea
&&E\left(\int_0^t R(ds,\infty) H_{n-1}\left(\iin \phi dX \right)\phi(s) f'(X_s)\right)\nn \\
\label{E441}\\
&& -\int_{\RR} \m(ds_1,ds_2)(\phi(s_1)-\phi(s_2))E\left( (f'(X_{s_1})\jeden(s_1)-f'(X_{s_2})\jeden(s_2))H_{n-1}\left( \iin \phi dX \right)\right) \nn \\
&& = E\left(H_n \left(\iin \phi dX \right)Z_f \right). \nn
\eea
We recall that by convention $H_{-1}$ is set to zero.   
We denote
\[
p(\s,y)= \frac{1}{\sqrt{2 \pi} \s}\exp\left( - \frac{y^2}{2 \s}\right), \s >0, y \in \R
\]
and we recall that 
\[
\frac{\pa p}{\pa \s}= \fr \frac{\pa^2 p}{\pa y^2}.
\]
Hence for all $n \in \N$ and $t>0$, we evaluate
\[
\frac{d}{dt}E\left( f^{(n)} (X_t)\right) = \frac{d}{dt}E\left( f^{(n)} (\sqrt{\gamma_t} N)\right), 
\]
where $N\sim N(0,1)$. We show that
\be \label{E442}
\frac{d}{dt}E\left( f^{(n)}(X_t)\right)=E\left(f^{(n+2)}(X_t) \right) \frac{d\gamma_t}{2}
\ee
in the sense of measures.
Let $\a \in C_0^{\infty}(\R)$ be a test function. With help of classical Lebesgue-Stieltjes integration theory, we have
\bea
&&\left< \frac{d}{dt}E(f^{(n)}(\sqrt{\gamma_t}  N)), \a \right> = \int_0^\infty \a (t) \left(\frac{d}{dt} \int_{\R} f^{(n)}(y)p(\gamma_t,y)dy\right) dt \nn \\
&&=  \iin \a(t)\left( \int_{\R}dy f^{(n)}(y)\frac{\pa p}{\pa \s}(\gamma_t,y))\right) d\gamma_t 
= \iin \a(t)\left( \int_{\R}dy f^{(n)}(y)\fr \frac{\pa^2}{\pa y^2}p(\gamma_t,y))\right) d\gamma_t \nn \\ 
&&= \iin \a(t)\left( \int_{\R}dy f^{(n+2)}(y) p(\gamma_t,y))\right) \frac{d\gamma_t}{2} = \iin \frac{d\gamma_t}{2}\a(t)E(f^{(n+2)}(\sqrt{\gamma_t} N)) \nn \\
&&=  \iin \frac{d\gamma_t}{2}\a(t)E(f^{(n+2)}(X_t)).  \nn
\eea
which proves (\ref{E442}). We prove now (\ref{E441}). The case of $n=0$ holds if we prove that $E(Z_f)=0$. This expectation gives
\[
E\left(f(X_t)-f(0)-\fr \iin f''(X_s)d\gamma_s \right),
\]
which vanishes applying (\ref{E442}) for $n=0$. It remains to prove (\ref{E441}) for $n \geq 1$.
Proposition \ref{pI2} implies
\begin{description}
	\item[a)] $E\left( H_{n-1}\left( \iin \phi dX\right) f'(X_s)\right)= \frac{1}{(n-1)!}E\left( f^{(n)}(X_s)\right) \left<1_{[0,s]},\phi \right>^{n-1}_{\cH} $
	
	\item[b)] $E\left( H_{n}\left( \iin \phi dX\right) f(X_s) \right)= \frac{1}{n!}E\left( f^{(n)}(X_s)\right) \left<1_{[0,s]},\phi \right>^{n}_{\cH}$
	
	\item[c)] $E\left( H_{n}\left( \iin \phi dX\right) f''(X_s)\right)= \frac{1}{n!}E\left( f^{(n+2)}(X_s)\right) \left<1_{[0,s]},\phi \right>^{n}_{\cH}$
\end{description}
Similarly to \cite{CN}, Lemma 4.3, we evaluate, as a measure,
\be
\frac{d}{dt}E\left( f(X_t)H_n\left( \iin \phi dX\right)\right). \label{E4119}
\ee
In Lemma \ref{lIto1} below we will show that $t \longmapsto \left<\phi,1_{[0,t]} \right>$ has bounded variation. By b), \eqref{E4119} gives
\bea
&& \frac{1}{n!} \frac{d}{dt}\left(E(f^{(n)}(X_t)\left< 1_{[0,t]},\phi\right>_{\cH}^n)\right)\nn  \\
&& = I_{1,n}(dt) + I_{2,n}(dt). \nn
\eea
where
\bea
I_{1,n}(t) &=& \frac{1}{n!} \frac{d}{dt}E\left(f^{(n)}(X_t)\right) \left< 1_{[0,t]},\phi\right>_{\cH}^n  \nn \\  I_{2,n}(t) &=& \frac{1}{(n-1)!} \int_0^t E(f^{(n)}(X_s)) \left<1_{[0,s]},\phi \right>_{\cH}^{n-1}d\left<1_{[0,t]},\phi\right>_{\cH} ds\nn
\eea
Using (\ref{E442}), we have
\[
I_{1,n}(t) = \int_0^t\frac{1}{n!}E(f^{(n+2)}(X_s))\left<1_{[0,s]},\phi\right>_{\cH}^n \frac{d\gamma_s}{2}
\]
which gives
\be
I_{1,n}(t) = \int_0^t E\left(f''(X_s)H_n\left( \iin \phi dX \right)  \right) \frac{d\gamma_s}{2} \label{E4221}
\ee
using c). Concerning the second term, a) implies
\be \label{E4222}
I_{2,n}(t) = E\left( H_{n-1}\left( \iin \phi dX\right)f'(X_s) d\left<1_{[0,s]},\phi \right>_{\cH}(s)\right).
\ee
By Remark \ref{RI1}, \eqref{E4221} and \eqref{E4222} we get
\bea
&&E\left((f(X_t)-f(0)) H_n\left(\iin \phi dX \right) \right) = E \left( f(X_t) H_n \left(\iin \phi dX \right) \right) \nn \\
&& = I_{1,n}(t) + I_{2,n}(t) = \int_0^t E \left(f''(X_s)H_n\left( \iin \phi dX \right) d\g_s  \right) + I_{2,n}(t).\nn
\eea
This implies that
\[
E \left( Z_f H_n \left(\iin \phi dX \right)\right) = I_{2,n}(t).
\]
It remains to prove that the left-hand side of (\ref{E441}) equals $I_{2,n}(t)$. Setting 
\newline $g(s):= E\left(f'(X_s)1_{[0,s]}H_{n-1}\left( \iin \phi dX\right) \right)$ we have to prove that
\bea
&&\int_0^t R(ds,\infty) g(s)\phi(s) - \fr \int_{\RR} \m(ds_1,ds_2)(\phi(s_1)- \phi(s_2))(g(s_1) - g(s_2))\nn\\
\label{E445} \\
&& = \int_0^t g(s) d\left<1_{[0, \cdot]},\phi \right>_{\cH}(s). \nn
\eea
This will the object of the following lemma.
\begin{lem}\label{lIto1} Let $\phi \in C_0^1(\R_+)$, $g: \R_+ \lra \R$ continuous and bounded.
\begin{enumerate}
	\item $t \longmapsto \left<1_{[0,t]},\phi\right>_{\cH}$ has bounded variation,
	\item (\ref{E445}) holds.
\end{enumerate}
\end{lem}
\textbf{Proof}: We denote by $\m_{\phi}$ the antisymmetric measure on $\RR$ defined by
\[
d\m_{\phi}(s_1,s_2)= \int_{\RR} d\bar{\m} (s_1,s_2) \frac{\phi(s_1)-\phi(s_2)}{s_1-s_2};
\]
the right-hand side is well-defined since $\phi \in C_0^1(\R_+)$. It follows
\bea
&&<1_{[0,t]},\phi>_{\cH}\nn \\
&&= \int_0^t R(ds,\infty)\phi(s) - \fr \int_{\RR}(1_{[0,t]}(s_1)-1_{[0,t]}(s_2))(\phi(s_1)-\phi(s_2))d\m(s_1,s_2) \nn \\
&& =  \int_0^t R(ds,\infty) \phi(s) - \fr \m_{\phi}([0,t]\times \R_+) + \fr \m_{\phi} (\R_+ \times [0,t]).
\eea
Therefore 
\be
\left<1_{[0,t]},\phi \right>_{\cH}= \int_0^t R(ds,\infty) \phi(s) + m_\phi([0,t]), \label{E446}
\ee
where $m_\phi([0,t])= \m_{\phi}([0,t] \times \R_+)$. This concludes the proof of 1).
\newline 2) The left-hand side of (\ref{E445}) gives
\[
\iin g(s)\phi(s)R(ds,\infty) + \iin g(s)dm_\phi(s) = \iin g(s) d<1_{[0, \cdot]},\phi>_{\cH}(s) 
\] \qed
\newline At this point (\ref{E440}) is established for every $F= H_n\left( \iin \phi dX \right)$, $\|\phi \|_{\cH}=1. $ If $\|\phi \|=0$ then (\ref{E441}) holds trivially. If $\| \phi \|_{\cH}= \s >0 $ then (\ref{E440}) with $\| \phi \|_{\cH}=1 $ can be extended to this case replacing $X$ with $\s X$. 
 
\section{Wiener and Skorohod integrals} \label{s4.7}

\setcounter{equation}{0}
If the integrand is deterministic, the Wiener integral equals Skorohod integral as Proposition \ref{p454} below shows. 
We list here some properties, whose proof is very close to the one of \cite{KRT}, where we supposed that $X$ has a covariance  measure structure. We suppose Assumptions (A), (B), (C) by default.

\begin{prop}\label{p454} Let $h \in L_R$. Then $h \in Dom(\d)$ and
\[
\iin h \d X = \iin h dX.
\]
\end{prop}
\textbf{Proof: } It follows from Proposition \ref{PR4211} and the definition of Skorohod integral.
\begin{prop} \label{p455}
Let $u \in Cyl(L_R)$. Then $u \in Dom(\d)$ and $\iin u \d X \in L^p(\Om)$, $\forall p\geq 1$.
\end{prop}
\textbf{Proof}: Let $u= G\psi$, $\psi \in L_R$, $G \in Cyl$. Proposition \ref{p454} says that $\psi \in Dom(\d)$. Applying Proposition \ref{p4.5.1} with $F = G$, $u = \psi$, it follows that $\psi G \in Dom(\d)$ and 
\[
\iin u \d X = G \iin \psi \d X - \left<\psi,DG\right>_{\cH}.
\]
Making explicit previous equality when $G=g(Y_1,\ldots, Y_n)$, where $g \in C^\infty_0(\R^n)$, $Y_j=\iin \f_j dX$, $1 \leq j \leq n$, then
\be
\iin u \d X= g(Y_1, \ldots, Y_n) \iin \psi dX - \sum_{j=1}^n \left<\f_j,\psi\right>_{\cH} \pa_j g(Y_1,\ldots,Y_n). \label{E455new}
\ee
The right-hand side belongs obviously to each $L^p$ since $Y_j$ is a Gaussian random variable and $g$, $\pa_j g$ are bounded. The final result for $u \in Cyl(L_R)$ follows by linearity. \qed
\begin{rem} \label{r457} (\ref{E455new}) provides an explicit expression of $\iin u \d X$, if $u \in Cyl(L_R)$.
\end{rem}
Next result concerns the commutation property of the derivative and Skorohod integral. First we observe that $(D_t F) \in Dom(\d)$ if $F \in Cyl$. Moreover, if $u \in Cyl(L_R)$, $(D_s u_t)$ belongs to $\DDd$. Closely to Proposition 7.3 of \cite{KRT} and to \cite{Nualart}, Chapter 1, (1.46) we can prove the following result.
\begin{prop} \label{p458} Let $u \in Cyl(L_R)$. Then
\[
\iin u \d X \in \DD
\]
and for every $t$
\[
D_t\left( \iin u \d X\right)= u_t + \iin (D_t u_s)\d X. 
\]
\end{prop}
We can now evaluate the $\LOm$-norm of the Skorohod integral.
\begin{prop} \label{p459} Let $u\in \DDd$. Then $u \in Dom(\d)$, $\iin u \d X \in \LOm$ and
\bea 
&&E\left( \iin u \d X \right)^2= E(\| u\|^2_{\cH})  \nn \\
&&\label{E459}\ \ \ \ -\fr E \left( \int_{\R_+} d\m(t_1,t_2) \left<D_{\cdot}(u_{t_1}-u_{t_2}),(D_{t_1}-D_{t_2})u_{\cdot}\right>_{\cH} \right)\\
&& \ \ \ \ \ \ \ \ + E\left( \iin R(dt,\infty)\|D_t u_{\cdot} \|^2_{\cH}\right). \nn
\eea
Moreover
\bea 
&&E\left( \iin u \d X \right)^2  \nn \\
&&\label{E460}\\
&& \leq E\left(\| u\|^2_{R} + \int_{\R_+} d|R|(dt,\infty)\| D_t u_{\cdot}^2\|_R + \fr \int_{\RR} d|\m|(s_1,s_2)\|D_{\cdot}u_{s_1}-D_{\cdot}u_{s_2} \|_R^2 \right). \nn
\eea
\begin{rem} \label{r460} We denote $(D^1u)(s,t)= D_t u_s$, $(Du)(s,t)=D_s u_t,
s, t \ge 0. $

\begin{description}
	\item[i)] The right-hand side of (\ref{E459}) can be written as
	\[
	E\left( \| u\|^2_{\cH} + \left<Du,D^1u\right>_{\cH \otimes \cH}\right),
	\] 
	\item[ii)] The right-hand side of (\ref{E460}) can be written as
	\[
	E \left( \| u\|^2_R + \|Du \|^2_{2,R} \right).
	\]
\end{description}
\end{rem}
\end{prop}
\textbf{Proof}(of Proposition \ref{p459}): Let $u \in Cyl(L_R)$. By Proposition \ref{p458}, $\iin u \d X \in \DD$ and we get
\bea
&&E \left( \iin u \d X\right)^2 = E\left(\left<u, D\iin u \d X \right>_{\cH} \right)\nn\\
&&= E \left( \iin u_s D_s \left(\iin u \d X \right)R(ds,\infty)\right) \nn\\
&& - \fr E \left( \int_{\RR} (u_{s_1}- u_{s_2})(D_{s_1}-D_{s_2})\left( \iin u \d X\right)\m(ds_1,ds_2) \right)\nn\\
&&= E_1 -\fr E_2,\nn
\eea 
where
\bea
E_1 &=& E\left(R(dt,\infty)u^2_t + \iin u_t\left(\iin (D_t u_r)\d X_r \right)R(dt,\infty) \right),\nn\\
E_2 &=& E\left(\int_{\RR} \m(dt_1,dt_2)(u_{t_1}-u_{t_2})^2 + \int_{\RR} \m(dt_1,dt_2)(u_{t_1}-u_{t_2})\iin (D_{t_1}-D_{t_2})u_r \d X_r\right).\nn
\eea
Consequently
\bea
&&E_1-\fr E_2 = E\left( \|u \|^2_{\cH}\right) + \iin R(dt,\infty)E\left(u_t \iin D_t u_r \d X_r \right) \nn\\
&& - \fr \int_{\RR} \m(dt_1,dt_2) E \left((u_{t_1}-u_{t_2}) \iin (D_{t_1}-D_{t_2})u_r \d X_r \right). \nn
\eea
Using again the duality relation of Skorohod integral, we obtain
\bea
&&E(\| u\|^2_{\cH}) + E\left( \iin R(dt,\infty) \left<D_{\cdot}u_t,D_t u_{\cdot}\right>_{\cH}\right. \nn\\ 
&& \ \ \ \ \ \ \ \ \ \left. -\fr \int_{\RR} d\m(t_1,t_2) \left<D_{\cdot}(u_{t_1}-u_{t_2}),(D_{t_1}-D_{t_2})u_{\cdot}\right>_{\cH} \right). \nn
\eea
This proves (\ref{E459}) and the result for $u \in Cyl(L_R)$. (\ref{E460}) is then a consequence of Cauchy-Schwarz with respect to the inner product $\left<\cdot,\cdot\right>_{\cH}$ and the fact that the norm $\|\cdot \|_R$ (resp. $\| \cdot\|_{{2,R}}$) dominates $\|\cdot\|_{\cH}$ (resp. $\| \cdot\|_{\cH \otimes \cH}$). The general result for $u \in \DDd$ follows because $Cyl(L_R)$ is dense in $\DDd$. \qed

\section{Link between symmetric and Skorohod integral}\label{s4.8}

\setcounter{equation}{0}
We wish now to establish a relation between the symmetric integral via regularization and Skorohod integral. We recall that, in most of our examples, forward integrals do not exist. If $X = B$ is a fractional Brownian motion with Hurst parameter $H < \fr$, then $X$ is not of finite quadratic variation, therefore $\int_0^{\cdot} X d^-X$ does not exist, see Introduction.
\par We suppose here again the validity of Assumptions (A), (B), (C) and $X_t=X_T$, $t\geq T$. We first need a technical lemma.
\begin{lem} \label{l461} Let $Y \in \DDd$ cadlag. For $\eps >0$, we set
\[
Y^{\eps}_t = \frac{1}{2\eps} \int_{(t-\eps)^+}^{(t+\eps)\wedge T} Y_s ds.
\]
We suppose next the validity of next hypothesis.
\newline \textbf{Hypothesis TR} We set $\n(dt)= |R|(dt,\infty)$
\begin{enumerate}
	\item For every $t \geq 0$, $Y_t\in \Dd$.
	\item $t \longmapsto Y_t$ is continuous and bounded in $\Dd$. In particular $t \longmapsto Y_t$ is continuous and bounded in $L^2$, $t \longmapsto DY_t$ is continuous, bounded in $L^2(\Om;L_R)$.
	\item For $r \in \R$ sufficiently small, $t \lra Y_{t+r}$ belongs to $\Ddd$ and $r \longmapsto Y_{\cdot + r}$ is continuous in $\Ddd$ at $r=0$. In particular $r \longmapsto Y_{\cdot +r}$ is continuous in $L^2(\Om;L_R)$ and $r \longmapsto D_{\cdot}Y_{\cdot+r}$ is continuous in $L^2(\Om;L_{2,R})$ at $r=0$.
\end{enumerate}
Then $Y^{\eps}$ converges to $Y$ in $\Ddd$.
\end{lem}
\textbf{Proof}: We set $a(s,t)= 1_{[(t-\eps)^+,(t+\eps) \wedge T]}(s)$, $\r_t = (t+\eps) \wedge T$. By Hypothesis TR and Corollary \ref{c3.72}, $a(s,\cdot) \in L_R$ for every $s\geq 0$. Hence Assumptions i), ii) of Proposition \ref{TPD101} are verified. By Hypothesis TR 2., Assumption iii) of the same Proposition \ref{TPD101} is also valid. Consequently $Y^{\eps}_t = \frac{1}{2\eps}\iin a(s,t)Y_s d\r_s$ defines a process in $\Ddd$. Moreover
\be
D_r Y_t^{\eps} = \frac{1}{2\eps} \int_{(t-\eps)^+}^{(t+\eps)\wedge T} ds D_r Y_s. \label{ETR300}
\ee 
We need to prove the following
\bea
E \left(\|Y-Y^{\eps} \|^2_R \right)&& \xrightarrow[\eps \to 0]{}0, \label{E462}\\
E \left(\| D(Y-Y^{\eps}) \|^2_{2,R}\right)&& \xrightarrow[\eps \to 0] {} 0. \label{E463}
\eea
1) The left hand-side of \eqref{E462}, using Bochner integration properties and Jensen's inequality, is bounded by
\be
\label{EEE1} \frac{1}{2 \eps}  \int_{(-\eps)^+}^{\eps} dr E \left(\|Y_{s+r}-Y_s \|^2_R \right).
\ee
By Hypothesis TR 3. $\lim_{r \to 0}$ $E\left(\|Y_{\cdot+r}-Y_{\cdot} \|_R \right)^2 = 0$ and \eqref{EEE1} converges to zero. \newline Again by Bochner integration properties and Jensen's inequality the left-hand side of \eqref{E463} is bounded by
\[
\frac{1}{2 \eps} \int_{(-\eps)^+}^{\eps} dr E \left(\| D_{\cdot}Y_{\cdot +r} -DY\|^2_{2,R} \right).
\]
Again this converges to zero since $r \longmapsto D_{\cdot}Y_{\cdot +r}$ is continuous in $L_{2,R}$. \qed
\par Hypothesis TR is quite technical. We provide a very important example, which is constituted by $Y= g(X)$, for suitable real functions $g$. Before treating this we need a preliminary lemma which looks similar, but is significantly different from Corollary \ref{PCYR}.
\begin{lem} \label{LTR20} Let $g \in C^2_b(\R)$. Then $Y_t = g(X_{t})$ belongs to $\Dd$ and
\be
\label{ELTR200} D_r Y_t = g'(X_t) 1_{[0,t]}(r), \ \ t \geq 0.
\ee
Moreover the assumptions of Lemma \ref{LD101} are verified with $\r(t)=t$.
\end{lem}
\textbf{Proof: } i) By Assumption (D) $1_{[0,t]} \in L_R$ because of Corollary \ref{c3.72}. Hence the first part and \eqref{ELTR200} follow by Proposition \ref{PPol}.
\newline ii) We continue verifying the assumption of Lemma \ref{LD101}. $Y$ is continuous in $L^2$ because $Y$ is pathwise continuous and $(Y_t)_{t \leq T}$ is uniformly integrable. Indeed $g$ has linear growth and $X$ is Gaussian, so there is constant $const.$ with
\[
\sup_{t \leq T} E(g(X_t))^4 \leq const. \left(1+ \left(\sup_{t \geq 0}Var X_t\right)^2\right).
\]
$Y$ is bounded in $L^2$ since $X$ by similar arguments as above $(g'(X_t))$ is continuous in $L^2$. Let now $t_2>t_1>0$. It follows that
\bea
\|D_{\cdot}Y_{t_2}-D_{\cdot}Y_{t_1} \|^2_R & \leq & (g'(X_{t_2})^2)\| 1_{]t_1,t_2]}\|^2_R + (g'(X_{t_2})-g'(X_{t_1}))^2\|1_{[0,t_1]} \|^2_R \nn \\
&=& 2 g'(X_{t_2})^2 Var(X_{t_2}-X_{t_1}) + 2\left(g'(X_{t_2})-g'(X_{t_1})\right)^2Var(X_{t_1}). \nn
\eea
Since $X$ and $g'(X)$ are continuous in $L^2$ and $g'$ has linear growth, we obtain that $t \longmapsto D_{\cdot}Y_t$ is continuous in $L^2$. By similar arguments $t \longmapsto \| DY_t\|_R^2$ is also bounded. This concludes the proof of the Lemma \ref{LTR20}. \qed

We go on with another step in the investigation between symmetric and Skorohod integral.
\begin{prop} \label{LTR21} Together with the assumptions mentioned
 at the beginning of Section \ref{s4.8}, we suppose
\be
\int_{\RR} \sup_{r \in{[-\eps_0,\eps_0]}} Var(X_{s_1+r}-X_{s_2+r}) d|\m|(s_1,s_2) <\infty \label{ELTRbis}
\ee
for some $\eps_0 > 0$.
Let $g \in C^2_b(\R)$. Then $Y = g(X)$ verifies Hypothesis TR. 
\end{prop} 
\textbf{Proof }: 
\newline 1) Hypothesis TR 1. was the objective of Lemma \ref{LTR20}.
\newline 2) We have
\[
E(g(X_{t_1})-g(X_{t_2}))^2 \leq \|g' \|_{\infty} Var(X_{t_2}-X_{t_1}).
\]
Since $X$ is continuous in $L^2$, then $g(X)$ is continuous in $L^2$. On the other hand \eqref{ELTR200} holds and $t \longmapsto 1_{[0,t]}$ is continuous from $\R_+$ to $L_R$ because 
\[
\|1_{[0,t]}-1_{[0,s]} \|^2_R = Var\left(X_t -X_s \right).
\]
taking into account Corollary \ref{c3.72}.
This implies that $t \longmapsto DY_{t+r}$ is continuous (and bounded) from $\R_+$ to $\Dd$ for $r$ small enough. Consequently Hypothesis TR 2. is valid.
\newline 3) i) By Proposition \ref{R4550} and Remark \ref{RLYR} we know that $X_{\cdot +r} \in \Ddd$ for $s$ small enough. Proposition \ref{PYR} implies, that $Y = g(X_{\cdot+r})$ belongs to $\Ddd$ and $D_sY_t = g'(X_{t+r})1_{[0,t+r]}(s)$.
\newline ii) To conclude the validity of Hypothesis TR 3. we need to show that
\begin{description}
	\item[a)] $r \longmapsto g(X_{\cdot +r})$ is continuous in $L^2(\Om;L_R)$ in a neighbourhood of 0.
	\item[b)] $r \longmapsto (s,t) \longmapsto g'(X_{t+r})1_{[0,t+r]}(s)$ is continuous in a neighbourhood of zero in $L^2(\Om;L_{2,R})$. 
\end{description}
Concerning a), by definition of $\|\cdot \|_R$,
\[
\|g(X_{\cdot +r})-g(X_{\cdot}) \|^2_R \leq \| g'\|_{\infty}\|X_{\cdot+r}-X_{\cdot} \|^2_R.
\]
Taking the expectation and since
\bea
E \left( \| X_{\cdot +r}-X_{\cdot}\|^2_R\right) &=& \iin R(ds,\infty) E\left(X_{s+r}-X_s \right)^2 \nn \\
&+& \fr \int_{\RR} d|\m|(s_1,s_2) E \left( X_{s_1+r}-X_{s_1}-X_{s_2+r}+X_{s_2}\right)^2 \nn.
\eea
Since $X$ is bounded and continuous in $L^2$, \eqref{ELTRbis} and Lebesgue's dominated convergence theorem imply that
\be \label{G333} \lim_{r \to 0} E \left( \|X_{\cdot+r}-X_{\cdot} \|^2_R\right)=0. 
\ee
Concerning b) we have to estimate,
\[
E\left( \|DY_{\cdot +r} - DY \|^2_{2,R}\right),
\]
where $D_sY_{t+r} = g'(X_{t+r})1_{[0,t+r]}(s)$. Previous expectation is bounded by
\[
2 \left(I_1(r)+I_2(r) \right),
\]
where
\bea
I_1(r)&=& E \left(g(X_{t+r}) \right)^2 \|1_{[0,t+r]}-1_{[0,t]} \|^2_R \nn \\
I_2(r)&=& E \left(\| g(X_{\cdot +r})-g(X_{\cdot}) \|^2_{R}\right) \|1_{[0,t]} \|^2_R. \nn
\eea
We clearly have
\[
I_1(r) \leq \left(g(0) + \|g' \|_{\infty}  E\left( X_{t+r}^2\right)\right) \|Var (X_{\cdot+r})-Var(X_{\cdot})  \|^2_R.
\]
Since $X$ is continuous and bounded in $L^2(\Om)$ and taking into account the definition of $\|\cdot \|_R$, \eqref{ELTRbis} and Lebesgue's dominated convergence theorem imply $\lim_{r \to 0} I_1(r)=0$. On the other hand
\[
\| 1_{[0,t]}\|^2_R = \iin Var(X_t) R(dt,\infty) + \fr\int_{\RR} Var\left(X_{t_2}-X_{t_1}\right) d|\m|(t_1,t_2).
\]
Since
\[
E\left(\|g(X_{\cdot+r})-g(X_{\cdot}) \|^2_R \right) \leq \|g'\|_{\infty} \|X_{\cdot + r} - X_{\cdot}\|_R^2,
\]
\eqref{G333} and \eqref{ELTRbis} imply $\lim_{r \to 0} I_2(r) =0.$ This concludes the proof of Proposition \ref{LTR21}. \qed
\begin{rem} \label{LTR22} If $X_t = \tilde{X}_{t \wedge T}$ and $\tilde{X}$ has stationary increments, then \eqref{ELTRbis} and \eqref{ELTR} are equivalent to
\[
\int_{0+} Q(r)|Q''|(dr) < \infty,
\]
whenever $Q(t) = Var \tilde{X}_t$.
\end{rem}

We are able now to state a theorem linking Skorohod integral and regularization integrals. We will introduce first a definition.
\begin{df} \label{DA1}
Let $Y \in \Ddd$. We say that $DY$ admits a \textbf{symmetric trace} if 
\[
\lim_{\eps \to 0} \int_0^{\t} \left< DY_t,1_{[t -\eps,t+ \eps]} \right>_{\cH} \frac{dt}{\eps}
\]
for every $\t >0$ in probability. We denote by $(Tr^0 DY)(\t)$ the mentioned
 quantity.
\end{df}
\begin{thm} \label{PA2} Let $Y$ be a process with the following assumptions
\begin{enumerate}
	\item Assumption (D).
	\item $Y \in \Ddd$.
	\item Hypothesis TR holds.
	\item $Y$ admits a symmetric trace.
\end{enumerate}
Then
\be
\label{ESymI}
\int_0^t Y d^0 X= \int_0^t Y \d X + (Tr^0 DY)(t).
\ee
\end{thm}
\textbf{Proof: } As in Lemma \ref{l461}, we denote
\[
Y^{\eps}_t = \frac{1}{2\eps} \int_{(t-\eps)^+}^{(t+\eps) \wedge T} ds Y_s.
\]
The $\eps$- approximation of the left-hand symmetric integral in \eqref{ESymI} gives
\[
\frac{1}{2 \eps} \int_0^t Y_s(X_{s+\eps}-X_{(s-\eps)^+})ds.
\]
Using Proposition \ref{p454} and Proposition \ref{p4.5.1} previous expression equals
\bea
&&\frac{1}{2 \eps} \int_0^t ds Y_s \iin \d X_u 1_{[s-\eps,s+\eps]}(u) = \nn \\
&& = \int_0^t ds \iin Y_s 1_{[s-\eps,s+\eps]}(u) - \frac{1}{2\eps} \int_0^t ds \left< DY_s,1_{[s-\eps,s+\eps]}\right>_{\cH}.
\eea
Using Fubini's theorem Proposition \ref{p4.5.2}, the last expression equals
\[
(I_1 + I_2) (\eps)
\]
where
\bea
I_1 (\eps)&=& \frac{1}{2 \eps} \iin \d X_u Y^{\eps}_u,\nn \\
I_2 (\eps)&=& \int_0^t \frac{ds}{2 \eps} \left< DY_s,1_{[s-\eps,s+ \eps]} \right>_{\cH}, \nn 
\eea
with
\[
Y^{\eps}_u = \frac{1}{2\eps} \int_{(u-\eps)^+}^{(u+\eps)\wedge t} Y_r dr.
\]
Lemma \ref{l461} implies that $I_1(\eps) \lra \int_0^t Y_u \d X_u$. The definition of symmetric trace implies that $\lim_{\eps \to 0} I_2(\eps)= (Tr^0 DY)(t)$. \qed

Next application will be an application to the case
 $Y = g(X)$, $g \in C^2_b(\R)$. 
\begin{cor} \label{TA3} We suppose the following
\begin{enumerate}
	\item Assumption (D), (\ref{ELTRbis}) are fulfilled. 
	\item $\forall T >0$
	\be
	\label{A41} \sup_{\eps > 0} \int_{\eps}^T |Var(X_{s+\eps}-X_s)-Var(X_s-X_{s-\eps})|ds < \infty.
	\ee
	\item $\g(t) = Var(X_t)$ has bounded variation. 
\end{enumerate}
Then, for $g \in C^3_{b}(\R)$ we have
	\[
	\int_0^t g(X)d^0X = \int_0^t g(X) \d X + \fr \int_0^t g'(X_s) d \g_s, \ \forall t\in[0,T].
	\]
\end{cor}
Finally we are able to state an Itô formula related to the symmetric (Stratonovich) integral via regularization.
\begin{cor} \label{TA4} Under the assumptions 1), 2), 3) of Corollary \ref{TA3}, if $f \in C^{\infty}_b(\R)$, then
\be
f(X_t) = f(X_0) + \int_0^tf'(X_s)d^0X_s.\label{A42}
\ee
\end{cor}
\begin{rem} \label{RTA4}
\begin{enumerate}
	\item The mentioned Ito's formula has to be considered as a side effect of Theorem \ref{PA2}. We do not aim in providing minimal assumptions. A refinement of this formula could concern the case when the paths of $X$ do not belong to $L_R$ and $X$ belongs only to $(Dom\d)^*$. Various techniques developed by \cite{GNRV, RT, CN} for the specific case of fractional (or bifractional) Brownian motion will probably help.
	\item If $X$ has stationary increments, then (\ref{A41}) holds.
\end{enumerate}
\end{rem}
\textbf{Proof }(of Corollary \ref{TA3}): The result follows as consequence of Theorem \ref{PA2}. Indeed $Y = g(X)$ belongs to $\Ddd$ because of Corollary \ref{PCYR}. Hypothesis TR holds because of Proposition \ref{LTR21}. So Hypotheses 1., 2., 3. of Theorem \ref{PA2} are verified. It remains to check that $Y$ admits a symmetric trace and
\be
Tr^0 Dg(X)_{\t} = \int_0^{\t} g'(X_s) d \g_s, \forall \t>0. \label{TTT10}
\ee 
As we said, Corollary \ref{PCYR} implies that
\[
D_r Y_t = g'(X_t) 1_{[0,t]}(r).
\]
Hence for $\t >0$ the left-hand side of \eqref{TTT10} is the limit when $\eps \to 0$ of
\be
\int_{0}^\tau \frac{1}{2 \eps} \left< DY_t, 1_{[t-\eps,t+\eps]}\right>_{\cH} dt = \frac{1}{2} \int_0^{\t} g'(X_t) \left<1_{[0,t]}, 1_{[t-\eps,t+\eps]} \right>_{\cH} dt. \label{EA42}
\ee
We consider the bounded variation function
\[
F_{\eps}(\t) = \int_0^{\t} \left< 1_{[0,t]},1_{[t -\eps, t+\eps]} \right>_{\cH} \frac{dt}{2 \eps}.
\]
Let $T > 0$. If we prove that
\be
\label{EA420} dF_{\eps}(\t) \Rightarrow \frac{d\g(\tau)}{2}, \ \ \t \in [0,T],
\ee
then \eqref{EA420} converges to $\fr \int_0^{\t}g'(X_t)d\g_t$ a.s. for every $\t \geq 0$ and the theorem would be established. To prove \eqref{EA420} we need to establish the following.
\begin{description}
	\item[i)] The total variation $d|F_{\eps}|(T)$ a.e. is bounded in $\eps > 0$.
	\item[ii)] $F_{\eps}(\t) \lra \frac{\g(\t)}{2}$, $\forall \t > 0$. 
\end{description}
Indeed we have
\bea
F_{\eps}(\t)&=& \int_{0}^{\t} \frac{dt}{2 \eps} \left( Cov(X_{t+\eps},X_t)-Cov(X_{t-\eps},X_t) \right) \nn \\
&=& \frac{1}{2 \eps} \int_{0}^{\t} dt \left( (\g(t+\eps) -\g(t)) - Var(X_{t+\eps}-X_t)\right) \nn \\
&-& \frac{1}{2 \eps} \int_0^{\t} dt
\left( \g(t)-\g(t-\eps) -Var(X_t-X_{t-\eps})\right). \nn
\eea
So
\[
F_{\eps}(\t) = I_{1,\eps}(\t) + I_{2,\eps}(\t),
\]
where
\bea
I_{1,\eps}(\t)&=& \frac{1}{2 \eps} \int_{\t-\eps}^{\t} dt \g(t), \nn \\
I_{2,\eps}(\t)&=& \frac{1}{2 \eps} \int_{\t-\eps}^{\t} dt \left(Var(X_{t+\eps}-X_t)-Var(X_t-X_{t-\eps}) \right)
\eea
i) The total variations of $I_{1,\eps}$ on $[0,T]$ are bounded by the total variation of $\g$. The total variation of $I_{2,\eps}$ are bounded because of \eqref{A41}.
\newline ii) Obviously
\bea
&& \lim_{\eps \to 0} I_{1,\eps}(\t) = \frac{\g(t)}{2} \nn \\
&& \lim_{\eps \to 0} I_{2,\eps}(\t) = 0.
\eea
Finally i) and ii) are proved. \qed

{\bf ACKNOWLEDGEMENTS:} 
Part of the work was done during the stay of the
first and third named authors at the Bielefeld University (SFB 701
and BiBoS). They are grateful to Prof. Michael R\"ockner for the invitation.

\end {document}